\documentclass[11pt]{article}
\usepackage{amssymb}

\usepackage{amsmath}
\usepackage{theorem}
\usepackage{here}
\usepackage[dvipdfmx]{color}

\newtheorem{thm}{Theorem}[section]
\newtheorem{lem}{Lemma}[section]

\newtheorem{prp}{Proposition}[section]
\theorembodyfont{\rmfamily}

\makeatletter

\@addtoreset{equation}{section}
\makeatother

\usepackage{comment} %
\usepackage{bm}
\usepackage[dvipdfmx]{graphicx} %

\def\ta{{\tau}}

\def\Xt{{\widetilde X}}

\def\Zt{{\widetilde Z}}

\def\pd{\partial}

\def\al{{\alpha}}
\def\be{{\beta}}

\def\de{{\delta}}

\def\la{{\lambda}}

\def\th{{\theta}}
\def\fai{{\varphi}}

\def\bla{{\text{\boldmath $\lambda$}}}

\def\bth{{\text{\boldmath $\theta$}}}

\def\bro{{\text{\boldmath $\rho$}}}

\def\lah{{\hat \la}}

\def\lat{{\tilde \la}}

\def\blah{{\widehat \bla}}

\def\blat{{\widetilde \bla}}

\def\De{{\Delta}}

\def\Ga{{\Gamma}}

\def\La{{\Lambda}}
\def\a{{\text{\boldmath $a$}}}

\def\d{{\text{\boldmath $d$}}}
\def\e{{\text{\boldmath $e$}}}

\def\j{{\text{\boldmath $j$}}}

\def\v{{\text{\boldmath $v$}}}

\def\x{{\text{\boldmath $x$}}}
\def\y{{\text{\boldmath $y$}}}

\def\A{{\text{\boldmath $A$}}}

\def\W{{\text{\boldmath $W$}}}
\def\X{{\text{\boldmath $X$}}}
\def\Y{{\text{\boldmath $Y$}}}
\def\Z{{\text{\boldmath $Z$}}}

\def\ph{{\hat p}}

\def\dbt{{\tilde \d}}
\def\dt{{\tilde d}}

\def\pt{{\tilde p}}

\def\diag{{\rm diag\,}}

\def\[{{\text{\boldmath $[$}}}
\def\]{{\text{\boldmath $]$}}}

\def\/{{\Bigr/\!\!}}

\def\1r{{\rm (1)}}
\def\2r{{\rm (2)}}
\def\3r{{\rm (3)}}
\def\4r{{\rm (4)}}
\def\5r{{\rm (5)}}

\def\non{{\nonumber}}
\begin{document}
\title{Bayesian Shrinkage Estimation for Stratified Count Data\footnote{This preprint has not undergone peer review (when applicable) or any post-submission improvements or corrections. 
The Version of Record of this article is published in Japanese Journal of Statistics and Data Science, and is available online at https://doi.org/10.1007/s42081-023-00224-z. }}
\author{
Yasuyuki Hamura\footnote{Graduate School of Economics, Kyoto University, 
Yoshida-Honmachi, Sakyo-ku, Kyoto, 606-8501, JAPAN. 
\newline{
E-Mail: yasu.stat@gmail.com}} \
}
\maketitle
\begin{abstract}
In this paper, we consider simultaneous estimation of Poisson parameters in situations where we can use side information in aggregated data. 
We use standardized squared error and entropy loss functions. 
Bayesian shrinkage estimators are derived based on conjugate priors. 
We compare the risk functions of direct estimators and Bayesian estimators with respect to different priors that are constructed based on different subsets of observations. 
We obtain conditions for domination and also prove minimaxity and admissibility in a simple setting. 

\par\vspace{4mm}
\end{abstract}

\section{Introduction}
\label{sec:introduction}
Simultaneous estimation of multiple parameters have been extensively considered in the literature and, in particular, various improved estimators of normal means have been proposed by many authors including Stein (1956) and James and Stein (1962); see, for example, Fourdrinier, Strawderman and Wells (2018). 
Since the normal distribution is closed under many operations, such results can be extended to more general regression models with unequal variances where the number of observations is greater than the number of parameters so that there is no one-to-one correspondence between observations and parameters. 
Because of this feature, simultaneous estimation of normal location parameters has a lot of application in many areas including that of small area estimation; see, for example, Fay and Herriot (1979) and Morris (1983). 

Since Clevenson and Zidek (1975), improved estimators of Poisson parameters which are analogous to those of normal means have also been proposed (e.g., Ghosh and Parsian 1981; Tsui 1979b; Tsui and Press 1982; Ghosh and Yang 1988; Komaki 2004; Komaki 2006; Komaki 2015; Chang and Shinozaki 2018; Stoltenberg and Hjort 2019; Hamura and Kubokawa 2019). 
Similar results have been obtained for general discrete distributions (e.g., Tsui 1979a; Hwang 1982; Ghosh, Hwang and Tsui 1983; Tsui 1984; Tsui 1986a; Tsui 1986b; Chou 1991; Dey and Chung 1992; Hamura and Kubokawa 2020). 
They mainly considered the case where the dimension of the sample space is equal to the dimension of the parameter space. 
Simultaneous estimation of Poisson parameters is important, for example, in disease mapping or mortality rate estimation and, in practice, count data often have more complicated structures (e.g., Tsutakawa, Shoop and Marienfeld 1985; Clayton and Kaldor 1987; Tsutakawa 1988; Manton et al. 1989; Heisterkamp, Doornbos and Gankema 1993). 
However, there have been few theoretical studies for such situations. 
In this paper, we try to fill in this gap as much as possible. 

Specifically, we note the property of the Poisson distribution that for any $\la _1 , \dots , \la _m > 0$, $m \in \mathbb{N}$, if $Y_i \sim {\rm{Po}} ( \la _i )$, $i = 1, \dots , m$, are independent, then $\sum_{i = 1}^{m} Y_i \sim {\rm{Po}} \big( \sum_{i = 1}^{m} \la _i )$. 
We first consider, for example, the model where we observe independent variables $X_1 \sim {\rm{Po}} ( \la _1 ), \dots , X_m \sim {\rm{Po}} ( \la _m )$ and $Y \sim {\rm{Po}} ( \la _1 + \dots + \la _m )$. 
Such a situation arises if $Y$ corresponds to a county and $X_1 , \dots , X_m$ correspond to the cities in it and they are sampled in preliminary and main surveys, respectively. 
We compare the maximum likelihood estimator and Bayesian shrinkage estimators based on $\X = ( X_1 , \dots , X_m )$ and $Y$, as well as those constructed based on $\X $ only. 
In particular, we prove the minimaxity and admissibility of a proposed estimator. 

Next, we consider two extensions: one to the case of multiple sets of observations $( \X _1 , Y_1 ), \dots , ( \X _m , Y_m )$ and the other to the case of entropy loss. 
We obtain conditions for domination in both cases. 
Furthermore, in the latter case, we also take into account an additional observation at the highest level, namely, $Z \sim {\rm{Po}} ( \la _{\cdot , \cdot } )$, where $\la _{\cdot , \cdot } = \sum_{i = 1}^{m} E[ Y_i ]$, under a balancedness assumption. 
Finally, we obtain results for a general hierarchy under an entropy loss function that is in a sense more ``balanced'' in the context of predictive density estimation.

The remainder of the paper is organized as follows. 
In Section \ref{sec:results}, we study the basic model where $\X $ and $Y$ are observed. 
In Section \ref{sec:extensions}, we consider the two extensions mentioned above. 
In Section \ref{sec:general}, we consider a general hierarchical model. 
Proofs and additional results are in the Appendix.

\section{Basic Results}
\label{sec:results} 
\subsection{The problem}
\label{subsec:setting} 
Let $m \in \mathbb{N}$ and suppose that $X_1 , \dots , X_m$ are independently distributed as ${\rm{Po}} ( \la _i )$, $i = 1, \dots , m$, where $\bla = (\la _1 , \dots , \la _m ) \in (0, \infty )^m$ are unknown parameters of interest. 
Suppose further that $Y \sim {\rm{Po}} \big( \sum_{i = 1}^{m} \la _i \big) $ is an additional independent variable. 
We consider the problem of estimating $\bla $ on the basis of $\X = ( X_1 , \dots , X_m )$ and $Y$ under the standardized squared loss 
\begin{align}
L( \d , \bla ) &= \sum_{i = 1}^{m} {( d_i - \la _i )^2 \over \la _i} \text{,} \label{eq:loss} 
\end{align}
$\d = ( d_1 , \dots , d_m ) \in \mathbb{R} ^m$.

\subsection{Shrinkage estimators} 
\label{subsec:estimators} 
In this section, we consider several estimators of $\bla $. 
Let $\La = \sum_{i = 1}^{m} \la _i$ and $\bth = ( \th _1 , \dots , \th _m ) = \bla / \La $. 
Then the likelihood function is 
\begin{align}
p( \X , Y | \bla ) = {\La ^Y \over Y !} e^{- \La } \prod_{i = 1}^{m} \Big( {{\la _i}^{X_i} \over X_i !} e^{- \la _i} \Big) = {1 \over Y !} \Big( \prod_{i = 1}^{m} {1 \over X_i !} \Big) \La ^{X_{\cdot } + Y} e^{- 2 \La } \prod_{i = 1}^{m} {\th _i}^{X_i} \text{,} \non 
\end{align}
where $X_{\cdot } = \sum_{i = 1}^{m} X_i$. 
Therefore, the maximum likelihood estimator of $\La $, $\bth $, and $\bla $ are $( X_{\cdot } + Y) / 2$, $\X / X_{\cdot }$, and 
\begin{align}
\blah ^{\rm{ML}} ( \X , Y) = {X_{\cdot } + Y \over 2} {\X \over X_{\cdot }} \text{,} \label{eq:ML} 
\end{align}
respectively. 

Let 
\begin{align}
\pi ^{\rm{O}} ( \bla ) d\bla &= 1 d\bla \quad \text{and} \quad \pi ( \bla ) d\bla = \Big\{ 1 / \Big( \sum_{i = 1}^{m} \la _i \Big) ^{m - 1} \Big\} d\bla \text{.} \non 
\end{align}
Using these priors for $\bla $ is equivalent to using 
\begin{align}
\pi ^{\rm{O}} ( \La , \bth ) d( \La , \bth ) &= \La ^{m - 1} d( \La , \bth ) \quad \text{and} \quad \pi ( \La , \bth ) d( \La , \bth ) = 1 d( \La , \bth ) \non 
\end{align}
for $( \La , \bth )$, respectively. 
The generalized Bayes estimator of $\bla $ against the prior $\pi ^{\rm{O}}$ is 
\begin{align}
\blah ^{\rm{O}} ( \X , Y) = {X_{\cdot } + Y + m - 1 \over 2} {\X \over X_{\cdot } + m - 1} \text{,} \non 
\end{align}
whereas that against the prior $\pi $ is 
\begin{align}
\blah ( \X , Y) = {X_{\cdot } + Y \over 2} {\X \over X_{\cdot } + m - 1} \text{.} \non 
\end{align}
The estimators $\blah ^{\rm{O}}$ and $\blah $ shrinks $\blah ^{\rm{ML}}$ and $\blah ^{\rm{O}}$, respectively, toward the origin; that is, 
\begin{align}
\blah ^{\rm{O}} ( \X , Y) &= \Big\{ 1 - {(m - 1) Y \over ( X_{\cdot } + Y) ( X_{\cdot } + m - 1)} \Big\} \blah ^{\rm{ML}} ( \X , Y) \label{eq:est_O_shrink} 
\end{align}
and 
\begin{align}
\blah ( \X , Y) &= \Big( 1 - {m - 1 \over X_{\cdot } + Y + m - 1} \Big) \blah ^{\rm{O}} ( \X , Y) \label{eq:est_shrink} \\
&= \Big( 1 - {m - 1 \over X_{\cdot } + Y + m - 1} \Big) \Big\{ 1 - {(m - 1) Y \over ( X_{\cdot } + Y) ( X_{\cdot } + m - 1)} \Big\} \blah ^{\rm{ML}} ( \X , Y) \text{.} \non 
\end{align}

Similarly, we can derive shrinkage estimators which depend only on $\X $. 
We write 
\begin{align}
\blah ^{\rm{O}} ( \X ) = {X_{\cdot } + m - 1 \over 1} {\X \over X_{\cdot } + m - 1} = \X \non 
\end{align}
and 
\begin{align}
\blah ( \X ) = X_{\cdot } {\X \over X_{\cdot } + m - 1} = \Big( 1 - {m - 1 \over X_{\cdot } + m - 1} \Big) \blah ^{\rm{O}} ( \X ) \text{.} \non 
\end{align}
The estimator $\blah ^{\rm{O}} ( \X )$ is the usual maximum likelihood estimator in the case that we observe $\X $ only. 
The shrinkage estimator $\blah ( \X )$ was proposed by Clevenson and Zidek (1975) for this case. 

In Sections \ref{subsec:dominance}, \ref{subsec:minimaxity}, and \ref{subsec:admissibility}, we investigate risk optimality of the estimators derived in this section. 
Our goal is to show that the estimator $\blah ( \X , Y)$ is admissible and minimax and dominates $\blah ^{\rm{ML}} ( \X , Y)$, $\blah ^{\rm{O}} ( \X , Y)$, and $\blah ( \X )$ if $m \ge 2$.

\subsection{Risk comparisons}
\label{subsec:dominance} 
Here, we compare the risk functions of the estimators we considered in the previous section. 
To this end, we use the following lemma, which is due to Hudson (1978). 

\begin{lem}
\label{lem:hudson} 
Let $\la \in (0, \infty )$ and $X \sim {\rm{Po}} ( \la )$. 
Let $\varphi \colon \mathbb{N} _0 \to [0, \infty )$. 
Then 
\begin{align}
E_{\la } [ \la \fai (X) ] = E_{\la } [ X \fai (X - 1) ] \text{,} \non 
\end{align}
where $X \fai (X - 1) = 0$ when $X = 0$. 
In particular, if $\fai (0) = 0$, then 
\begin{align}
E_{\la } \Big[ {1 \over \la } \fai (X) \Big] = E_{\la } \Big[ {1 \over X + 1} \fai (X + 1) \Big] \text{.} \non 
\end{align}
\end{lem}

We first compare the three estimators $\blah ^{\rm{ML}} ( \X , Y)$, $\blah ^{\rm{O}} ( \X , Y)$, and $\blah ( \X , Y)$. 

\begin{thm}
\label{thm:dominance} 
Suppose that $m \ge 2$. 
\begin{itemize}
\item[{\rm{(i)}}]
The maximum likelihood estimator $\blah ^{\rm{ML}} ( \X , Y)$ is dominated by the shrinkage estimator $\blah ^{\rm{O}} ( \X , Y)$. 
\item[{\rm{(ii)}}]
The shrinkage estimator $\blah ^{\rm{O}} ( \X , Y)$ is dominated by the second shrinkage estimator $\blah ( \X , Y)$. 
\end{itemize}
\end{thm}

Theorem \ref{thm:dominance} corresponds to the result of Clevenson and Zidek (1975) that $\blah ( \X )$ dominates $\blah ^{\rm{O}} ( \X )$ if $m \ge 2$. 
In the present case where we observe $Y$ also, it can be shown that $\blah ( \X , Y)$ dominates $\blah ( \X )$. 

\begin{thm}
\label{thm:dominance_2} 
Suppose that $m \ge 2$. 
Then $\blah ( \X )$ is dominated by $\blah ( \X , Y)$. 
\end{thm}

\subsection{Minimaxity}
\label{subsec:minimaxity} 
Here, we consider whether $\blah ^{\rm{ML}} ( \X , Y)$, $\blah ^{\rm{O}} ( \X , Y)$, $\blah ( \X , Y)$, $\blah ^{\rm{O}} ( \X )$, or $\blah ( \X )$ is minimax. 

\begin{lem}
\label{lem:minimaxity} 
\hfill
\begin{itemize}
\item[{\rm{(i)}}]
Let $\pi ^{( \be )} ( \bla ) d\bla = \be ^m \big( \prod_{i = 1}^{m} e^{- \be \la _i} \big) d\bla $ for $\be > 0$. 
Then for any $\be > 0$, we have $\int_{(0, \infty )^m} \pi ^{( \be )} ( \bla ) d\bla = 1$ and 
\begin{align}
\limsup_{\be \to 0} \int_{(0, \infty )^m} E_{\bla } [ L( \blah ^{( \be )} ( \X , Y), \bla ) ] \pi ^{( \be )} ( \bla ) d\bla \ge m - {1 \over 2} \text{,} \non 
\end{align}
where $\blah ^{( \be )} ( \X , Y)$ is the Bayes estimator of $\bla $ under the prior $\pi ^{( \be )}$. 
\item[{\rm{(ii)}}]
If $m \ge 2$, then 
\begin{align}
\sup_{\bla \in (0, \infty )^m} E_{\bla } [ L( \blah ^{\rm{O}} ( \X , Y), \bla ) ] = m - {1 \over 2} \text{.} \non 
\end{align}
If $m = 1$, then 
\begin{align}
\sup_{\bla \in (0, \infty )^m} E_{\bla } [ L( \blah ^{\rm{O}} ( \X , Y), \bla ) ] > m - {1 \over 2} \text{.} \non 
\end{align}
\end{itemize}
\end{lem}

\begin{thm}
\label{thm:minimaxity} 
\hfill
\begin{itemize}
\item[{\rm{(i)}}]
If $m \ge 2$, the shrinkage estimators $\blah ^{\rm{O}} ( \X , Y)$ and $\blah ( \X , Y)$ are minimax. 
If $m = 1$, they are not minimax and the usual estimator $( X_1 + Y) / 2$ is minimax with constant risk $1 / 2$. 
\item[{\rm{(ii)}}]
The maximum likelihood estimator $\blah ^{\rm{ML}} ( \X , Y)$ is not minimax. 
\item[{\rm{(iii)}}]
The estimators $\blah ^{\rm{O}} ( \X )$ and $\blah ( \X )$ depending only on $\X $ are not minimax. 
\end{itemize}
\end{thm}

\subsection{Admissibility}
\label{subsec:admissibility} 
We can prove the admissibility of the second shrinkage estimator $\blah ( \X , Y)$ by using Blyth's method. 
In particular, we use a sequence of smooth functions considered by Maruyama and Strawderman (2020, 2021b). 

\begin{thm}
\label{thm:admissibility} 
The estimator $\blah ( \X , Y)$ is admissible. 
\end{thm}

\section{Extensions}
\label{sec:extensions} 
\subsection{Multiple sets of observations}
\label{subsec:multiple} 
In this section, we consider the following model where there are multiple sets of observations: 
\begin{align}
&X_{i, j} \sim {\rm{Po}} ( \la _{i, j} ) \text{,} \quad j = 1, \dots , n_i \text{,} \quad i = 1, \dots , m \text{,} \non \\
&Y_i \sim {\rm{Po}} ( \la _{i, \cdot } ) \text{,} \quad i = 1, \dots , m \text{.} \non 
\end{align}
Here, $n_1 , \dots , n_m > 0$ are known constants and $\bla = (( \la _{1, 1} , \dots , \la _{1, n_1} ), \dots , ( \la _{m, 1} , \dots , \la _{m, n_m} )) \in (0, \infty )^{n_1} \times \dots \times (0, \infty )^{n_m}$ are unknown parameters and we write $\la _{i, \cdot } = \sum_{j = 1}^{n_i} \la _{i, j}$ for $i = 1, \dots , m$. 
We want to estimate $\bla $ under the loss function 
\begin{align}
L( \d, \bla ) &= \sum_{i = 1}^{m} \sum_{j = 1}^{n_i} {1 \over \la _{i, j}} ( d_{i, j} - \la _{i, j} )^2 \text{,} \non 
\end{align}
$\d = (( d_{i, j} )_{j = 1, \dots , n_i} )_{i = 1, \dots , m}$, on the basis of the observations $\X = (( X_{i, j} )_{j = 1, \dots , n_i} )_{i = 1, \dots , m}$ and $\Y = ( Y_i )_{i = 1, \dots , m}$. 

The likelihood function is 
\begin{align}
p( \X , \Y | \bla ) &= \Big\{ \prod_{i = 1}^{m} \Big( {1 \over Y_i !} \prod_{j = 1}^{n_i} {1 \over X_{i, j} !} \Big) \Big\} \La ^{X_{\cdot , \cdot } + Y_{\cdot }} e^{- 2 \La } \prod_{i = 1}^{m} \Big( {\th _i}^{X_{i, \cdot } + Y_i} \prod_{j = 1}^{n_i} {\rho _{i, j}}^{X_{i, j}} \Big) \text{,} \non 
\end{align}
where $X_{i, \cdot } = \sum_{j = 1}^{n_i} X_{i, j}$ for $j = 1, \dots , n_i$ and $X_{\cdot , \cdot } = \sum_{i = 1}^{m} X_{i, \cdot }$ and $Y_{\cdot } = \sum_{i = 1}^{m} Y_i$ and where $\La = \sum_{i = 1}^{m} \la _{i, \cdot }$, $\th _i = \la _{i, \cdot } / \La $ for $i = 1, \dots , m$, and $\rho _{i, j} = \la _{i, j} / \la _{i, \cdot }$ for $j = 1, \dots , n_i$ for $i = 1, \dots , m$. 
The maximum likelihood estimator is 
\begin{align}
\blah ^{\rm{ML}} ( \X , \Y ) = \Big( \Big( {X_{i, \cdot } + Y_i \over 2} {X_{i, j} \over X_{i, \cdot }} \Big) _{j = 1, \dots , n_i} \Big) _{i = 1, \dots , m} \text{.} \non 
\end{align}

We consider the priors 
\begin{align}
\pi ^{\rm{O}} ( \bla ) d\bla = 1 d\bla \quad \text{and} \quad \pi ( \La , \bth , \bro ) d( \La , \bth , \bro ) = 1 d( \La , \bth , \bro ) \text{.} \non 
\end{align}
The first prior implies 
\begin{align}
\pi ^{\rm{O}} ( \La , \bth , \bro ) d( \La , \bth , \bro ) = \La ^{n_{\cdot } - 1} \Big( \prod_{i = 1}^{m} {\th _i}^{n_i - 1} \Big) d( \La , \bth , \bro ) \text{,} \non 
\end{align}
where $n_{\cdot } = \sum_{i = 1}^{m} n_i$. 
The generalized Bayes estimators under these priors are 
\begin{align}
&\blah ^{\rm{O}} ( \X , \Y ) = \Big( \Big( {X_{i, \cdot } + Y_i + n_i - 1 \over 2} {X_{i, j} \over X_{i, \cdot } + n_i - 1} \Big) _{j = 1, \dots , n_i} \Big) _{i = 1, \dots , m} \non 
\end{align}
and 
\begin{align}
&\blah ( \X , \Y ) = \Big( \Big( {X_{\cdot , \cdot } + Y_{\cdot } \over 2} {X_{i, \cdot } + Y_i \over X_{\cdot , \cdot } + Y_{\cdot } + m - 1} {X_{i, j} \over X_{i, \cdot } + n_i - 1} \Big) _{j = 1, \dots , n_i} \Big) _{i = 1, \dots , m} \text{,} \non 
\end{align}
respectively.

\begin{thm}
\label{thm:multiple} 
\hfill
\begin{itemize}
\item[{\rm{(i)}}]
Suppose that $\max_{1 \le i \le m} n_i \ge 2$. 
Then $\blah ^{\rm{ML}} ( \X , \Y )$ is dominated by $\blah ^{\rm{O}} ( \X , \Y )$. 
\item[{\rm{(ii)}}]
Suppose that $m \ge 2$ and $\min_{1 \le i \le m} n_i \ge 4$. 
Suppose that 
\begin{align}
{(m - 1)^2 \over 2 m} + {m - 1 \over m^2} + {2 n_i - 3 \over m} \le {( n_i - 1) ( n_i - 3) \over 2} \non 
\end{align}
for all $i = 1, \dots , m$. 
Then $\blah ^{\rm{O}} ( \X , \Y )$ is dominated by $\blah ( \X , \Y )$. 
\end{itemize}
\end{thm}

\subsection{Entropy loss} 
\label{subsec:prediction} 
In this section, we consider Bayesian point estimation under the entropy loss, which is known to be directly related to Bayesian predictive density estimation under the Kullback-Leibler divergence (e.g., Komaki (2006, 2015)). 
In this case, if $n_1 , \dots , n_m$ are balanced, we can obtain an analytical result even when we take into account an additional observation, $Z \sim {\rm{Po}} ( \La )$, as well as $\X $ and $\Y $ of Section \ref{subsec:multiple}. 

The entropy loss we consider here is 
\begin{align}
L^{\rm{KL}} ( \dbt , \bla ) &= \sum_{i = 1}^{m} \sum_{j = 1}^{n_i} \Big( \dt _{i, j} - \la _{i, j} - \la _{i, j} \log {\dt _{i, j} \over \la _{i, j}} \Big) \text{,} \non 
\end{align}
where $\dbt = (( \dt _{i, j} )_{j = 1, \dots , n_i} )_{i = 1, \dots , m} \in (0, \infty )^{n_1} \times \dots \times (0, \infty )^{n_m}$. 
Since the elements of $\blah ^{\rm{ML}} ( \X , \Y )$ and $\blah ^{\rm{O}} ( \X , \Y )$ of Section \ref{subsec:multiple} are zero with positive probability, their risk functions are not well defined. 
Instead, we consider the Bayes estimator against the Jeffreys prior as a benchmark. 
Note that the likelihood function corresponding to $( \X , \Y , Z)$ is 
\begin{align}
p( \X , \Y , Z | \bla ) &= \Big\{ {1 \over Z !} \prod_{i = 1}^{m} \Big( {1 \over Y_i !} \prod_{j = 1}^{n_i} {1 \over X_{i, j} !} \Big) \Big\} \La ^{X_{\cdot , \cdot } + Y_{\cdot } + Z} e^{- 3 \La } \prod_{i = 1}^{m} \Big( {\th _i}^{X_{i, \cdot } + Y_i} \prod_{j = 1}^{n_i} {\rho _{i, j}}^{X_{i, j}} \Big) \text{.} \non 
\end{align}

\begin{lem}
\label{lem:jeff} 
The Jeffreys prior corresponding to the set of observations $( \X , \Y )$ and that corresponding to the set of observations $( \X , \Y , Z)$ are identical and given by 
\begin{align}
\pi ^{\rm{J}} ( \La , \bth , \bro ) d( \La , \bth , \bro ) \propto \La ^{n_{\cdot } / 2 - 1} \prod_{i = 1}^{m} \Big( {\th _i}^{n_i / 2 - 1} \prod_{j = 1}^{n_i} {\rho _{i, j}}^{1 / 2 - 1} \Big) d( \La , \bth , \bro ) \text{.} \non 
\end{align}
\end{lem}

We consider the class of stick-breaking priors given by 
\begin{align}
\pi ^{( \al ; \a )} d( \La , \bth , \bro ) &= \La ^{\al - 1} \Big\{ \prod_{i = 1}^{m} \Big( {\th _i}^{a_i - 1} \prod_{j = 1}^{n_i} {\rho _{i, j}}^{1 / 2 - 1} \Big) \Big\} d( \La , \bth , \bro ) \text{,} \non 
\end{align}
where $\al > 0$ and $\a = ( a_i )_{i = 1, \dots , m} \in (0, \infty )^m$. 
When we observe $Z$, the generalized Bayes estimator under this prior is 
\begin{align}
\blah ^{( \al ; \a )} ( \X , \Y , Z) &= \Big( \Big( {X_{\cdot , \cdot } + Y_{\cdot } + Z + \al \over 3} {X_{i, \cdot } + Y_i + a_i \over X_{\cdot , \cdot } + Y_{\cdot } + a_{\cdot }} {X_{i, j} + 1 / 2 \over X_{i, \cdot } + n_i / 2} \Big) _{j = 1, \dots , n_i} \Big) _{i = 1, \dots , m} \text{,} \non 
\end{align}
where $a_{\cdot } = \sum_{i = 1}^{m} a_i$. 
In particular, that under the Jeffreys prior is 
\begin{align}
\blah ^{\rm{J}} ( \X , \Y , Z) &= \Big( \Big( {X_{\cdot , \cdot } + Y_{\cdot } + Z + n_{\cdot } / 2 \over 3} {X_{i, \cdot } + Y_i + n_i / 2 \over X_{\cdot , \cdot } + Y_{\cdot } + n_{\cdot } / 2} {X_{i, j} + 1 / 2 \over X_{i, \cdot } + n_i / 2} \Big) _{j = 1, \dots , n_i} \Big) _{i = 1, \dots , m} \text{.} \non 
\end{align}
They are 
\begin{align}
\blah ^{( \al ; \a )} ( \X , \Y ) &= \Big( \Big( {X_{\cdot , \cdot } + Y_{\cdot } + \al \over 2} {X_{i, \cdot } + Y_i + a_i \over X_{\cdot , \cdot } + Y_{\cdot } + a_{\cdot }} {X_{i, j} + 1 / 2 \over X_{i, \cdot } + n_i / 2} \Big) _{j = 1, \dots , n_i} \Big) _{i = 1, \dots , m} \non 
\end{align}
and 
\begin{align}
\blah ^{\rm{J}} ( \X , \Y ) &= \Big( \Big( {X_{i, \cdot } + Y_i + n_i / 2 \over 2} {X_{i, j} + 1 / 2 \over X_{i, \cdot } + n_i / 2} \Big) _{j = 1, \dots , n_i} \Big) _{i = 1, \dots , m} \non 
\end{align}
when we do not observe $Z$. 

First, we compare $\blah ^{\rm{J}} ( \X , \Y )$ and $\blah ^{( \al ; \a )} ( \X , \Y )$. 

\begin{thm}
\label{thm:KL_0} 
Suppose that $n_i / 2 \ge a_i \ge 1$ for all $i = 1, \dots , m$ and that $a_{\cdot } \ge \al \ge 1$. 
Suppose that $\al < n_{\cdot } / 2$. 
Then $\blah ^{\rm{J}} ( \X , \Y )$ is dominated by $\blah ^{( \al ; \a )} ( \X , \Y )$. 
\end{thm}

Next, we derive the remarkable result that, under the unbalancedness assumption that $n_1 = \dots = n_m$, we can take into account the additional observation at the highest level, namely, $Z \sim {\rm{Po}} ( \La )$. 

\begin{thm}
\label{thm:KL} 
Assume that $n_1 = \dots = n_m > 2$. 
Suppose that $1 < \al < n_{\cdot } / 2$ and $(m + 1) / m \le a_1 = \dots = a_m < n_1 / 2$. 
Suppose that 
\begin{align}
{2 \over 3} ( n_{\cdot } / 2 - \al ) ( \al - 1) \ge {3 \over 2} (m - 1) ( n_{\cdot } / 2 - a_{\cdot } ) \label{eq:KL_A1} 
\end{align}
and that 
\begin{align}
{( n_{\cdot } / 2 - \al ) ( n_{\cdot } / 2 - 1) \over 3} \ge {(m - 1) ( n_{\cdot } / 2 - a_{\cdot } ) \over 2} \text{.} \label{eq:KL_A2} 
\end{align}
Then $\blah ^{\rm{J}} ( \X , \Y , Z)$ is dominated by $\blah ^{( \al ; \a )} ( \X , \Y , Z)$ under the entropy loss. 
\end{thm}

Assumption (\ref{eq:KL_A1}) is satisfied if $(9 / 4) (m - 1) \le \al - 1 \le a_{\cdot } - 1$. 
Assumption (\ref{eq:KL_A2}) is satisfied if, in addition, $(3 / 2) (m - 1) \le n_{\cdot } / 2 - 1$.

\section{A General Hierarchical Model}
\label{sec:general} 
Here, we consider a general hierarchy with an arbitrary number of levels. 
Theorems \ref{thm:general_1} and \ref{thm:general_2} are not direct extensions of the results of the previous sections. 

Let $D \in \mathbb{N}$ and $n_1 , \dots , n_D \in \mathbb{N}$. 
Suppose that $\bla = ( \cdots ( \la _{i_1 , \dots , i_D} )_{i_D = 1, \dots , n_D} \cdots )_{i_1 = 1, \dots , n_1} \in \prod_{i_1 = 1}^{n_1} \dots \prod_{i_D = 1}^{n_D} (0, \infty ) = (0, \infty )^{n_1 \dotsm n_D}$ are $( n_1 \dotsm n_D )$-dimensional unknown parameters and that we observe $\sum_{d = 0}^{D} \prod_{d' = 1}^{d} n_{d'}$ independent Poisson variables 
\begin{align}
\X = (( \cdots ( X_{( i_1 , \dots , i_d )} )_{i_d = 1, \dots , n_d} \cdots )_{i_1 = 1, \dots , n_1} )_{d = 0, 1, \dots , D} \non 
\end{align}
with means 
\begin{align}
\Big( \Big( \cdots \Big( \sum_{i_{d + 1} = 1}^{n_{d + 1}} \cdots \sum_{i_D = 1}^{n_D} \la _{i_1 , \dots , i_D} \Big) _{i_d = 1, \dots , n_d} \cdots \Big) _{i_1 = 1, \dots , n_1} \Big) _{d = 0, 1, \dots , D} \text{,} \non 
\end{align}
where the first variable is $X$ with mean $\la = \sum_{i_1 = 1}^{n_1} \cdots \sum_{i_D = 1}^{n_D} \la _{i_1 , \dots , i_D}$ and the last variables are $( \cdots ( X_{( i_1 , \dots , i_D )} )_{i_D = 1, \dots , n_D} \cdots )_{i_1 = 1, \dots , n_1}$ with means $\bla $. 

For $d = 1, \dots , D$, let 
\begin{align}
( \cdots ( \la _{( i_1 , \dots , i_d )} )_{i_d = 1, \dots , n_d} \cdots )_{i_1 = 1, \dots , n_1} = \Big( \cdots \Big( \sum_{i_{d + 1} = 1}^{n_{d + 1}} \cdots \sum_{i_D = 1}^{n_D} \la _{i_1 , \dots , i_D} \Big) _{i_d = 1, \dots , n_d} \cdots \Big) _{i_1 = 1, \dots , n_1} \non 
\end{align}
and 
\begin{align}
( \cdots ( \th _{( i_1 , \dots , i_d )} )_{i_d = 1, \dots , n_d} \cdots )_{i_1 = 1, \dots , n_1} = ( \cdots ( \la _{( i_1 , \dots , i_d )} / \la _{( i_1 , \dots , i_{d - 1} )} )_{i_d = 1, \dots , n_d} \cdots )_{i_1 = 1, \dots , n_1} \text{.} \non 
\end{align}
Note that for any $d = 1, \dots , D$, we have $E_{\bla } [ X_{( i_1 , \dots , i_d )} ] = \la _{( i_1 , \dots , i_d )} = \la \prod_{d' = 1}^{d} \th _{( i_1 , \dots , i_{d'} )}$ for all $( i_1 , \dots , i_d ) \in \prod_{d' = 1}^{d} \{ 1, \dots , n_{d'} \} $. 
Then the likelihood function is 
\begin{align}
p( \X | \bla ) &= \prod_{d = 0}^{D} \prod_{i_1 = 1}^{n_1} \cdots \prod_{i_d = 1}^{n_d} {{\la _{( i_1 , \dots , i_d )}}^{X_{( i_1 , \dots , i_d )}} \over X_{( i_1 , \dots , i_d )} !} e^{- \la _{( i_1 , \dots , i_d )}} \non \\
&= \Big( \prod_{d = 0}^{D} \prod_{i_1 = 1}^{n_1} \cdots \prod_{i_d = 1}^{n_d} {1 \over X_{( i_1 , \dots , i_d )} !} \Big) \la ^{X_{(+)}} e^{- (1 + D) \la } \prod_{d = 1}^{D} \prod_{i_1 = 1}^{n_1} \cdots \prod_{i_d = 1}^{n_d} {\th _{( i_1 , \dots , i_d )}}^{X_{( i_1 , \dots , i_d , +)}} \text{,} \non 
\end{align}
where 
\begin{align}
X_{( i_1 , \dots , i_d , +)} = \sum_{d' = d}^{D} \sum_{i_{d + 1} = 1}^{n_{d + 1}} \cdots \sum_{i_{d'} = 1}^{n_{d'}} X_{( i_1 , \dots , i_d , i_{d + 1} , \dots , i_{d'} )} \non 
\end{align}
for $( i_1 , \dots , i_d ) \in \prod_{d' = 1}^{d} \{ 1, \dots , n_{d'} \} $ for $d = 0, 1, \dots , D$. 

We consider the problem of estimating $\bla $ under the entropy loss 
\begin{align}
\widetilde{L} ^{\rm{KL}} ( \dbt , \bla ) &= \sum_{d = 0}^{D} \sum_{i_1 = 1}^{n_1} \cdots \sum_{i_d = 1}^{n_d} \la _{( i_1 , \dots , i_d )} \Big( {\sum_{i_{d + 1} = 1}^{n_{d + 1}} \cdots \sum_{i_D = 1}^{n_D} \dt _{i_1 , \dots , i_D} \over \la _{( i_1 , \dots , i_d )}} - 1 - \log {\sum_{i_{d + 1} = 1}^{n_{d + 1}} \cdots \sum_{i_D = 1}^{n_D} \dt _{i_1 , \dots , i_D} \over \la _{( i_1 , \dots , i_d )}} \Big) \text{,} \label{eq:general_loss} 
\end{align}
$\dbt =  ( \cdots ( \dt _{i_1 , \dots , i_D} )_{i_D = 1, \dots , n_D} \cdots )_{i_1 = 1, \dots , n_1} \in (0, \infty )^{n_1 \dotsm n_D}$, on the basis of $\X $. 
Note that this entropy loss is different from that considered in Section \ref{subsec:prediction} and is an entropy analogue of the standardized squared error loss considered by Stoltenberg and Hjort (2019). 
Moreover, the above entropy loss corresponds to the predictive density estimation for an independent copy of $\X $ while the entropy loss of Section \ref{subsec:prediction} corresponds to the predictive density estimation for an independent copy of $( \cdots ( X_{( i_1 , \dots , i_D )} )_{i_D = 1, \dots , n_D} \cdots )_{i_1 = 1, \dots , n_1}$; see Proposition \ref{prp:connection} in the Appendix. 

As in Section \ref{subsec:prediction}, we consider the Bayes estimator against the Jeffreys prior as a benchmark. 
Let $\bth = (( \cdots ( \th _{( i_1 , \dots , i_d )} )_{i_d = 1, \dots , n_d} \cdots )_{i_1 = 1, \dots , n_1} )_{d = 0, 1, \dots , D}$. 
For $D' = 0, 1, \dots , D$, let $\X ( D' ) = (( \cdots ( X_{( i_1 , \dots , i_d )} )_{i_d = 1, \dots , n_d} \cdots )_{i_1 = 1, \dots , n_1} )_{d = D' , \dots , D}$. 

\begin{lem}
\label{lem:general_Jeff} 
For any $D' = 0, 1, \dots , D$, the Jeffreys prior corresponding to $\X ( D' )$ is given by 
\begin{align}
\pi ^{\rm{J}} ( \la , \bth ) d( \la , \bth ) \propto \la ^{n_1 \dotsm n_D / 2 - 1} \Big( \prod_{d = 1}^{D} \prod_{i_1 = 1}^{n_1} \cdots \prod_{i_d = 1}^{n_d} {\th _{( i_1 , \dots , i_d )}}^{n_{d + 1} \dotsm n_D / 2 - 1} \Big) d( \la , \bth ) \text{.} \non 
\end{align}
\end{lem}

In general, let 
\begin{align}
\pi ^{( \a )} ( \la , \bth ) d( \la , \bth ) \propto \la ^{a_0 - 1} \Big( \prod_{d = 1}^{D} \prod_{i_1 = 1}^{n_1} \cdots \prod_{i_d = 1}^{n_d} {\th _{( i_1 , \dots , i_d )}}^{a_d - 1} \Big) d( \la , \bth ) \non 
\end{align}
for $\a = ( a_d )_{d = 0, 1, \dots , D} \in (0, \infty )^{1 + D}$. 
For $D' = 0, 1, \dots , D$, let $\Xt _{( i_1 , \dots , i_d )} ( D' ) = 1(d \ge D' ) X_{( i_1 , \dots , i_d )}$ for $( i_1 , \dots , i_d ) \in \prod_{d' = 1}^{d} \{ 1, \dots , n_{d'} \} $ for $d = 0, 1, \dots , D$. 
Then, for $D' = 0, 1, \dots , D$, let $\Xt _{( i_1 , \dots , i_d , +)} ( D' ) = \sum_{d' = d}^{D} \sum_{i_{d + 1} = 1}^{n_{d + 1}} \cdots \sum_{i_{d'} = 1}^{n_{d'}} \Xt _{( i_1 , \dots , i_d , i_{d + 1} , \dots , i_{d'} )} ( D' )$ for $( i_1 , \dots , i_d ) \in \prod_{d' = 1}^{d} \{ 1, \dots , n_{d'} \} $ for $d = 0, 1, \dots , D$. 

\begin{lem}
\label{lem:general_estimator} 
Let $\a = ( a_d )_{d = 0, 1, \dots , D} \in (0, \infty )^{1 + D}$. 
Under the entropy loss (\ref{eq:general_loss}), the generalized Bayes estimator of $\bla $ with respect to $\pi ^{( \a )}$ and $\X ( D' )$ is given by 
\begin{align}
\blah ^{( \a )} ( \X ( D' )) &= \Big( \cdots \Big( {\Xt _{(+)} ( D' ) + a_0 \over 1 + D - D'} \prod_{d = 1}^{D} {\Xt _{( i_1 , \dots , i_d , +)} ( D' ) + a_d \over \sum_{{i_d}' = 1}^{n_d} \Xt _{( i_1 , \dots , i_{d - 1} , {i_d}' , +)} ( D' ) + n_d a_d} \Big) _{i_D = 1, \dots , n_D} \cdots \Big) _{i_1 = 1, \dots , n_1} \non 
\end{align}
for any $D' = 0, 1, \dots , D$. 
\end{lem}

Let 
\begin{align}
\a ^{( D_0 )} ( D' ) = \Big( \Big( \Big\{ \prod_{d' = d + 1}^{D} {n_{d'}}^{1( d' \in (- \infty , D' ] \cup [ D_0 + 1, \infty ))} \Big\} {1 \over 2} \Big) _{d = 0, \dots , D_0 - 1} , \Big( \Big( \prod_{d' = d + 1}^{D} n_{d'} \Big) {1 \over 2} \Big) _{d = D_0 , \dots , D} \Big) \non 
\end{align}
for $D' = 0, 1, \dots , D_0$ for $D_0 = 1, \dots , D$. 
Note that $\pi ^{\rm{J}} = \pi ^{( \a ^{( D_0 )} ( D_0 ))}$ for any $D_0 = 1, \dots , D$

We first compare the risk functions of $\blah ^{( \a ^{( D_0 )} ( D_0 ))} ( \X (D))$ and $\blah ^{( \a ^{( D_0 )} ( D_0 ))} ( \X (0))$. 

\begin{thm}
\label{thm:general_1} 
Suppose that $n_D \ge 2$. 
Then for any $D' = 1, \dots , D$, $\blah ^{( \a ^{( D_0 )} ( D_0 ))} ( \X ( D' ))$ is dominated by $\blah ^{( \a ^{( D_0 )} ( D_0 ))} ( \X ( D' - 1))$. 
In particular, the Bayes estimator with respect to the Jeffreys prior and $\X (D)$ is dominated by that with respect to the Jeffreys prior and $\X $. 
\end{thm}

Next, we compare the risk functions of $\blah ^{( \a ^{( D_0 )} ( D_0 ))} ( \X )$ and $\blah ^{( \a ^{( D_0 )} (0))} ( \X )$ for some $D_0 = 1, \dots , D$. 

\begin{thm}
\label{thm:general_2} 
Let $D_0 = 1, \dots , D$ and let 
\begin{align}
a_{D_0}^{( D_0 )} = \Big( \prod_{d' = D_0 + 1}^{D} n_{d'} \Big) {1 \over 2} \text{.} \non 
\end{align}
Suppose that $n_{D'} \ge 2$ for all $D' = 1, \dots , D_0$. 
Suppose that $a_{D_0}^{( D_0 )} \ge 2$ and that 
\begin{align}
{2 + D - D_0 \over D_0 - 1} {a_{D_0}^{( D_0 )} - 1 \over a_{D_0}^{( D_0 )}} \ge n_{D'} {n_{D' - 1} a_{D_0}^{( D_0 )} \over n_{D' - 1} a_{D_0}^{( D_0 )} - 2} \non 
\end{align}
for all $D' = 2, \dots , D_0$. 
Then for any $D' = 1, \dots , D_0$, $\blah ^{( \a ^{( D_0 )} ( D' ))} ( \X )$ is dominated by $\blah ^{( \a ^{( D_0 )} ( D' - 1))} ( \X )$. 
In particular, the Bayes estimator with respect to the Jeffreys prior and $\X $ is dominated by that with respect to $\pi ^{( \a ^{( D_0 )} (0))}$ and $\X $. 
\end{thm}

\section{Appendix}
\subsection{Proofs}
\label{subsec:proofs} 
Here, we prove the results of the main text. 

\bigskip

\noindent
{\bf Proof of Lemma \ref{lem:hudson}.} \ \ We have 
\begin{align}
E_{\la } [ \la \fai (X) ] &= \sum_{x = 0}^{\infty } \fai (x) {\la ^{x + 1} \over x !} e^{- \la } = \sum_{x = 1}^{\infty } x \fai (x - 1) {\la ^x \over x !} e^{- \la } = E_{\la } [ X \fai (X - 1) ] \text{,} \non 
\end{align}
which proves Lemma \ref{lem:hudson}. 
\hfill$\Box$

\bigskip

\noindent
{\bf Proof of Theorem \ref{thm:dominance}.} \ \ For part (i), let $\De _1 ( \bla ) = E_{\bla } [ L( \blah ^{\rm{O}} ( \X , Y), \bla ) ] - E_{\bla } [ L( \blah ^{\rm{ML}} ( \X , Y), \bla ) ]$ and $( \lah _{1}^{\rm{ML}} ( \X , Y), \dots , \lah _{m}^{\rm{ML}} ( \X , Y)) = \blah ^{\rm{ML}} ( \X , Y)$. 
Then, by (\ref{eq:est_O_shrink}), 
\begin{align}
\De _1 ( \bla ) &= E_{\bla } \Big[ \sum_{i = 1}^{m} {1 \over \la _i} \Big[ \Big\{ \lah _{i}^{\rm{ML}} ( \X , Y) - \la _i - {(m - 1) Y \over ( X_{\cdot } + Y) ( X_{\cdot } + m - 1)} \lah _{i}^{\rm{ML}} ( \X , Y) \Big\} ^2 - \{ \lah _{i}^{\rm{ML}} ( \X , Y) - \la _i \} ^2 \Big] \Big] \non \\
&= E_{\bla } \Big[ \sum_{i = 1}^{m} {1 \over \la _i} \Big[ \{ \lah _{i}^{\rm{ML}} ( \X , Y) \} ^2 \Big\{ {(m - 1)^2 Y^2 \over ( X_{\cdot } + Y)^2 ( X_{\cdot } + m - 1)^2} - 2 {(m - 1) Y \over ( X_{\cdot } + Y) ( X_{\cdot } + m - 1)} \Big\} \non \\
&\quad + 2 \la _i \lah _{i}^{\rm{ML}} ( \X , Y) {(m - 1) Y \over ( X_{\cdot } + Y) ( X_{\cdot } + m - 1)} \Big] \Big] \non \\
&= E_{\bla } \Big[ \sum_{i = 1}^{m} \Big[ {X_i + 1 \over ( X_{\cdot } + 1)^2} {( X_{\cdot } + Y + 1)^2 \over 4} \Big\{ {(m - 1)^2 Y^2 \over ( X_{\cdot } + Y + 1)^2 ( X_{\cdot } + m)^2} - 2 {(m - 1) Y \over ( X_{\cdot } + Y + 1) ( X_{\cdot } + m)} \Big\} \non \\
&\quad + 2 {X_{\cdot } + Y \over 2} {X_i \over X_{\cdot }} {(m - 1) Y \over ( X_{\cdot } + Y) ( X_{\cdot } + m - 1)} \Big] \Big] \text{,} \non 
\end{align}
where the last equality follows from Lemma \ref{lem:hudson}. 
Noting that $\blah ^{\rm{ML}} = \bm{0}^{(m)}$ when $\X = \bm{0}^{(m)}$, we have 
\begin{align}
\De _1 ( \bla ) &= E_{\bla } \Big[ {X_{\cdot } + m \over ( X_{\cdot } + 1)^2} {1 \over 4} \Big\{ {(m - 1)^2 Y^2 \over ( X_{\cdot } + m)^2} - 2 {(m - 1) Y ( X_{\cdot } + Y + 1) \over X_{\cdot } + m} \Big\} + 1( X_{\cdot } \ge 1) {(m - 1) Y \over X_{\cdot } + m - 1} \Big] \non \\
&= E_{\bla } \Big[ {(m - 1) / 4 \over ( X_{\cdot } + 1)^2} \Big\{ {(m - 1) Y^2 \over X_{\cdot } + m} - 2 y ( X_{\cdot } + Y + 1) \Big\} + 1( X_{\cdot } \ge 1) {(m - 1) Y \over X_{\cdot } + m - 1} \Big] \non \\
&= E_{\bla } \Big[ (m - 1) Y \Big[ {1 / 4 \over ( X_{\cdot } + 1)^2} \Big\{ {(m - 1) Y \over X_{\cdot } + m} - 2 ( X_{\cdot } + Y + 1) \Big\} + 1( X_{\cdot } \ge 1) {1 \over X_{\cdot } + m - 1} \Big] \Big] \text{.} \non 
\end{align}
Now, by the covariance inequality, 
\begin{align}
{\De _1 ( \bla ) \over (m - 1) \La } &\le E_{\bla } \Big[ {E_{\bla } [ (m - 1) Y | X_{\cdot } ] \over (m - 1) \La } E_{\bla } \Big[ {1 / 4 \over ( X_{\cdot } + 1)^2} \Big\{ {(m - 1) Y \over X_{\cdot } + m} - 2 ( X_{\cdot } + Y + 1) \Big\} + 1( X_{\cdot } \ge 1) {1 \over X_{\cdot } + m - 1} \Big| X_{\cdot } \Big] \Big] \non \\
&= E_{\bla } \Big[ {1 / 4 \over ( X_{\cdot } + 1)^2} \Big\{ {(m - 1) \La \over X_{\cdot } + m} - 2 ( X_{\cdot } + \La + 1) \Big\} + 1( X_{\cdot } \ge 1) {1 \over X_{\cdot } + m - 1} \Big] \non \\
&= E_{\bla } \Big[ - {1 / 2 \over X_{\cdot } + 1} - {1 / 2 \over ( X_{\cdot } + 1)^2} {X_{\cdot } + (m + 1) / 2 \over X_{\cdot } + m} \La + 1( X_{\cdot } \ge 1) {1 \over X_{\cdot } + m - 1} \Big] \text{,} \non 
\end{align}
where $X_{\cdot } \sim {\rm{Po}} ( \La )$ is independent of $Y$. 
Therefore, by Lemma \ref{lem:hudson}, 
\begin{align}
{\De _1 ( \bla ) \over (m - 1) \La } &\le E_{\bla } \Big[ - {1 / 2 \over X_{\cdot } + 1} - {1 / 2 \over ( X_{\cdot } )^2} {X_{\cdot } + (m - 1) / 2 \over X_{\cdot } + m - 1} X_{\cdot } + 1( X_{\cdot } \ge 1) {1 \over X_{\cdot } + m - 1} \Big] = E_{\bla } [ I( X_{\cdot } ) ] \text{,} \non 
\end{align}
where 
\begin{align}
I(x) &= - {1 / 2 \over x + 1} - 1(x \ge 1) {1 / 2 \over x} {x + (m - 1) / 2 \over x + m - 1} + 1(x \ge 1) {1 \over x + m - 1} \non 
\end{align}
for $x \in \mathbb{N} _0$. 
Finally, direct calculation shows that $I(x) < 0$ for all $x \in \mathbb{N} _0$ if $m \ge 2$. 

For part (ii), let $\De _2 ( \bla ) = E_{\bla } [ L( \blah ( \X , Y), \bla ) ] - E_{\bla } [ L( \blah ^{\rm{O}} ( \X , Y), \bla ) ]$ and $( \lah _{1}^{\rm{O}} ( \X , Y), \dots , \lah _{m}^{\rm{O}} ( \X , Y)) = \blah ^{\rm{O}} ( \X , Y)$. 
Then, by (\ref{eq:est_shrink}), 
\begin{align}
\De _2 ( \bla ) &= E_{\bla } \Big[ \sum_{i = 1}^{m} {1 \over \la _i} \Big[ \Big\{ \lah _{i}^{\rm{O}} ( \X , Y) - \la _i - {m - 1 \over X_{\cdot } + Y + m - 1} \lah _{i}^{\rm{O}} ( \X , Y) \Big\} ^2 - \{ \lah _{i}^{\rm{O}} ( \X , Y) - \la _i \} ^2 \Big] \Big] \non \\
&= E_{\bla } \Big[ \sum_{i = 1}^{m} {1 \over \la _i} \Big[ \{ \lah _{i}^{\rm{O}} ( \X , Y) \} ^2 \Big\{ \Big( {m - 1 \over X_{\cdot } + Y + m - 1} \Big) ^2 - 2 {m - 1 \over X_{\cdot } + Y + m - 1} \Big\} \non \\
&\quad + 2 \la _i \lah _{i}^{\rm{O}} ( \X , Y) {m - 1 \over X_{\cdot } + Y + m - 1} \Big] \Big] \text{.} \non 
\end{align}
By Lemma \ref{lem:hudson}, 
\begin{align}
\De _2 ( \bla ) &= E_{\bla } \Big[ \sum_{i = 1}^{m} \Big[ {X_i + 1 \over ( X_{\cdot } + m)^2} {( X_{\cdot } + Y + m)^2 \over 4} \Big\{ {(m - 1)^2 \over ( X_{\cdot } + Y + m)^2} - 2 {m - 1 \over X_{\cdot } + Y + m} \Big\} \non \\
&\quad + 2 {X_{\cdot } + Y + m - 1 \over 2} {X_i \over X_{\cdot } + m - 1} {m - 1 \over X_{\cdot } + Y + m - 1} \Big] \Big] \non \\
&= E_{\bla } \Big[ {1 / 4 \over X_{\cdot } + m} \{ (m - 1)^2 - 2 (m - 1) ( X_{\cdot } + Y + m) \} + {X_{\cdot } \over X_{\cdot } + m - 1} (m - 1) \Big] \non \\
&= (m - 1) E_{\bla } \Big[ {1 / 4 \over X_{\cdot } + m} \{ m - 1 - 2 ( X_{\cdot } + \La + m) \} + {X_{\cdot } \over X_{\cdot } + m - 1} \Big] \non \\
&= (m - 1) E_{\bla } \Big[ {(m - 1) / 4 \over X_{\cdot } + m} - {1 \over 2} - {\La / 2 \over X_{\cdot } + m} + {X_{\cdot } \over X_{\cdot } + m - 1} \Big] \text{.} \non 
\end{align}
Again by Lemma \ref{lem:hudson}, 
\begin{align}
\De _2 ( \bla ) &= (m - 1) E_{\bla } \Big[ {(m - 1) / 4 \over X_{\cdot } + m} - {1 \over 2} - {X_{\cdot } / 2 \over X_{\cdot } + m - 1} + {X_{\cdot } \over X_{\cdot } + m - 1} \Big] \non \\
&= (m - 1) E_{\bla } \Big[ {(m - 1) / 4 \over X_{\cdot } + m} - {1 \over 2} + {X_{\cdot } / 2 \over X_{\cdot } + m - 1} \Big] \text{.} \non 
\end{align}
The right-hand side of the above equality is negative if $m \ge 2$. 
This completes the proof. 
\hfill$\Box$

\bigskip

\noindent
{\bf Proof of Theorem \ref{thm:dominance_2}.} \ \ Let $\De _3 ( \bla ) = E_{\bla } [ L( \blah ( \X , Y), \bla ) ] - E_{\bla } [ L( \blah ( \X ), \bla ) ]$. 
Then 
\begin{align}
\De _3 ( \bla ) &= E_{\bla } \Big[ \sum_{i = 1}^{m} {1 \over \la _i} \Big\{ \Big( {X_{\cdot } + Y \over 2} {X_i \over X_{\cdot } + m - 1} - \la _i \Big) ^2 - \Big( X_{\cdot } {X_i \over X_{\cdot } + m - 1} - \la _i \Big) ^2 \Big\} \Big] \non \\
&= E_{\bla } \Big[ \sum_{i = 1}^{m} {1 \over \la _i} \Big[ {{X_i}^2 \over ( X_{\cdot } + m - 1)^2} \Big\{ {( X_{\cdot } + Y)^2 \over 4} - {X_{\cdot }}^2 \Big\} - 2 \la _i {X_i \over X_{\cdot } + m - 1} \Big( {X_{\cdot } + Y \over 2} - X_{\cdot } \Big) \Big] \Big] \non \\
&= E_{\bla } \Big[ \sum_{i = 1}^{m} \Big[ {X_i + 1 \over ( X_{\cdot } + m)^2} \Big\{ {( X_{\cdot } + Y + 1)^2 \over 4} - ( X_{\cdot } + 1)^2 \Big\} - 2 {X_i \over X_{\cdot } + m - 1} \Big( {X_{\cdot } + Y \over 2} - X_{\cdot } \Big) \Big] \Big] \non \\
&= E_{\bla } \Big[ {1 \over X_{\cdot } + m} \Big\{ {( X_{\cdot } + Y + 1)^2 \over 4} - ( X_{\cdot } + 1)^2 \Big\} - 2 {X_{\cdot } \over X_{\cdot } + m - 1} \Big( {X_{\cdot } + Y \over 2} - X_{\cdot } \Big) \Big] \text{.} \non 
\end{align}
Since $E_{\bla } [ Y ] = \La $ and $E_{\bla } [ Y^2 ] = \La + \La ^2$, 
\begin{align}
\De _3 ( \bla ) &= E \Big[ {1 \over X_{\cdot } + m} \Big\{ {( X_{\cdot } + 1)^2 + 2 ( X_{\cdot } + 1) \La + \La + \La ^2 \over 4} - ( X_{\cdot } + 1)^2 \Big\} - 2 {X_{\cdot } \over X_{\cdot } + m - 1} \Big( {X_{\cdot } + \La \over 2} - X_{\cdot } \Big) \Big] \non \\
&= E \Big[ {1 \over X_{\cdot } + m} \Big\{ - {3 \over 4} ( X_{\cdot } + 1)^2 + {\La \over 2} ( X_{\cdot } + 1) + {(1 + \La ) \La \over 4} \Big\} + {{X_{\cdot }}^2 \over X_{\cdot } + m - 1} - \La {X_{\cdot } \over X_{\cdot } + m - 1} \Big] \text{.} \non 
\end{align}
By Lemma \ref{lem:hudson}, 
\begin{align}
\De _3 ( \bla ) &= E_{\bla } \Big[ - {3 \over 4} {( X_{\cdot } + 1)^2 \over X_{\cdot } + m} + {1 \over 2} {{X_{\cdot }}^2 \over X_{\cdot } + m - 1} + {1 + \La \over 4} {X_{\cdot } \over X_{\cdot } + m - 1} + {{X_{\cdot }}^2 \over X_{\cdot } + m - 1} - {X_{\cdot } ( X_{\cdot } - 1) \over X_{\cdot } + m - 2} \Big] \non \\
&= E_{\bla } \Big[ - {3 \over 4} {( X_{\cdot } + 1)^2 \over X_{\cdot } + m} + {3 \over 2} {{X_{\cdot }}^2 \over X_{\cdot } + m - 1} + {1 \over 4} {X_{\cdot } \over X_{\cdot } + m - 1} + {1 \over 4} {X_{\cdot } ( X_{\cdot } - 1) \over X_{\cdot } + m - 2} - {X_{\cdot } ( X_{\cdot } - 1) \over X_{\cdot } + m - 2} \Big] = E_{\bla } [ J( X_{\cdot } ) ] \text{,} \non 
\end{align}
where 
\begin{align}
J(x) &= {1 \over 4} {x \over x + m - 1} + {3 \over 2} {x^2 \over x + m - 1} - {3 \over 4} {(x + 1)^2 \over x + m} - {3 \over 4} {x (x - 1) \over x + m - 2} \non 
\end{align}
for $x \in \mathbb{N} _0$. 
Note that 
\begin{align}
{x \over x + m - 1} - {(x + 1)^2 \over x + m} = {- x^2 - (2 m - 1) x - (m - 1) \over (x + m - 1) (x + m)} \non 
\end{align}
and 
\begin{align}
{x^2 \over x + m - 1} - {x (x - 1) \over x + m - 2} = {(m - 1) x \over (x + m - 1) (x + m - 2)} \non 
\end{align}
for all $x \in \mathbb{N} _0$. 
Then 
\begin{align}
J(x) &= {1 \over 4} {x \over x + m - 1} - {3 \over 4} {x^2 + (2 m - 1) x + m - 1 \over (x + m - 1) (x + m)} + {3 \over 4} {(m - 1) x \over (x + m - 1) (x + m - 2)} \non 
\end{align}
for all $x \in \mathbb{N} _0$. 
In particular, $J(0) < 0$. 

Now, fix $x \in \mathbb{N}$. 
Suppose first that $m = 2$. 
Then 
\begin{align}
J(x) &= {1 \over 4} {x \over x + 1} - {3 \over 4} {x^2 + 3 x + 1 \over (x + 1) (x + 2)} + {3 \over 4} {x \over (x + 1) x} = - {2 x^2 + 4 x - 3 \over 4 (x + 1) (x + 2)} < 0 \text{.} \non 
\end{align}
Next, suppose that $m \ge 3$. 
Then 
\begin{align}
J(x) &= {1 \over 4} {x \over x + m - 1} {m \over x + m} - {1 \over 2} {x^2 \over (x + m - 1) (x + m)} - {3 \over 4} {(2 m - 1) x + m - 1 \over (x + m - 1) (x + m)} + {3 \over 4} {(m - 1) x \over (x + m - 1) (x + m - 2)} \non \\
&\le {1 \over 4} {x \over x + m - 1} {m \over x + m} - {1 \over 2} {x^2 \over (x + m - 1) (x + m)} - {3 \over 4} {(2 m - 1) x + m - 1 \over (x + m - 1) (x + m)} + {3 \over 4} {(m + 1) x \over (x + m - 1) (x + m)} \non \\
&= {1 \over 4} {x \over x + m - 1} {m \over x + m} - {1 \over 2} {x^2 \over (x + m - 1) (x + m)} - {3 \over 4} {(m - 2) x + m - 1 \over (x + m - 1) (x + m)} \non \\
&\le {1 \over 4} {x \over x + m - 1} {m \over x + m} - {3 \over 4} {(m - 2) x + m - 1 \over (x + m - 1) (x + m)} < 0 \text{.} \non 
\end{align}
This completes the proof. 
\hfill$\Box$

\bigskip

\noindent
{\bf Proof of Lemma \ref{lem:minimaxity}.} \ \ Note that 
\begin{align}
\blah ^{( \be )} ( \X , Y) = {X_{\cdot } + Y + m - 1 \over 2 + \be } {\X \over X_{\cdot } + m - 1} \non 
\end{align}
as in Section \ref{subsec:estimators}. 
Let $\blah ^{(0)} ( \X , Y) = \{ ( X_{\cdot } + Y) / 2 \} \X / X_{\cdot } = \blah ^{\rm{O}} ( \X , Y)$. 
Fix $\be \ge 0$. 
Then 
\begin{align}
&E_{\bla } [ L( \blah ^{( \be )} ( \X , Y), \bla ) ] \non \\
&= E_{\bla } \Big[ \sum_{i = 1}^{m} {1 \over \la _i} \Big\{ \Big( {X_{\cdot } + Y + m - 1 \over 2 + \be } {X_i \over X_{\cdot } + m - 1} \Big) ^2 - 2 \la _i {X_{\cdot } + Y + m - 1 \over 2 + \be } {X_i \over X_{\cdot } + m - 1} + {\la _i}^2 \Big\} \Big] \non \\
&= E_{\bla } \Big[ \sum_{i = 1}^{m} \Big\{ {X_i + 1 \over ( X_{\cdot } + m)^2} {( X_{\cdot } + Y + m)^2 \over (2 + \be )^2} - 2 {X_i \over X_{\cdot } + m - 1} {X_{\cdot } + Y + m - 1 \over 2 + \be } + \la _{i} \Big\} \Big] \non \\
&= E_{\bla } \Big[ {1 \over X_{\cdot } + m} {( X_{\cdot } + Y + m)^2 \over (2 + \be )^2} - 2 {X_{\cdot } \over X_{\cdot } + m - 1} {X_{\cdot } + Y + m - 1 \over 2 + \be } + \La \Big] \text{,} \non 
\end{align}
where the second equality follows from Lemma \ref{lem:hudson}. 
Therefore, 
\begin{align}
&E_{\bla } [ L( \blah ^{( \be )} ( \X , Y), \bla ) ] \non \\
&= E_{\bla } \Big[ {1 \over (2 + \be )^2} \Big( X_{\cdot } + m + 2 Y + {Y^2 \over X_{\cdot } + m} \Big) - 2 {X_{\cdot } \over X_{\cdot } + m - 1} {Y \over 2 + \be } - 2 {X_{\cdot } \over 2 + \be } + \La \Big] \non \\
&= E_{\bla } \Big[ {1 \over (2 + \be )^2} \Big( X_{\cdot } + m + 2 \La + {\La \over X_{\cdot } + m} + {\La ^2 \over X_{\cdot } + m} \Big) - 2 {X_{\cdot } \over X_{\cdot } + m - 1} {\La \over 2 + \be } - 2 {X_{\cdot } \over 2 + \be } + \La \Big] \non \\
&= E_{\bla } \Big[ {1 \over (2 + \be )^2} \Big( X_{\cdot } + m + 2 \La + {\La \over X_{\cdot } + m} + {\La X_{\cdot } \over X_{\cdot } + m - 1} \Big) - 2 {X_{\cdot } \over X_{\cdot } + m - 1} {\La \over 2 + \be } - 2 {X_{\cdot } \over 2 + \be } + \La \Big] \non \\
&= E_{\bla } \Big[ {1 \over (2 + \be )^2} \Big( \La + m + 2 \La + {\La \over X_{\cdot } + m} \Big) - \Big\{ {2 \over 2 + \be } - {1 \over (2 + \be )^2} \Big\} {\La X_{\cdot } \over X_{\cdot } + m - 1} + \Big( 1 - {2 \over 2 + \be } \Big) \La \Big] \text{,} \label{lminimaxityp1} 
\end{align}
where the third equality follows from Lemma \ref{lem:hudson}. 
By Jensen's inequality, 
\begin{align}
&E_{\bla } [ L( \blah ^{( \be )} ( \X , Y), \bla ) ] \non \\
&\ge {1 \over (2 + \be )^2} \Big( \La + m + 2 \La + {\La \over \La + m} \Big) - \Big\{ {2 \over 2 + \be } - {1 \over (2 + \be )^2} \Big\} {\La ^2 \over \La + m - 1} + \Big( 1 - {2 \over 2 + \be } \Big) \La \non \\
&= {1 \over (2 + \be )^2} \Big( 3 \La + m + 1 - {m \over \La + m} \Big) - \Big\{ {2 \over 2 + \be } - {1 \over (2 + \be )^2} \Big\} {\La ^2 \over \La + m - 1} + {\be \over 2 + \be } \La \text{.} \non 
\end{align}
Note that 
\begin{align}
{\La ^2 \over \La + m - 1} &= \La - (m - 1) \Big( 1 - {m - 1 \over \La + m - 1} \Big) \text{.} \non 
\end{align}
Then 
\begin{align}
&E_{\bla } [ L( \blah ^{( \be )} ( \X , Y), \bla ) ] \non \\
&\ge {1 \over (2 + \be )^2} (3 \La + m + 1) - \Big\{ {2 \over 2 + \be } - {1 \over (2 + \be )^2} \Big\} \{ \La - (m - 1) \} + {\be \over 2 + \be } \La \non \\
&\quad - {1 \over (2 + \be )^2} {m \over \La + m} - \Big\{ {2 \over 2 + \be } - {1 \over (2 + \be )^2} \Big\} {(m - 1)^2 \over \La + m - 1} \text{.} \label{lminimaxityp2} 
\end{align}
Now, suppose that $\be > 0$. 
Then, since $\bla \sim \pi ^{( \be )}$ implies $\La \sim {\rm{Ga}} (m, \be )$, the Bayes risk under the proper prior $\pi ^{( \be )}$ is 
\begin{align}
&\int_{(0, \infty )^m} E_{\bla } [ L( \blah ^{( \be )} ( \X , Y), \bla ) ] \pi ^{( \be )} ( \bla ) d\bla \non \\
&\ge {1 \over (2 + \be )^2} \Big( {3 m \over \be } + m + 1 \Big) - \Big\{ {2 \over 2 + \be } - {1 \over (2 + \be )^2} \Big\} \Big\{ {m \over \be } - (m - 1) \Big\} + {m \over 2 + \be } \non \\
&\quad - \int_{(0, \infty )^m} \Big[ {1 \over (2 + \be )^2} {m \over \La + m} + \Big\{ {2 \over 2 + \be } - {1 \over (2 + \be )^2} \Big\} {(m - 1)^2 \over \La + m - 1} \Big] \pi ^{( \be )} ( \bla ) d\bla \text{.} \non 
\end{align}
If $m - 1 \neq 0$, by the dominated convergence theorem, 
\begin{align}
\int_{(0, \infty )^m} {\pi ^{( \be )} ( \bla ) \over \La + m - 1} d\bla &= \int_{0}^{\infty } {\{ \be ^m / \Ga (m) \} \La ^{m - 1} e^{- \be \La } \over \La + m - 1} d\La = \int_{0}^{\infty } {\{ 1 / \Ga (m) \} t^{m - 1} e^{- t} \over t / \be + m - 1} dt \to 0 \non 
\end{align}
as $\be \to 0$. 
Similarly, 
\begin{align}
\int_{(0, \infty )^m} {\pi ^{( \be )} ( \bla ) \over \La + m} d\bla \to 0 \non 
\end{align}
as $\be \to 0$. 
Also, 
\begin{align}
&{1 \over (2 + \be )^2} \Big( {3 m \over \be } + m + 1 \Big) - \Big\{ {2 \over 2 + \be } - {1 \over (2 + \be )^2} \Big\} \Big\{ {m \over \be } - (m - 1) \Big\} + {m \over 2 + \be } \non \\
&= {m + 1 \over (2 + \be )^2} + {3 m / \be \over (2 + \be )^2} - \Big\{ {2 \over 2 + \be } - {1 \over (2 + \be )^2} \Big\} {m \over \be } + (m - 1) \Big\{ {2 \over 2 + \be } - {1 \over (2 + \be )^2} \Big\} + {m \over 2 + \be } \non \\
&= {m \over \be } \Big\{ {4 \over (2 + \be )^2} - {4 + 2 \be \over (2 + \be )^2} \Big\} + {m + 1 \over (2 + \be )^2} + (m - 1) \Big\{ {2 \over 2 + \be } - {1 \over (2 + \be )^2} \Big\} + {m \over 2 + \be } \non \\
&= - {2 m \over (2 + \be )^2} + {m + 1 \over (2 + \be )^2} + (m - 1) \Big\{ {2 \over 2 + \be } - {1 \over (2 + \be )^2} \Big\} + {m \over 2 + \be } \text{,} \non 
\end{align}
which converges to $m - 1 / 2$ as $\be \to 0$. 
Thus, 
\begin{align}
\limsup_{\be \to 0} \int_{(0, \infty )^m} E_{\bla } [ L( \blah ^{( \be )} ( \X , Y), \bla ) ] \pi ^{( \be )} ( \bla ) d\bla &\ge m - {1 \over 2} \text{,} \non 
\end{align}
and this completes the proof of part (i). 

By (\ref{lminimaxityp2}), 
\begin{align}
&E_{\bla } [ L( \blah ^{(0)} ( \X , Y), \bla ) ] \ge {1 \over 4} (3 \La + m + 1) - {3 \over 4} \{ \La - (m - 1) \} - {1 \over 4} {m \over \La + m} - {3 \over 4} {(m - 1)^2 \over \La + m - 1} %
\to m - {1 \over 2} \non 
\end{align}
as $\La \to \infty $. 
Meanwhile, by (\ref{lminimaxityp1}), 
\begin{align}
E_{\bla } [ L( \blah ^{(0)} ( \X , Y), \bla ) ] &= E \Big[ {1 \over 4} \Big( \La + m + 2 \La + {\La \over X_{\cdot } + m} \Big) - {3 \over 4} {\La X_{\cdot } \over X_{\cdot } + m - 1} \Big] \label{lminimaxityp3} \\
&= E \Big[ {1 \over 4} \Big( \La + m + 2 \La + {X_{\cdot } \over X_{\cdot } + m - 1} \Big) - {3 \over 4} {X_{\cdot } ( X_{\cdot } - 1) \over X_{\cdot } + m - 2} \Big] \text{,} \non 
\end{align}
where the second inequality follows from Lemma \ref{lem:hudson}. 
Now, suppose first that $m \ge 3$. 
Then, since 
\begin{align}
{X_{\cdot } ( X_{\cdot } - 1) \over X_{\cdot } + m - 2} &= X_{\cdot } - (m - 1) \Big( 1 - {m - 2 \over X_{\cdot } + m - 2} \Big) \text{,} \non 
\end{align}
we have 
\begin{align}
E_{\bla } [ L( \blah ^{(0)} ( \X , Y), \bla ) ] &= E \Big[ {1 \over 4} \Big( \La + m + 2 \La + {X_{\cdot } \over X_{\cdot } + m - 1} \Big) - {3 \over 4} \{ X_{\cdot } - (m - 1) \} - {3 \over 4} {(m - 1) (m - 2) \over X_{\cdot } + m - 2} \Big] \non \\
&= E \Big[ {1 \over 4} \Big( 3 \La + m + {X_{\cdot } \over X_{\cdot } + m - 1} \Big) - {3 \over 4} \{ \La - (m - 1) \} - {3 \over 4} {(m - 1) (m - 2) \over X_{\cdot } + m - 2} \Big] \text{.} \non 
\end{align}
By Jensen's inequality, 
\begin{align}
E_{\bla } [ L( \blah ^{(0)} ( \X , Y), \bla ) ] &\le {1 \over 4} \Big( 3 \La + m + {\La \over \La + m - 1} \Big) - {3 \over 4} \{ \La - (m - 1) \} - {3 \over 4} {(m - 1) (m - 2) \over \La + m - 2} \non \\
&= {1 \over 4} \Big( m + {\La \over \La + m - 1} \Big) + {3 \over 4} (m - 1) - {3 \over 4} {(m - 1) (m - 2) \over \La + m - 2} \non \\
&= m - {1 \over 2} - {(m - 1) / 4 \over \La + m - 1} - {3 \over 4} {(m - 1) (m - 2) \over \La + m - 2} \text{,} \non 
\end{align}
which is less than $m - 1 / 2$ and converges to $m - 1 / 2$ as $\La \to \infty $. 
Thus, 
\begin{align}
\sup_{\La \in (0, \infty )} E_{\bla } [ L( \blah ^{(0)} ( \X , Y), \bla ) ] &= m - {1 \over 2} \text{.} \label{lminimaxityp4} 
\end{align}
Next, suppose that $m = 2$. 
Then, by (\ref{lminimaxityp3}) and Lemma \ref{lem:hudson}, 
\begin{align}
E_{\bla } [ L( \blah ^{(0)} ( \X , Y), \bla ) ] &= E \Big[ {1 \over 4} \Big( \La + 2 + 2 \La + {X_{\cdot } \over X_{\cdot } + 1} \Big) - {3 \over 4} \La \Big( 1 - {1 \over X_{\cdot } + 1} \Big) \Big] \non \\
&= E \Big[ {1 \over 2} + {1 \over 4} {X_{\cdot } \over X_{\cdot } + 1} + {3 \over 4} {\La \over X_{\cdot } + 1} \Big] \non \\
&= E \Big[ {1 \over 2} + {1 \over 4} {X_{\cdot } \over X_{\cdot } + 1} + {3 \over 4} 1( X_{\cdot } \ge 1) \Big] \non \\
&= E \Big[ {1 \over 2} + {1 \over 4} {X_{\cdot } \over X_{\cdot } + 1} + {3 \over 4} (1 - e^{- \La } ) \Big] \text{.} \non 
\end{align}
By Jensen's inequality, 
\begin{align}
E_{\bla } [ L( \blah ^{(0)} ( \X , Y), \bla ) ] &\le {1 \over 2} + {1 \over 4} {\La \over \La + 1} + {3 \over 4} (1 - e^{- \La } ) \text{.} \non 
\end{align}
The right-hand side of the above inequality is less than $3 / 2$ and converges to $3 / 2$ as $\La \to \infty $. 
Thus, we have (\ref{lminimaxityp4}). 
Finally, suppose that $m = 1$. 
Then, by (\ref{lminimaxityp3}), 
\begin{align}
E_{\bla } [ L( \blah ^{(0)} ( \X , Y), \bla ) ] &= E \Big[ {1 \over 4} \{ \La + 1 + 2 \La + 1( X_{\cdot } \ge 1) \} - {3 \over 4} \La 1( X_{\cdot } \ge 1) \Big] \non \\
&= {3 \La \over 4} + {1 \over 4} + {1 - 3 \La \over 4} (1 - e^{- \La} ) = {1 \over 2} + {3 \La - 1 \over 4} e^{- \La} \text{,} \non 
\end{align}
which is greater than $1 / 2 = m - 1 / 2$ if $\La > 1 / 3$. 
This completes the proof of part (ii). 
\hfill$\Box$

\bigskip

\noindent
{\bf Proof of Theorem \ref{thm:minimaxity}.} \ \ For part (i), suppose first that $m \ge 2$. 
Then, by Lemma \ref{lem:minimaxity}, $\blah ^{\rm{O}} ( \X , Y)$ is minimax. 
Therefore, by Theorem \ref{thm:dominance}, $\blah ( \X , Y)$ is minimax. 
Next, suppose that $m = 1$. 
In this case, it is well known that the estimator $( X_1 + Y) / 2$ is minimax with constant risk $1 / 2 = m - 1 / 2$. 
Therefore, by Lemma \ref{lem:minimaxity}, $\blah ^{\rm{O}} ( \X , Y)$ is not minimax. 
Since $\blah ^{\rm{O}} ( \X , Y)$ and $\blah ( \X , Y)$ are identical in this case, $\blah ( \X , Y)$ is not minimax. 

For part (ii), we have 
\begin{align}
E_{\bla } [ L( \blah ^{\rm{ML}} ( \X , Y), \bla ) ] &= %
E_{\bla } \Big[ \sum_{i = 1}^{m} \Big\{ {1 \over \la _i} {{X_i}^2 \over {X_{\cdot }}^2} {( X_{\cdot } + Y)^2 \over 4} - ( X_{\cdot } + Y) {X_i \over X_{\cdot }} + \la _i \Big\} \Big] \non \\
&= E_{\bla } \Big[ \sum_{i = 1}^{m} \Big\{ {X_i + 1 \over ( X_{\cdot } + 1)^2} {( X_{\cdot } + Y + 1)^2 \over 4} - ( X_{\cdot } + Y) {X_i \over X_{\cdot }} + \la _i \Big\} \Big] \non \\
&= E_{\bla } \Big[ {X_{\cdot } + m \over ( X_{\cdot } + 1)^2} {( X_{\cdot } + Y + 1)^2 \over 4} - ( X_{\cdot } + Y) 1( X_{\cdot } \ge 1) + \La \Big] \text{,} \non 
\end{align}
where the second inequality follows from Lemma \ref{lem:hudson}. 
Then 
\begin{align}
E_{\bla } [ L( \blah ^{\rm{ML}} ( \X , Y), \bla ) ] &= E_{\bla } \Big[ {X_{\cdot } + m \over 4} \Big\{ 1 + {2 Y \over X_{\cdot } + 1} + {Y^2 \over ( X_{\cdot } + 1)^2} \Big\} - ( X_{\cdot } + Y) 1( X_{\cdot } \ge 1) + \La \Big] \non \\
&= E_{\bla } \Big[ {X_{\cdot } + m \over 4} \Big\{ 1 + {2 \La \over X_{\cdot } + 1} + {\La + \La ^2 \over ( X_{\cdot } + 1)^2} \Big\} - ( X_{\cdot } + \La ) 1( X_{\cdot } \ge 1) + \La \Big] \non \\
&= E_{\bla } \Big[ {\La + m \over 4} + {1 \over 4} \Big( 2 \La + {\La + \La ^2 \over X_{\cdot } + 1} \Big) + {m - 1 \over 4} \Big\{ {2 \La \over X_{\cdot } + 1} + {\La + \La ^2 \over ( X_{\cdot } + 1)^2} \Big\} \non \\
&\quad - ( X_{\cdot } + \La ) 1( X_{\cdot } \ge 1) + \La \Big] \non \\
&= E_{\bla } \Big[ {3 \La + m \over 4} + {1 \over 4} {(2 m - 1) \La + \La ^2 \over X_{\cdot } + 1} + {m - 1 \over 4} {\La + \La ^2 \over ( X_{\cdot } + 1)^2} \non \\
&\quad - ( X_{\cdot } + \La ) 1( X_{\cdot } \ge 1) + \La \Big] \text{.} \non 
\end{align}
Now, by Section 4.3 of Johnson, Kemp and Kotz (2005), 
\begin{align}
E_{\bla } \Big[ {1 \over X_{\cdot } + 1} \Big] = {1 - e^{- \La } \over \La } \label{tminimaxityp1} 
\end{align}
and 
\begin{align}
E_{\bla } \Big[ {1 \over ( X_{\cdot } + 1)^2} \Big] \ge E_{\bla } \Big[ {1 \over ( X_{\cdot } + 1) ( X_{\cdot } + 2)} \Big] = {1 - (1 + \La ) e^{- \La } \over \La ^2} \text{.} \non 
\end{align}
Therefore, 
\begin{align}
&E_{\bla } [ L( \blah ^{\rm{ML}} ( \X , Y), \bla ) ] \non \\
&\ge {3 \La + m \over 4} + {(2 m - 1) \La + \La ^2 \over 4} {1 - e^{- \La } \over \La } + {m - 1 \over 4} ( \La + \La ^2 ) {1 - (1 + \La ) e^{- \La } \over \La ^2} - \La - \La (1 - e^{- \La } ) + \La \non \\
&= \underline{R} ( \La ) \text{,} \non 
\end{align}
where 
\begin{align}
&\underline{R} ( \La ) = {2 m - 1 \over 2} + {m - 1 \over 4} {1 \over \La } - {1 \over 4} \Big\{ 4 m - 3 + (m - 4) \La + (m - 1) {1 \over \La } \Big\} e^{- \La } \text{,} \non \\
&{\underline{R}}' ( \La ) = - {m - 1 \over 4} {1 \over \La ^2} - {1 \over 4} \Big\{ m - 4 - (m - 1) {1 \over \La ^2} \Big\} e^{- \La } + {1 \over 4} \Big\{ 4 m - 3 + (m - 4) \La + (m - 1) {1 \over \La } \Big\} e^{- \La } \text{.} \non 
\end{align}
We have that $\underline{R} ( \La ) \to m - 1 / 2$ as $\La \to \infty $ and that ${\underline{R}}' ( \La ) \sim - \{ (m - 1) / 4 \} / \La ^2 < 0$ as $\La \to \infty $ if $m \ge 2$. 
If $m = 1$, then 
\begin{align}
&{\underline{R}}' ( \La ) = {3 \over 4} e^{- \La } + {1 \over 4} \Big\{ 1 - 3 \La \Big\} e^{- \La } = \Big( 1 - {3 \over 4} \La \Big) e^{- \La } < 0 \non 
\end{align}
when $\La > 4 / 3$. 
Thus, $E_{\bla } [ L( \blah ^{\rm{ML}} ( \X , Y), \bla ) ] \ge \underline{R} ( \La ) > m - 1 / 2$ for some large $\La $, which proves part (ii). 

For part (iii), it is clear that $\blah ^{\rm{O}} ( \X ) = \X $ has constant risk $m > m - 1 / 2$ and is not minimax. 
Suppose first that $m = 1$. 
Then $E_{\bla } [ L( \blah ( \X ), \bla ) ] = E_{\bla } [ L( \X , \bla ) ] = 1 > m - 1 / 2$. 
Next, suppose that $m \ge 2$. 
Then 
\begin{align}
&E_{\bla } [ L( \blah ( \X ), \bla ) ] - E_{\bla } [ L( \X , \bla ) ] \non \\
&= E_{\bla } \Big[ \sum_{i = 1}^{m} {1 \over \la _i} \Big\{ \Big( X_i - \la _i - X_i {m - 1 \over X_{\cdot } + m - 1} \Big) - ( X_i - \la _i )^2 \Big\} \Big] \non \\
&= E_{\bla } \Big[ \sum_{i = 1}^{m} {1 \over \la _i} \Big[ {X_i}^2 \Big\{ {(m - 1)^2 \over ( X_{\cdot } + m - 1)^2} - 2 {m - 1 \over X_{\cdot } + m - 1} \Big\} + 2 \la _i X_i {m - 1 \over X_{\cdot } + m - 1} \Big] \Big] \non \\
&= E_{\bla } \Big[ \sum_{i = 1}^{m} \Big[ ( X_i + 1) \Big\{ {(m - 1)^2 \over ( X_{\cdot } + m)^2} - 2 {m - 1 \over X_{\cdot } + m} \Big\} + 2 X_i {m - 1 \over X_{\cdot } + m - 1} \Big] \Big] \non \\
&= E_{\bla } \Big[ {(m - 1)^2 \over X_{\cdot } + m} - 2 (m - 1) + 2 X_{\cdot } {m - 1 \over X_{\cdot } + m - 1} \Big] = E_{\bla } \Big[ {(m - 1)^2 \over X_{\cdot } + m} - 2 {(m - 1)^2 \over X_{\cdot } + m - 1} \Big] \text{,} \non 
\end{align}
where the third equality follows from Lemma \ref{lem:hudson}. 
By Jensen's inequaliy, 
\begin{align}
E_{\bla } \Big[ {1 \over X_{\cdot } + m} \Big] \ge {1 \over \La + m} \text{.} \non 
\end{align}
By (\ref{tminimaxityp1}), 
\begin{align}
E_{\bla } \Big[ {1 \over X_{\cdot } + m - 1} \Big] \le E_{\bla } \Big[ {1 \over X_{\cdot } + 1} \Big] \le {1 \over \La } \text{.} \non 
\end{align}
Therefore, 
\begin{align}
&E_{\bla } [ L( \blah ( \X ), \bla ) ] - m = E_{\bla } [ L( \blah ( \X ), \bla ) ] - E_{\bla } [ L( \X , \bla ) ] \ge {(m - 1)^2 \over \La + m} - 2 {(m - 1)^2 \over \La } \to 0 \non 
\end{align}
as $\La \to \infty $. 
Thus, $E_{\bla } [ L( \blah ( \X ), \bla ) ] > m - 1 / 2$ for sufficiently large $\La $, and the result follows. 
\hfill$\Box$

\bigskip

\noindent
{\bf Proof of Theorem \ref{thm:admissibility}.} \ \ Let 
\begin{align}
h_k ( \La ) &= 1 - {\log (1 + \La ) \over \log (1 + k + \La )} \non 
\end{align}
for $k \in \mathbb{N}$. 
Then, by Lemma C.2 of Maruyama and Strawderman (2021b), $h_k ( \La ) \in (0, 1)$, $h_k ( \La ) \le h_{k + 1} ( \La )$, and 
\begin{align}
{h_k}' ( \La ) &= - {k / (1 + \La ) + h_k ( \La ) \over (1 + k + \La ) \log (1 + k + \La )} \label{tadmissibilityp1} 
\end{align}
for any $k \in \mathbb{N}$ and 
\begin{align}
\lim_{k \to \infty } h_k ( \La ) = 1 \text{.} \non 
\end{align}
Also, for any $k \in \mathbb{N}$, 
\begin{align}
\int_{\Xi } \{ h_k ( \La ) \} ^2 \pi ( \La , \bth ) d( \La , \bth ) < \infty \non 
\end{align}
where $\Xi = (0, \infty ) \times \big\{ ( {\la _1}' , \dots , {\la _m}' ) / \sum_{i = 1}^{m} {\la _i}' \big| ( {\la _1}' , \dots , {\la _m}' ) \in \mathbb{R} ^m \big\} $, and the proper Bayes estimator under the prior $\{ h_k ( \La ) \} ^2 \pi ( \La , \bth ) d( \La , \bth )$ is 
\begin{align}
\blah ^{(k)} ( \X , Y) &= ( \lah _{1}^{(k)} ( \X , Y), \dots , \lah _{m}^{(k)} ( \X , Y)) = {\X \over X_{\cdot } + m - 1} \frac{ \int_{0}^{\infty } \{ h_k ( \La ) \} ^2 \La ^{X_{\cdot } + Y} e^{- 2 \La } d\La }{ \int_{0}^{\infty } \{ h_k ( \La ) \} ^2 \La ^{X_{\cdot } + Y - 1} e^{- 2 \La } d\La } \text{.} \non 
\end{align}
To establish the admissibility of $\blah ( \X , Y)$, it is sufficient to show that 
\begin{align}
\lim_{k \to \infty } \De _k = 0 \text{,} \non 
\end{align}
where 
\begin{align}
0 \le \De _k = \int_{\Xi } E_{\La \bth } [ L( \blah ( \X , Y), \La \bth ) - L( \blah ^{(k)} ( \X , Y), \La \bth ) ] \{ h_k ( \La ) \} ^2 \pi ( \La , \bth ) d( \La , \bth ) \non 
\end{align}
for $k \in \mathbb{N}$. 

Let $( \lah _1 ( \X , Y), \dots , \lah _m ( \X , Y)) = \blah ( \X , Y)$. 
As in the proof of Lemma 4 of Maruyama and Strawderman (2021a), for any $k \in \mathbb{N}$, 
\begin{align}
\De _k &= \int_{\Xi } E_{\La \bth } \Big[ \sum_{i = 1}^{m} {1 \over \La \th _i} [ \{ \lah _{i} ( \X , Y) - \La \th _i \} ^2 - \{ \lah _{i}^{(k)} ( \X , Y) - \La \th _i \} ^2 ] \Big] d{\pi _k} ( \La , \bth ) \non \\
&= \int_{\Xi } \sum_{\x \in {\mathbb{N} _0}^m} \sum_{y = 0}^{\infty } \sum_{i = 1}^{m} \Big[ {\{ \lah _{i} ( \x , y) \} ^2 - \{ \lah _{i}^{(k)} ( \x , y) \} ^2 \over \La \th _i} - 2 \{ \lah _{i} ( \x , y) - \lah _{i}^{(k)} ( \x , y) \} \Big] p( \x ,y | \La , \bth ) d{\pi _k} ( \La , \bth ) \non \\
&= \sum_{i = 1}^{m} \sum_{\substack{\x \in {\mathbb{N} _0}^m \\ \e _i \cdot \x \ge 1}} \sum_{y = 0}^{\infty } \Big( [ \{ \lah _{i} ( \x , y) \} ^2 - \{ \lah _{i}^{(k)} ( \x , y) \} ^2 ] \int_{\Xi } {1 \over \La \th _i} p( \x , y | \La \bth ) d{\pi _k} ( \La , \bth ) \non \\
&\quad - 2 \{ \lah _{i} ( \x , y) - \lah _{i}^{(k)} ( \x , y) \} \int_{\Xi } p( \x , y | \La \bth ) d{\pi _k} ( \La , \bth ) \Big) \non \\
&= \sum_{i = 1}^{m} \sum_{\substack{\x \in {\mathbb{N} _0}^m \\ \e _i \cdot \x \ge 1}} \sum_{y = 0}^{\infty } \int_{\Xi } {1 \over \La \th _i} p( \x , y | \La \bth ) d{\pi _k} ( \La , \bth ) [ \{ \lah _{i} ( \x , y) \} ^2 - \{ \lah _{i}^{(k)} ( \x , y) \} ^2 \non \\
&\quad - 2 \lah _{i}^{(k)} ( \x , y) \{ \lah _{i} ( \x , y) - \lah _{i}^{(k)} ( \x , y) \} ] \non \\
&= \sum_{i = 1}^{m} \sum_{\substack{\x \in {\mathbb{N} _0}^m \\ \e _i \cdot \x \ge 1}} \sum_{y = 0}^{\infty } \int_{\Xi } {1 \over \La \th _i} p( \x , y | \La \bth ) d{\pi _k} ( \La , \bth ) \{ \lah _{i} ( \x , y) - \lah _{i}^{(k)} ( \x , y) \} ^2 \text{,} \non 
\end{align}
where $d{\pi _k} ( \La , \bth ) = \{ h_k ( \La ) \} ^2 \pi ( \La , \bth ) d( \La , \bth )$ for $k \in \mathbb{N}$ and where $\e _i$ denotes the $i$th row of the $m \times m$ identity matrix. 
For any $k \in \mathbb{N}$, $i = 1, \dots , m$, $\x = ( x_1 , \dots , x_m ) \in {\mathbb{N} _0}^m$, and $y \in \mathbb{N} _0$, if $x_i \ge 1$, we have 
\begin{align}
&\int_{\Xi } {1 \over \La \th _i} p( \x , y | \La \bth ) d{\pi _k} ( \La , \bth ) \{ \lah _{i} ( \x , y) - \lah _{i}^{(k)} ( \x , y) \} ^2 \non \\
&= \Big( \prod_{j = 1}^{m} {1 \over x_j !} \Big) {1 \over y !} {\prod_{j = 1}^{m} \Ga ( x_j + 1) \over x_i \Ga ( x_{\cdot } + m - 1)} \int_{0}^{\infty } \{ h_k ( \La ) \} ^2 \La ^{x_{\cdot } + y - 1} e^{- 2 \La } d\La \non \\
&\quad \times \Big( {x_i \over x_{\cdot } + m - 1} \Big) ^2 \Big[ {x_{\cdot } + y \over 2} - {\int_{0}^{\infty } \{ h_k ( \La ) \} ^2 \La ^{x_{\cdot } + y} e^{- 2 \La } d\La \over \int_{0}^{\infty } \{ h_k ( \La ) \} ^2 \La ^{x_{\cdot } + y - 1} e^{- 2 \La } d\La } \Big] ^2 \non \\
&= {1 \over y !} {x_i \over \Ga ( x_{\cdot } + m - 1)} {1 \over ( x_{\cdot } + m - 1)^2} \non \\
&\quad \times \int_{0}^{\infty } \{ h_k ( \La ) \} ^2 \La ^{x_{\cdot } + y - 1} e^{- 2 \La } d\La \Big[ {x_{\cdot } + y \over 2} - {\int_{0}^{\infty } \{ h_k ( \La ) \} ^2 \La ^{x_{\cdot } + y} e^{- 2 \La } d\La \over \int_{0}^{\infty } \{ h_k ( \La ) \} ^2 \La ^{x_{\cdot } + y - 1} e^{- 2 \La } d\La } \Big] ^2 \non \\
&= {1 \over y !} {x_i \over \Ga ( x_{\cdot } + m - 1)} {1 \over ( x_{\cdot } + m - 1)^2} {\big\{ \int_{0}^{\infty } h_k ( \La ) {h_k}' ( \La ) \La ^{x_{\cdot } + y} e^{- 2 \La } d\La \big\} ^2 \over \int_{0}^{\infty } \{ h_k ( \La ) \} ^2 \La ^{x_{\cdot } + y - 1} e^{- 2 \La } d\La } \text{,} \non 
\end{align}
where $x_{\cdot } = \sum_{i = 1}^{m} x_i$ and where the last equality follows since 
\begin{align}
&\int_{0}^{\infty } \{ h_k ( \La ) \} ^2 \La ^{x_{\cdot } + y} e^{- 2 \La } d\La \non \\
&= \Big[ - {1 \over 2} \{ h_k ( \La ) \} ^2 \La ^{x_{\cdot } + y} e^{- 2 \La } \Big] _{0}^{\infty } + {1 \over 2} \int_{0}^{\infty } [2 h_k ( \La ) {h_k}' ( \La ) \La ^{x_{\cdot } + y} + \{ h_k ( \La ) \} ^2 ( x_{\cdot } + y) \La ^{x_{\cdot } + y - 1} ] e^{- 2 \La } d\La \non \\
&= {x_{\cdot } + y \over 2} \int_{0}^{\infty } \{ h_k ( \La ) \} ^2 \La ^{x_{\cdot } + y - 1} e^{- 2 \La } d\La + \int_{0}^{\infty } h_k ( \La ) {h_k}' ( \La ) \La ^{x_{\cdot } + y} e^{- 2 \La } d\La \text{.} \non 
\end{align}
Therefore, 
\begin{align}
\De _k &= \sum_{i = 1}^{m} \sum_{\substack{\x \in {\mathbb{N} _0}^m \\ \e _i \cdot \x \ge 1}} \sum_{y = 0}^{\infty } {1 \over y !} {x_i \over \Ga ( x_{\cdot } + m - 1)} {1 \over ( x_{\cdot } + m - 1)^2} {\big\{ \int_{0}^{\infty } h_k ( \La ) {h_k}' ( \La ) \La ^{x_{\cdot } + y} e^{- 2 \La } d\La \big\} ^2 \over \int_{0}^{\infty } \{ h_k ( \La ) \} ^2 \La ^{x_{\cdot } + y - 1} e^{- 2 \La } d\La } \non \\
&\le \sum_{\x \in {\mathbb{N} _0}^m} \sum_{y = 0}^{\infty } {1 \over y !} {1( x_{\cdot } \ge 1) x_{\cdot } \over \Ga ( x_{\cdot } + m - 1)} {1 \over ( x_{\cdot } + m - 1)^2} {\big\{ \int_{0}^{\infty } h_k ( \La ) {h_k}' ( \La ) \La ^{x_{\cdot } + y} e^{- 2 \La } d\La \big\} ^2 \over \int_{0}^{\infty } \{ h_k ( \La ) \} ^2 \La ^{x_{\cdot } + y - 1} e^{- 2 \La } d\La } \non \\
&= \sum_{y = 0}^{\infty } \sum_{x_{\cdot } = 1}^{\infty } {1 \over y !} {x_{\cdot } \over \Ga ( x_{\cdot } + m - 1)} {1 \over ( x_{\cdot } + m - 1)^2} {\big\{ \int_{0}^{\infty } h_k ( \La ) {h_k}' ( \La ) \La ^{x_{\cdot } + y} e^{- 2 \La } d\La \big\} ^2 \over \int_{0}^{\infty } \{ h_k ( \La ) \} ^2 \La ^{x_{\cdot } + y - 1} e^{- 2 \La } d\La } \sum_{\substack{( {x_1}' , \dots , {x_m}' ) \in \{ 0, \dots , x_{\cdot } \} ^m \\ \sum_{i = 1}^{m} {x_i}' = x_{\cdot }}} 1 \non \\
&\le \sum_{y = 0}^{\infty } \sum_{x_{\cdot } = 1}^{\infty } {1 \over y !} {x_{\cdot } \over \Ga ( x_{\cdot } + m - 1)} {1 \over ( x_{\cdot } + m - 1)^2} {\big\{ \int_{0}^{\infty } h_k ( \La ) {h_k}' ( \La ) \La ^{x_{\cdot } + y} e^{- 2 \La } d\La \big\} ^2 \over \int_{0}^{\infty } \{ h_k ( \La ) \} ^2 \La ^{x_{\cdot } + y - 1} e^{- 2 \La } d\La } ( x_{\cdot } + 1)^{m - 1} \text{.} \non 
\end{align}
Note that for any $y \in \mathbb{N} _0$ and $x_{\cdot } \in \mathbb{N}$, 
\begin{align}
{( x_{\cdot } + 1)^{m - 1} \over \Ga ( x_{\cdot } + m - 1)} \le {x_{\cdot } + 1 \over x_{\cdot } !} \text{.} \non 
\end{align}
Then 
\begin{align}
\De _k &\le \sum_{y = 0}^{\infty } \sum_{x_{\cdot } = 1}^{\infty } {1 \over y !} {1 \over x_{\cdot } !} {x_{\cdot } ( x_{\cdot } + 1) \over ( x_{\cdot } + m - 1)^2} {\big\{ \int_{0}^{\infty } h_k ( \La ) {h_k}' ( \La ) \La ^{x_{\cdot } + y} e^{- 2 \La } d\La \big\} ^2 \over \int_{0}^{\infty } \{ h_k ( \La ) \} ^2 \La ^{x_{\cdot } + y - 1} e^{- 2 \La } d\La } \non \\
&\le \sum_{\substack{( x_{\cdot } , y) \in {\mathbb{N} _0}^2 \\ x_{\cdot } + y \ge 1}} {1 \over y !} {1 \over x_{\cdot } !} {x_{\cdot } ( x_{\cdot } + 1) \over ( x_{\cdot } + m - 1)^2} {\big\{ \int_{0}^{\infty } h_k ( \La ) {h_k}' ( \La ) \La ^{x_{\cdot } + y} e^{- 2 \La } d\La \big\} ^2 \over \int_{0}^{\infty } \{ h_k ( \La ) \} ^2 \La ^{x_{\cdot } + y - 1} e^{- 2 \La } d\La } \non \\
&= \sum_{w = 1}^{\infty } {\big\{ \int_{0}^{\infty } h_k ( \La ) {h_k}' ( \La ) \La ^w e^{- 2 \La } d\La \big\} ^2 \over \int_{0}^{\infty } \{ h_k ( \La ) \} ^2 \La ^{w - 1} e^{- 2 \La } d\La } \sum_{x = 0}^{w} {1 \over (w - x) !} {1 \over x !} {x (x + 1) \over (x + m - 1)^2} \non \\
&\le \sum_{w = 1}^{\infty } {\big\{ \int_{0}^{\infty } h_k ( \La ) {h_k}' ( \La ) \La ^w e^{- 2 \La } d\La \big\} ^2 \over \int_{0}^{\infty } \{ h_k ( \La ) \} ^2 \La ^{w - 1} e^{- 2 \La } d\La } {2^w \over w !} \sum_{x = 0}^{w} \binom{w}{x} \Big( {1 \over 2} \Big) ^w \{ 1 + 1(m = 1) \} \non \\
&\le 2 \sum_{w = 1}^{\infty } {2^w \over \Ga (w + 1)} {\big\{ \int_{0}^{\infty } h_k ( \La ) {h_k}' ( \La ) \La ^w e^{- 2 \La } d\La \big\} ^2 \over \int_{0}^{\infty } \{ h_k ( \La ) \} ^2 \La ^{w - 1} e^{- 2 \La } d\La } \text{.} \label{tadmissibilityp1.5} 
\end{align}
Here, by (\ref{tadmissibilityp1}) and by the dominated convergence theorem, 
\begin{align}
\lim_{k \to \infty } {\big\{ \int_{0}^{\infty } h_k ( \La ) {h_k}' ( \La ) \La ^w e^{- 2 \La } d\La \big\} ^2 \over \int_{0}^{\infty } \{ h_k ( \La ) \} ^2 \La ^{w - 1} e^{- 2 \La } d\La } = 0 \label{tadmissibilityp2} 
\end{align}
for any $w \in \mathbb{N}$. 

Fix $k \in \mathbb{N} \cap [2, \infty )$ and $w \in \mathbb{N}$. 
By the Cauchy-Schwarz inequality, 
\begin{align}
{\big\{ \int_{0}^{\infty } h_k ( \La ) {h_k}' ( \La ) \La ^w e^{- 2 \La } d\La \big\} ^2 \over \int_{0}^{\infty } \{ h_k ( \La ) \} ^2 \La ^{w - 1} e^{- 2 \La } d\La } &\le \int_{0}^{\infty } \{ {h_k}' ( \La ) \} ^2 \La ^{w + 1} e^{- 2 \La } d\La = {1 \over 2^{w + 2}} \int_{0}^{\infty } \{ {h_k}' (u / 2) \} ^2 u^{w + 1} e^{- u} du \text{.} \non 
\end{align}
By (\ref{tadmissibilityp1}), 
\begin{align}
| {h_k}' ( \La )| &\le {|k / (1 + \La )| + | h_k ( \La )| \over (1 + k + \La ) \log (1 + k + \La )} \le {2 \over (1 + \La ) \log (3 + \La )} \text{.} \non 
\end{align}
Therefore, 
\begin{align}
{\big\{ \int_{0}^{\infty } h_k ( \La ) {h_k}' ( \La ) \La ^w e^{- 2 \La } d\La \big\} ^2 \over \int_{0}^{\infty } \{ h_k ( \La ) \} ^2 \La ^{w - 1} e^{- 2 \La } d\La } &\le {1 \over 2^w} \int_{0}^{\infty } {u^{w + 1} e^{- u} \over (1 + u / 2)^2 \{ \log (3 + u / 2) \} ^2} du \le {4 \over 2^w} \int_{0}^{\infty } {u^{w - 1} e^{- u} \over \{ \log (3 + u / 2) \} ^2} du \text{.} \non 
\end{align}
Since 
\begin{align}
&\int_{0}^{\infty } {u^{w - 1} e^{- u} \over \{ \log (3 + u / 2) \} ^2} du \non \\
&= \Big[ {1 \over w} {{u^w} e^{- u} \over \{ \log (3 + u / 2) \} ^2} \Big] _{0}^{\infty } - {1 \over w} \int_{0}^{\infty } {u^w} \Big[ {- e^{- u} \over \{ \log (3 + u / 2) \} ^2} + {- 2 e^{- u} \over \{ \log (3 + u / 2) \} ^3} {1 / 2 \over 3 + u / 2} \Big] du \non \\
&= {1 \over w} \int_{0}^{\infty } {u^w} \Big[ {e^{- u} \over \{ \log (3 + u / 2) \} ^2} + {e^{- u} \over \{ \log (3 + u / 2) \} ^3} {1 \over 3 + u / 2} \Big] du \le {4 \over 3} {1 \over w} \int_{0}^{\infty } {{u^w} e^{- u} \over \{ \log (3 + u / 2) \} ^2} du \text{,} \non 
\end{align}
it follows that 
\begin{align}
{\big\{ \int_{0}^{\infty } h_k ( \La ) {h_k}' ( \La ) \La ^w e^{- 2 \La } d\La \big\} ^2 \over \int_{0}^{\infty } \{ h_k ( \La ) \} ^2 \La ^{w - 1} e^{- 2 \La } d\La } &\le {4 \over 2^w} {4 \over 3} {1 \over w} \int_{0}^{\infty } {{u^w} e^{- u} \over \{ \log (3 + u / 2) \} ^2} du \text{.} \non 
\end{align}
Furthermore, letting $M = e + 1$, we have 
\begin{align}
\int_{0}^{\infty } {{u^w} e^{- u} \over \{ \log (3 + u / 2) \} ^2} du &= \int_{0}^{w / M} {u^w e^{- u} \over \{ \log (3 + u / 2) \} ^2} du + \int_{w / M}^{\infty } {{u^w} e^{- u} \over \{ \log (3 + u / 2) \} ^2} du \non \\
&\le \int_{0}^{w / M} (w / M)^w du + \int_{w / M}^{\infty } {{u^w} e^{- u} \over [ \log \{ 3 + w / (2 M) \} ]^2} du \non \\
&\le {w^{w + 1} \over M^{w + 1}} + {\Ga (w + 1) \over [ \log \{ 3 + w / (2 M) \} ]^2} \text{.} \non 
\end{align}
Also, 
\begin{align}
\Ga (w + 1) &\ge (2 \pi )^{1 / 2} {(w + 1)^{w + 1} \over (w + 1)^{1 / 2}} {1 \over e^{w + 1}} \text{.} \non 
\end{align}
Thus, 
\begin{align}
&{2^w \over \Ga (w + 1)} {\big\{ \int_{0}^{\infty } h_k ( \La ) {h_k}' ( \La ) \La ^w e^{- 2 \La } d\La \big\} ^2 \over \int_{0}^{\infty } \{ h_k ( \La ) \} ^2 \La ^{w - 1} e^{- 2 \La } d\La } \non \\
&\le {16 / 3 \over \Ga (w + 1)} {1 \over w} \int_{0}^{\infty } {{u^w} e^{- u} \over \{ \log (3 + u / 2) \} ^2} du \non \\
&\le {16 \over 3} {1 \over w} \Big( {1 \over (2 \pi )^{1 / 2}} {(w + 1)^{1 / 2} \over (w + 1)^{w + 1}} e^{w + 1} {w^{w + 1} \over M^{w + 1}} + {1 \over [ \log \{ 3 + w / (2 M) \} ]^2} \Big) \non \\
&\le {16 \over 3} {1 \over w} \Big( {(w + 1)^{1 / 2} \over (2 \pi )^{1 / 2}} {e^{w + 1} \over M^{w + 1}} + {1 \over [ \log \{ 3 + w / (2 M) \} ]^2} \Big) \text{,} \label{tadmissibilityp3} 
\end{align}
and we note that 
\begin{align}
\sum_{w = 1}^{\infty } {16 \over 3} {1 \over w} \Big( {(w + 1)^{1 / 2} \over (2 \pi )^{1 / 2}} {e^{w + 1} \over M^{w + 1}} + {1 \over [ \log \{ 3 + w / (2 M) \} ]^2} \Big) < \infty \text{.} \label{tadmissibilityp4} 
\end{align}

Finally, by the dominated convergence theorem, we conclude from (\ref{tadmissibilityp1.5}), (\ref{tadmissibilityp2}), (\ref{tadmissibilityp3}), and (\ref{tadmissibilityp4}) that $\lim_{k \to \infty } \De _k = 0$. 
This completes the proof. 
\hfill$\Box$

\bigskip

\noindent
{\bf Proof of Theorem \ref{thm:multiple}.} \ \ For part (i), let $\De _1 ( \bla ) = E_{\bla } [ L( \blah ^{\rm{O}} ( \X , \Y ), \bla ) ] - E_{\bla } [ L( \blah ^{\rm{ML}} ( \X , \Y ), \bla ) ]$. 
Let $(( \lah _{i, j}^{\rm{O}} ( \X , \Y ))_{j = 1, \dots , n_i} )_{i = 1, \dots , m} = \blah ^{\rm{O}} ( \X , \Y )$ and $(( \lah _{i, j}^{\rm{ML}} ( \X , \Y ))_{j = 1, \dots , n_i} )_{i = 1, \dots , m} = \blah ^{\rm{ML}} ( \X , \Y )$. 
Then, since 
\begin{align}
\lah _{i, j}^{\rm{O}} ( \X , \Y ) = \Big\{ 1 - {( n_i - 1) Y_i \over ( X_{i, \cdot } + Y_i ) ( X_{i, \cdot } + n_i - 1)} \Big\} \lah _{i, j}^{\rm{ML}} ( \X , \Y ) \non 
\end{align}
for all $j = 1, \dots , n_i$ for all $i = 1, \dots , m$, 
\begin{align}
\De _1 ( \bla ) &= E_{\bla } \Big[ \sum_{i = 1}^{m} \sum_{j = 1}^{n_i} {1 \over \la _{i, j}} \Big[ \{ \lah _{i, j}^{\rm{ML}} ( \X , \Y ) \} ^2 \Big\{ {( n_i - 1)^2 {Y_i}^2 \over ( X_{i, \cdot } + Y_i )^2 ( X_{i, \cdot } + n_i - 1)^2} - 2 {( n_i - 1) Y_i \over ( X_{i, \cdot } + Y_i ) ( X_{i, \cdot } + n_i - 1)} \Big\} \non \\
&\quad + 2 \la _{i, j} \lah _{i, j}^{\rm{ML}} ( \X , \Y ) {( n_i - 1) Y_i \over ( X_{i, \cdot } + Y_i ) ( X_{i, \cdot } + n_i - 1)} \Big] \Big] \non \\
&= E_{\bla } \Big[ \sum_{i = 1}^{m} \sum_{j = 1}^{n_i} \Big[ {X_{i, j} + 1 \over ( X_{i, \cdot } + 1)^2} {( X_{i, \cdot } + Y_i + 1)^2 \over 4} \Big\{ {( n_i - 1)^2 {Y_i}^2 \over ( X_{i, \cdot } + Y_i + 1)^2 ( X_{i, \cdot } + n_i )^2} - 2 {( n_i - 1) Y_i \over ( X_{i, \cdot } + Y_i + 1) ( X_{i, \cdot } + n_i )} \Big\} \non \\
&\quad + 2 {X_{i, j} \over X_{i, \cdot }} {X_{i, \cdot } + Y_i \over 2} {( n_i - 1) Y_i \over ( X_{i, \cdot } + Y_i ) ( X_{i, \cdot } + n_i - 1)} \Big] \Big] \non \\
&= E_{\bla } \Big[ \sum_{i = 1}^{m} \Big[ {X_{i, \cdot } + n_i \over ( X_{i, \cdot } + 1)^2} {( X_{i, \cdot } + Y_i + 1)^2 \over 4} \Big\{ {( n_i - 1)^2 {Y_i}^2 \over ( X_{i, \cdot } + Y_i + 1)^2 ( X_{i, \cdot } + n_i )^2} - 2 {( n_i - 1) Y_i \over ( X_{i, \cdot } + Y_i + 1) ( X_{i, \cdot } + n_i )} \Big\} \non \\
&\quad + 2 {X_{i, \cdot } \over X_{i, \cdot }} {X_{i, \cdot } + Y_i \over 2} {( n_i - 1) Y_i \over ( X_{i, \cdot } + Y_i ) ( X_{i, \cdot } + n_i - 1)} \Big] \Big] \non \\
&= E_{\bla } \Big[ \sum_{i = 1}^{m} ( n_i - 1) Y_i \Big[ {1 \over ( X_{i, \cdot } + 1)^2} {1 \over 4} \Big\{ {( n_i - 1) Y_i \over X_{i, \cdot } + n_i} - 2 ( X_{i, \cdot } + Y_i + 1) \Big\} + {X_{i, \cdot } \over X_{i, \cdot }} {1 \over X_{i, \cdot } + n_i - 1} \Big] \Big] \text{,} \non 
\end{align}
where the second equality follows from Lemma \ref{lem:hudson}. 
By the covariance inequality, 
\begin{align}
\De _1 ( \bla ) %
&\le \sum_{i = 1}^{m} E_{\bla } [ ( n_i - 1) Y_i ] E_{\bla } \Big[ {1 \over ( X_{i, \cdot } + 1)^2} {1 \over 4} \Big\{ {( n_i - 1) Y_i \over X_{i, \cdot } + n_i} - 2 ( X_{i, \cdot } + Y_i + 1) \Big\} + {X_{i, \cdot } \over X_{i, \cdot }} {1 \over X_{i, \cdot } + n_i - 1} \Big] \text{.} \non 
\end{align}
Now, fix $i = 1, \dots , m$ such that $n_i \ge 2$. 
Then 
\begin{align}
&E_{\bla } \Big[ {1 \over ( X_{i, \cdot } + 1)^2} {1 \over 4} \Big\{ {( n_i - 1) Y_i \over X_{i, \cdot } + n_i} - 2 ( X_{i, \cdot } + Y_i + 1) \Big\} + {X_{i, \cdot } \over X_{i, \cdot }} {1 \over X_{i, \cdot } + n_i - 1} \Big] \non \\
&= E_{\bla } \Big[ {1 \over ( X_{i, \cdot } + 1)^2} {1 \over 4} \Big\{ {( n_i - 1) \la _{i, \cdot } \over X_{i, \cdot } + n_i} - 2 ( X_{i, \cdot } + \la _{i, \cdot } + 1) \Big\} + {X_{i, \cdot } \over X_{i, \cdot }} {1 \over X_{i, \cdot } + n_i - 1} \Big] \non \\
&= E_{\bla } \Big[ {1 \over 4} {1 \over {X_{i, \cdot }}^2} {( n_i - 1) X_{i, \cdot } \over X_{i, \cdot } + n_i - 1} - {1 \over 2} {1 \over X_{i, \cdot } + 1} - {1 \over 2} {X_{i, \cdot } \over {X_{i, \cdot }}^2} + {X_{i, \cdot } \over X_{i, \cdot }} {1 \over X_{i, \cdot } + n_i - 1} \Big] \non \\
&\le E_{\bla } \Big[ {1 \over 4} 1( X_{i, \cdot } \ge 1) \Big( {1 \over X_{i, \cdot }} {n_i - 1 \over X_{i, \cdot } + n_i - 1} - 2 {1 \over X_{i, \cdot } + 1} - 2 {1 \over X_{i, \cdot }} + 4 {1 \over X_{i, \cdot } + n_i - 1} \Big) \Big] \text{,} \label{tmultiplep1} 
\end{align}
where the second equality follows from Lemma \ref{lem:hudson}. 
Furthermore, 
\begin{align}
&{1 \over X_{i, \cdot }} {n_i - 1 \over X_{i, \cdot } + n_i - 1} - 2 {1 \over X_{i, \cdot } + 1} - 2 {1 \over X_{i, \cdot }} + 4 {1 \over X_{i, \cdot } + n_i - 1} \non \\
&= - {1 \over X_{i, \cdot }} {n_i - 1 \over X_{i, \cdot } + n_i - 1} - 2 {1 \over X_{i, \cdot } + 1} + 2 {1 \over X_{i, \cdot } + n_i - 1} < 0 \non 
\end{align}
when $X_{i, \cdot } \ge 1$. 
Thus, (\ref{tmultiplep1}) is negative and this proves part (i). 

For part (ii), let $\De _2 ( \bla ) = E_{\bla } [ L( \blah ( \X , \Y ), \bla ) ] - E_{\bla } [ L( \blah ^{\rm{O}} ( \X , \Y ), \bla ) ]$, let $(( \lah _{i, j} ( \X , \Y ))_{j = 1, \dots , n_i} )_{i = 1, \dots , m} = \blah ( \X , \Y )$, and let 
\begin{align}
\blat ( \X , \Y ) = (( \lat _{i, j} ( \X , \Y ))_{j = 1, \dots , n_i} )_{i = 1, \dots , m} = \Big( \Big( {X_{i, \cdot } + Y_i \over 2} {X_{i, j} \over X_{i, \cdot } + n_i - 1} \Big) _{j = 1, \dots , n_i} \Big) _{i = 1, \dots , m} \text{.} \non 
\end{align}
Then, since 
\begin{align}
&\lat _{i, j} ( \X , \Y ) = \Big( 1 - {n_i - 1 \over X_{i, \cdot } + Y_i + n_i - 1} \Big) \lah _{i, j}^{\rm{O}} ( \X , \Y ) \non 
\end{align}
and 
\begin{align}
&\lah _{i, j} ( \X , \Y ) = \Big( 1 - {m - 1 \over X_{\cdot , \cdot } + Y_{\cdot } + m - 1} \Big) \lat _{i, j} ( \X , \Y ) \non 
\end{align}
for all $j = 1, \dots , n_i$ for all $i = 1, \dots , m$, we have 
\begin{align}
\De _2 ( \bla ) &= E_{\bla } [ L( \blah ( \X , \Y ), \bla ) ] - E_{\bla } [ L( \blat ( \X , \Y ), \bla ) ] + E_{\bla } [ L( \blat ( \X , \Y ), \bla ) ] - E_{\bla } [ L( \blah ^{\rm{O}} ( \X , \Y ), \bla ) ] \non \\
&= E \Big[ \sum_{i = 1}^{m} \sum_{j = 1}^{n_i} {1 \over \la _{i, j}} \Big[ \{ \lat _{i, j} ( \X , \Y ) \} ^2 \Big\{ \Big( {m - 1 \over X_{\cdot , \cdot } + Y_{\cdot } + m - 1} \Big) ^2 - 2 {m - 1 \over X_{\cdot , \cdot } + Y_{\cdot } + m - 1} \Big\} \non \\
&\quad + 2 \la _{i, j} \lat _{i, j} ( \X , \Y ) {m - 1 \over X_{\cdot , \cdot } + Y_{\cdot } + m - 1} \non \\
&\quad + \{ \lah _{i, j}^{\rm{O}} ( \X , \Y ) \} ^2 \Big\{ \Big( {n_i - 1 \over X_{i, \cdot } + Y_i + n_i - 1} \Big) ^2 - 2 {n_i - 1 \over X_{i, \cdot } + Y_i + n_i - 1} \Big\} \non \\
&\quad + 2 \la _{i, j} \lah _{i, j}^{\rm{O}} ( \X , \Y ) {n_i - 1 \over X_{i, \cdot } + Y_i + n_i - 1} \Big] \Big] \non \\
&= E \Big[ \sum_{i = 1}^{m} \sum_{j = 1}^{n_i} \Big[ {X_{i, j} + 1 \over ( X_{i, \cdot } + n_i )^2} {( X_{i, \cdot } + Y_i + 1)^2 \over 4} \Big\{ \Big( {m - 1 \over X_{\cdot , \cdot } + Y_{\cdot } + m} \Big) ^2 - 2 {m - 1 \over X_{\cdot , \cdot } + Y_{\cdot } + m} \Big\} \non \\
&\quad + 2 {X_{i, \cdot } + Y_i \over 2} {X_{i, j} \over X_{i, \cdot } + n_i - 1} {m - 1 \over X_{\cdot , \cdot } + Y_{\cdot } + m - 1} \non \\
&\quad + {X_{i, j} + 1 \over ( X_{i, \cdot } + n_i )^2} {( X_{i, \cdot } + Y_i + n_i )^2 \over 4} \Big\{ \Big( {n_i - 1 \over X_{i, \cdot } + Y_i + n_i} \Big) ^2 - 2 {n_i - 1 \over X_{i, \cdot } + Y_i + n_i} \Big\} \non \\
&\quad + 2 {X_{i, \cdot } + Y_i + n_i - 1 \over 2} {X_{i, j} \over X_{i, \cdot } + n_i - 1} {n_i - 1 \over X_{i, \cdot } + Y_i + n_i - 1} \Big] \Big] \text{,} \non 
\end{align}
where the last equality follows from Lemma \ref{lem:hudson}. 
Therefore, 
\begin{align}
\De _2 ( \bla ) &= E \Big[ \sum_{i = 1}^{m} \Big[ {1 \over X_{i, \cdot } + n_i} {( X_{i, \cdot } + Y_i + 1)^2 \over 4} \Big\{ \Big( {m - 1 \over X_{\cdot , \cdot } + Y_{\cdot } + m} \Big) ^2 - 2 {m - 1 \over X_{\cdot , \cdot } + Y_{\cdot } + m} \Big\} \non \\
&\quad + 2 {X_{i, \cdot } + Y_i \over 2} {X_{i, \cdot } \over X_{i, \cdot } + n_i - 1} {m - 1 \over X_{\cdot , \cdot } + Y_{\cdot } + m - 1} \non \\
&\quad + {1 \over X_{i, \cdot } + n_i} {( X_{i, \cdot } + Y_i + n_i )^2 \over 4} \Big\{ \Big( {n_i - 1 \over X_{i, \cdot } + Y_i + n_i} \Big) ^2 - 2 {n_i - 1 \over X_{i, \cdot } + Y_i + n_i} \Big\} \non \\
&\quad + 2 {X_{i, \cdot } + Y_i + n_i - 1 \over 2} {X_{i, \cdot } \over X_{i, \cdot } + n_i - 1} {n_i - 1 \over X_{i, \cdot } + Y_i + n_i - 1} \Big] \Big] \non \\
&= E \Big[ \sum_{i = 1}^{m} \Big[ {(m - 1) / 4 \over X_{i, \cdot } + n_i} {( X_{i, \cdot } + Y_i + 1)^2 \over ( X_{\cdot , \cdot } + Y_{\cdot } + m)^2} \{ m - 1 - 2 ( X_{\cdot , \cdot } + Y_{\cdot } + m) \} + {X_{i, \cdot } \over X_{i, \cdot } + n_i - 1} {(m - 1) ( X_{i, \cdot } + Y_i ) \over X_{\cdot , \cdot } + Y_{\cdot } + m - 1} \non \\
&\quad + {( n_i - 1) / 4 \over X_{i, \cdot } + n_i} \{ n_i - 1 - 2 ( X_{i, \cdot } + Y_i + n_i ) \} + {( n_i - 1) X_{i, \cdot } \over X_{i, \cdot } + n_i - 1} \Big] \Big] \text{.} \non 
\end{align}

Let $W_i = X_{i, \cdot } + Y_i \sim {\rm{Po}} (2 \la _{i, \cdot } )$ for $i = 1, \dots , m$ and let $\W = ( W_i )_{i = 1, \dots , m}$ and $W_{\cdot } = \sum_{i = 1}^{m} W_i$. 
Then, since $X_{i, \cdot } | W_i \sim {\rm{Bin}} ( W_i , 1 / 2)$ for all $i = 1, \dots , m$, by Jensen's inequality, 
\begin{align}
\De _2 ( \bla ) &\le E \Big[ \sum_{i = 1}^{m} I_{1, i} ( \W ) + \sum_{i = 1}^{m} I_{2, i} ( \W ) \Big] \text{,} \label{tmultiplep1.5} 
\end{align}
where 
\begin{align}
I_{1, i} ( \W ) &= {(m - 1) / 4 \over W_i / 2 + n_i} {( W_i + 1)^2 \over ( W_{\cdot } + m)^2} \{ m - 1 - 2 ( W_{\cdot } + m) \} + {W_i / 2 \over W_i / 2 + n_i - 1} {(m - 1) W_i \over W_{\cdot } + m - 1} \non 
\end{align}
and 
\begin{align}
I_{2, i} ( \W ) &= {( n_i - 1) / 4 \over W_i / 2 + n_i} \{ n_i - 1 - 2 ( W_i + n_i ) \} + {( n_i - 1) W_i / 2 \over W_i / 2 + n_i - 1} \non 
\end{align}
for $i = 1, \dots , m$. 

First, 
\begin{align}
I_{1, i} ( \W ) %
&= {(m - 1)^2 / 2 \over W_i + 2 n_i} {( W_i + 1)^2 \over ( W_{\cdot } + m)^2} - {m - 1 \over W_{\cdot } + m} {( W_i + 1)^2 \over W_i + 2 n_i} + {m - 1 \over W_{\cdot } + m - 1} {{W_i}^2 \over W_i + 2 n_i - 2} \non 
\end{align}
for all $i = 1, \dots , m$. 
Note that 
\begin{align}
&{( W_i + 1)^2 \over W_i + 2 n_i} = ( W_i + 1) \Big( 1 - {2 n_i - 1 \over W_i + 2 n_i} \Big) \quad \text{and} \quad {{W_i}^2 \over W_i + 2 n_i - 2} = W_i \Big( 1 - {2 n_i - 2 \over W_i + 2 n_i - 2} \Big) \non 
\end{align}
for all $i = 1, \dots , m$. 
Then 
\begin{align}
&\sum_{i = 1}^{m} I_{1, i} ( \W ) \non \\
&= \sum_{i = 1}^{m} \Big\{ {(m - 1)^2 / 2 \over W_i + 2 n_i} {( W_i + 1)^2 \over ( w_{\cdot } + m)^2} - {m - 1 \over W_{\cdot } + m} ( W_i + 1) \Big( 1 - {2 n_i - 1 \over W_i + 2 n_i} \Big) + {m - 1 \over W_{\cdot } + m - 1} W_i \Big( 1 - {2 n_i - 2 \over W_i + 2 n_i - 2} \Big) \Big\} \non \\
&= \sum_{i = 1}^{m} {(m - 1)^2 / 2 \over W_i + 2 n_i} {( W_i + 1)^2 \over ( W_{\cdot } + m)^2} \non \\
&\quad - (m - 1) + \sum_i {m - 1 \over W_{\cdot } + m} ( W_i + 1) {2 n_i - 1 \over W_i + 2 n_i} + {m - 1 \over W_{\cdot } + m - 1} W_{\cdot } - \sum_{i = 1}^{m} {m - 1 \over W_{\cdot } + m - 1} W_i {2 n_i - 2 \over W_i + 2 n_i - 2} \non \\
&= (m - 1)^2 J_1 ( \W ) + (m - 1) J_2 ( \W ) \text{,} \non 
\end{align}
where 
\begin{align}
J_1 ( \W ) = \sum_{i = 1}^{m} {1 / 2 \over W_i + 2 n_i} {( W_i + 1)^2 \over ( W_{\cdot } + m)^2} - {1 \over W_{\cdot } + m - 1} \non 
\end{align}
and 
\begin{align}
J_2 ( \W ) = \sum_{i = 1}^{m} \Big( {W_i + 1 \over W_{\cdot } + m} {2 n_i - 1 \over W_i + 2 n_i} - {W_i \over W_{\cdot } + m - 1} {2 n_i - 2 \over W_i + 2 n_i - 2} \Big) \text{.} \non 
\end{align}
By the covariance inequality, 
\begin{align}
J_1 ( \W ) %
&\le m \Big\{ {1 \over m} \sum_{i = 1}^{m} {1 / 2 \over W_i + 2 n_i} \Big\} {1 \over m} \sum_{i = 1}^{m} {( W_i + 1)^2 \over ( W_{\cdot } + m)^2} - {1 \over W_{\cdot } + m - 1} \non \\
&\le {1 \over m} \Big\{ \sum_{i = 1}^{m} {1 / 2 \over W_i + 2 n_i} \Big\} \sum_{i = 1}^{m} {W_i + 1 \over W_{\cdot } + m} - {1 \over W_{\cdot } + m - 1} \non \\
&= \sum_{i = 1}^{m} {1 / (2 m) \over W_i + 2 n_i} - {1 \over W_{\cdot } + m - 1} \text{.} \non 
\end{align}
Also, 
\begin{align}
J_2 ( \W ) &\le \sum_{i = 1}^{m} \Big( {W_i + 1 \over W_{\cdot } + m} {2 n_i - 1 \over W_i + 2 n_i} - {W_i \over W_{\cdot } + m} {2 n_i - 2 \over W_i + 2 n_i} \Big) \non \\
&= \sum_{i = 1}^{m} {W_i + 2 n_i - 1 \over ( W_{\cdot } + m) ( W_i + 2 n_i )} \le {m \over W_{\cdot } + m} \text{.} \non 
\end{align}
Therefore, 
\begin{align}
\sum_{i = 1}^{m} I_{1, i} ( \W ) &\le \sum_{i = 1}^{m} {(m - 1)^2 / (2 m) \over W_i + 2 n_i} - {(m - 1)^2 \over W_{\cdot } + m - 1} + {(m - 1) m \over W_{\cdot } + m} \non \\
&= (m - 1) \sum_{i = 1}^{m} {(m - 1) / (2 m) \over W_i + 2 n_i} + (m - 1) {W_{\cdot } \over ( W_{\cdot } + m - 1) ( W_{\cdot } + m)} \text{.} \label{tmultiplep2} 
\end{align}
Next, it can be shown that 
\begin{align}
I_{2, i} ( \W ) %
&= - ( n_i - 1) {n_i - 3 \over 2} {W_i + 2 ( n_i + 1) ( n_i - 1) / ( n_i - 3) \over ( W_i + 2 n_i ) ( W_i + 2 n_i - 2)} \non 
\end{align}
for all $i = 1, \dots , m$. 
Therefore, 
\begin{align}
I_{2, i} ( \W ) &\le - ( n_i - 1) {n_i - 3 \over 2} {1 \over W_i + 2 n_i - 2} \label{tmultiplep3} 
\end{align}
for all $i = 1, \dots , m$. 
Thus, by (\ref{tmultiplep1.5}), (\ref{tmultiplep2}), and (\ref{tmultiplep3}), we have 
\begin{align}
\De _2 ( \bla ) %
&\le E \Big[ \sum_{i = 1}^{m} {(m - 1)^2 / (2 m) \over W_i + 2 n_i} + {(m - 1) W_{\cdot } \over ( W_{\cdot } + m - 1) ( W_{\cdot } + m)} - \sum_{i = 1}^{m} {( n_i - 1) ( n_i - 3 ) / 2 \over W_i + 2 n_i - 2} \Big] \text{.} \non 
\end{align}

Finally, since 
\begin{align}
{(m - 1) W_{\cdot } \over ( W_{\cdot } + m - 1) ( W_{\cdot } + m)} &\le {m - 1 \over W_{\cdot } + m} = {(m - 1) / m \over W_{\cdot } / m + 1} \le {1 \over m} \sum_{i = 1}^{m} {(m - 1) / m \over W_i + 1} \le {1 \over m} \sum_{i = 1}^{m} {(m - 1) / m + 2 n_i - 3 \over W_i + 2 n_i - 2} \text{,} \non 
\end{align}
it follows that 
\begin{align}
\De &< E \Big[ \sum_{i = 1}^{m} \Big\{ {(m - 1)^2 / (2 m) \over W_i + 2 n_i - 2} + {(m - 1) / m^2 + (2 n_i - 3) / m \over W_i + 2 n_i - 2} - {( n_i - 1) ( n_i - 3) / 2 \over W_i + 2 n_i - 2} \Big\} \Big] \text{.} \non 
\end{align}
The right-hand side of the above inequality is nonpositive by assumption. 
This completes the proof. 
\hfill$\Box$

\bigskip

\noindent
{\bf Proof of Lemma \ref{lem:jeff}.} \ \ Here, we consider the case where we observe $\X $, $\Y $, and $Z$. 
Note that 
\begin{align}
&{\pd \over \pd \La } \log p( \X , \Y , Z | \bla ) = {X_{\cdot , \cdot } + Y_{\cdot } + Z \over \La } - 3 \text{,} \non \\
&{\pd \over \pd \th _i} \log p( \X , \Y , Z | \bla ) = {X_{i, \cdot } + Y_i \over \th _i} - {X_{m, \cdot } + Y_m \over \th _m} \quad \text{for $i = 1, \dots , m - 1$} \text{,} \quad \text{and} \non \\
&{\pd \over \pd \rho _{i, j}} \log p( \X , \Y , Z | \bla ) = {X_{i, j} \over \rho _{i, j}} - {X_{i, n_i} \over \rho _{i, n_i}} \quad \text{for $j = 1, \dots , n_i - 1$ for $i = 1, \dots , m$} \text{.} \non 
\end{align}
Then 
\begin{align}
&{\pd ^2 \over ( \pd \La )^2} \log p( \X , \Y , Z | \bla ) = - {X_{\cdot , \cdot } + Y_{\cdot } + Z \over \La ^2} \text{,} \non 
\end{align}
\begin{align}
&{\pd ^2 \over ( \pd \th _i ) ( \pd \th _{i'} )} \log p( \X , \Y , Z | \bla ) = - 1(i = i' ) {X_{i, \cdot } + Y_i \over {\th _i}^2} - {X_{m, \cdot } + Y_m \over {\th _m}^2} \non 
\end{align}
for $i, i' = 1, \dots , m - 1$, and 
\begin{align}
&{\pd ^2 \over ( \pd \rho _{i, j} ) ( \pd \rho _{i, j'} )} \log p( \X , \Y , Z | \bla ) = - 1(j = j' ) {X_{i, j} \over {\rho _{i, j}}^2} - {X_{i, n_i} \over {\rho _{i, n_i}}^2} \non 
\end{align}
for $j, j' = 1, \dots , n_i - 1$ for $i = 1, \dots , m$. 
Therefore, 
\begin{align}
\pi ^{\rm{J}} ( \La , \bth , \bro ) &\propto \sqrt{E_{\bla } \Big[ {X_{\cdot , \cdot } + Y_{\cdot } + Z \over \La ^2} \Big] } \non \\
&\quad \times \sqrt{\Big| E_{\bla } \Big[ \diag \Big( {X_{1, \cdot } + Y_1 \over {\th _1}^2}, \dots , {X_{m - 1, \cdot } + Y_{m - 1} \over {\th _{m - 1}}^2} \Big) + \j ^{(m - 1)} ( \j ^{(m - 1)} )^{\top } {X_{m, \cdot } + Y_m \over {\th _m}^2} \Big] \Big| } \non \\
&\quad \times \prod_{i = 1}^{m} \sqrt{\Big| E_{\bla } \Big[ \diag \Big( {X_{i, 1} \over {\rho _{i, 1}}^2}, \dots , {X_{i, n_i - 1} \over {\rho _{i, n_i - 1}}^2} \Big) + \j ^{( n_i - 1)} ( \j ^{( n_i - 1)} )^{\top } {X_{i, n_i} \over {\rho _{i, n_i}}^2} \Big] \Big| } \non \\
&\propto \La ^{n_{\cdot } / 2 - 1} \sqrt{\Big| \diag \Big( {1 \over \th _1}, \dots , {1 \over \th _{m - 1}} \Big) + \j ^{(m - 1)} ( \j ^{(m - 1)} )^{\top } {1 \over \th _m} \Big| } \Big\{ \prod_{i = 1}^{m} {\th _i}^{( n_i - 1) / 2} \Big\} \non \\
&\quad \times \prod_{i = 1}^{m} \sqrt{\Big| \diag \Big( {1 \over \rho _{i, 1}}, \dots , {1 \over \rho _{i, n_i - 1}} \Big) + \j ^{( n_i - 1)} ( \j ^{( n_i - 1)} )^{\top } {1 \over \rho _{i, n_i}} \Big| } \text{.} \non 
\end{align}
Since $| \A + \v \v ^{\top } | = | \A | (1 + \v ^{\top } \A ^{- 1} \v ^{\top } )$ for any $n \times n$ nonsingular matrix $\A $ and any $\v \in \mathbb{R} ^n$ for any $n \in \mathbb{N}$, it follows that 
\begin{align}
\pi ^{\rm{J}} ( \La , \bth , \bro ) &\propto \La ^{n_{\cdot } / 2 - 1} \sqrt{\Big( \prod_{i = 1}^{n_i - 1} {1 \over \th _i} \Big) \Big( 1 + {1 \over \th _m} \sum_{i = 1}^{m - 1} \th _i \Big) } \Big\{ \prod_{i = 1}^{m} {\th _i}^{( n_i - 1) / 2} \Big\} \non \\
&\quad \times \prod_{i = 1}^{m} \sqrt{ \Big( \prod_{j = 1}^{n_i} {1 \over \rho _{i, j}} \Big) \Big( 1 + {1 \over \rho _{i, n_i}} \sum_{j = 1}^{n_i - 1} \rho _{i, j} \Big) } \non \\
&= \La ^{n_{\cdot } / 2 - 1} \prod_{i = 1}^{m} \Big( {\th _i}^{n_i / 2 - 1} \prod_{j = 1}^{n_i} {\rho _{i, j}}^{1 / 2 - 1} \Big) \text{,} \non 
\end{align}
which is the desired result. 
\hfill$\Box$

\bigskip

\noindent
{\bf Proof of Theorem \ref{thm:KL_0}.} \ \ Let $\De _{3}^{( \al ; \a )} ( \bla ) = E_{\bla } [ L^{\rm{KL}} ( \blah ^{( \al ; \a )} ( \X , \Y ), \bla ) ] - E_{\bla } [ L^{\rm{KL}} ( \blah ^{\rm{J}} ( \X , \Y ), \bla ) ]$. 
Let $(( \lah _{i, j}^{\rm{J}} ( \X , \Y ))_{j = 1, \dots , n_i} )_{i = 1, \dots , m} = \blah ^{\rm{J}} ( \X , \Y )$ and $(( \lah _{i, j}^{( \al ; \a )} ( \X , \Y ))_{j = 1, \dots , n_i} )_{i = 1, \dots , m} = \blah ^{( \al ; \a )} ( \X , \Y )$. 
Let 
\begin{align}
\blat ^{( \a )} ( \X , \Y ) &= (( \lat _{i, j}^{( \a )} ( \X , \Y ))_{j = 1, \dots , n_i} )_{i = 1, \dots , m} = \Big( \Big( {X_{i, \cdot } + Y_i + a_i \over 2} {X_{i, j} + 1 / 2 \over X_{i, \cdot } + n_i / 2} \Big) _{j = 1, \dots , n_i} \Big) _{i = 1, \dots , m} \text{.} \non 
\end{align}
Then 
\begin{align}
\De _{3}^{( \al ; \a )} ( \bla ) &= E_{\bla } [ L^{\rm{KL}} ( \blah ^{( \al ; \a )} ( \X , \Y ), \bla ) ] - E_{\bla } [ L^{\rm{KL}} ( \blat ^{( \a )} ( \X , \Y ), \bla ) ] \non \\
&\quad + E_{\bla } [ L^{\rm{KL}} ( \blat ^{( \a )} ( \X , \Y ), \bla ) ] - E_{\bla } [ L^{\rm{KL}} ( \blah ^{\rm{J}} ( \X , \Y ), \bla ) ] \non \\
&= E_{\bla } \Big[ \sum_{i = 1}^{m} \sum_{j = 1}^{n_i} \Big\{ \lah _{i, j}^{( \al ; \a )} ( \X , \Y ) - \lat _{i, j}^{( \a )} ( \X , \Y ) - \la _{i, j} \log {\lah _{i, j}^{( \al ; \a )} ( \X , \Y ) \over \lat _{i, j}^{( \a )} ( \X , \Y )} \Big\} \Big] \non \\
&\quad + E_{\bla } \Big[ \sum_{i = 1}^{m} \sum_{j = 1}^{n_i} \Big\{ \lat _{i, j}^{( \a )} ( \X , \Y ) - \lah _{i, j}^{\rm{J}} ( \X , \Y ) - \la _{i, j} \log {\lat _{i, j}^{( \al ; \a )} ( \X , \Y ) \over \lah _{i, j}^{\rm{J}} ( \X , \Y )} \Big\} \Big] \non \\
&= E_{\bla } \Big[ - {a_{\cdot } - \al \over 2} + \La \log {X_{\cdot , \cdot } + Y_{\cdot } + a_{\cdot } \over X_{\cdot , \cdot } + Y_{\cdot } + \al } \Big\} \Big] + E_{\bla } \Big[ - {n_{\cdot } / 2 - a_{\cdot } \over 2} + \sum_{i = 1}^{m} \la _{i, \cdot } \log {X_{i, \cdot } + Y_i + n_i / 2 \over X_{i, \cdot } + Y_i + a_i} \Big] \text{.} \non 
\end{align}
By Lemma \ref{lem:hudson} and by assumption, it follows that 
\begin{align}
\De _{3}^{( \al ; \a )} ( \bla ) &= E_{\bla } \Big[ - {a_{\cdot } - \al \over 2} + {X_{\cdot , \cdot } + Y_{\cdot } \over 2} \log {X_{\cdot , \cdot } + Y_{\cdot } + a_{\cdot } - 1 \over X_{\cdot , \cdot } + Y_{\cdot } + \al - 1} \Big\} \Big] \non \\
&\quad + E_{\bla } \Big[ - {n_{\cdot } / 2 - a_{\cdot } \over 2} + \sum_{i = 1}^{m} {X_{i, \cdot } + Y_i \over 2} \log {X_{i, \cdot } + Y_i + n_i / 2 - 1 \over X_{i, \cdot } + Y_i + a_i - 1} \Big] \non \\
&< E_{\bla } \Big[ - {a_{\cdot } - \al \over 2} + {X_{\cdot , \cdot } + Y_{\cdot } \over 2} {a_{\cdot } - \al \over X_{\cdot , \cdot } + Y_{\cdot } + \al - 1} \Big\} \Big] \non \\
&\quad + E_{\bla } \Big[ - {n_{\cdot } / 2 - a_{\cdot } \over 2} + \sum_{i = 1}^{m} {X_{i, \cdot } + Y_i \over 2} {n_i / 2 - a_i \over X_{i, \cdot } + Y_i + a_i - 1} \Big] \le 0 \text{,} \non 
\end{align}
where the first inequality is strict since $a_{\cdot } - \al + \sum_{i = 1}^{m} ( n_i / 2 - a_i ) > 0$. 
This proves the theorem. 
\hfill$\Box$

\bigskip

\noindent
{\bf Proof of Theorem \ref{thm:KL}.} \ \ Let $\De _{4}^{( \al ; \a )} ( \bla ) = E_{\bla } [ L^{\rm{KL}} ( \blah ^{( \al ; \a )} ( \X , \Y , Z) , \bla ) ] - E_{\bla } [ L^{\rm{KL}} ( \blah ^{\rm{J}} ( \X , \Y , Z) , \bla ) ]$. 
Let $(( \lah _{i, j}^{\rm{J}} ( \X , \Y , Z))_{j = 1, \dots , n_i} )_{i = 1, \dots , m} = \blah ^{\rm{J}} ( \X , \Y , Z)$ and $(( \lah _{i, j}^{( \al ; \a )} ( \X , \Y , Z))_{j = 1, \dots , n_i} )_{i = 1, \dots , m} = \blah ^{( \al ; \a )} ( \X , \Y , Z)$. 
Then 
\begin{align}
\De _{4}^{( \al ; \a )} ( \bla ) &= E \Big[ \sum_{i = 1}^{m} \sum_{j = 1}^{n_i} \Big\{ \lah _{i, j}^{( \al ; \a )} ( \X , \Y , Z) - \lah _{i, j}^{\rm{J}} ( \X , \Y , Z) - \la _{i, j} \log {\lah _{i, j}^{( \al ; \a )} ( \X , \Y , Z) \over \lah _{i, j}^{\rm{J}} ( \X , \Y , Z)} \Big\} \Big] \non \\
&= E \Big[ - {n_{\cdot } / 2 - \al \over 3} + \sum_{i = 1}^{m} \sum_{j = 1}^{n_i} \Big\{ \la _{i, j} \log {\lah _{i, j}^{\rm{J}} ( \X , \Y , Z) \over \lah _{i, j}^{( \al ; \a )} ( \X , \Y , Z)} \Big\} \Big] \non \\
&= E \Big[ - {n_{\cdot } / 2 - \al \over 3} + \La \Big( \log {X_{\cdot , \cdot } + Y_{\cdot } + Z + n_{\cdot } / 2 \over X_{\cdot , \cdot } + Y_{\cdot } + Z + \al } + \log {X_{\cdot , \cdot } + Y_{\cdot } + a_{\cdot } \over X_{\cdot , \cdot } + Y_{\cdot } + n_{\cdot } / 2} \Big) \non \\
&\quad + \sum_{i = 1}^{m} \la _{i, \cdot } \log {X_{i, \cdot } + Y_i + n_i / 2 \over X_{i, \cdot } + Y_i + a_i} \Big] \text{.} \non 
\end{align}
By Lemma \ref{lem:hudson}, 
\begin{align}
\De _{4}^{( \al ; \a )} ( \bla ) &= \De _{4, 1}^{( \al ; \a )} ( \bla ) + \De _{4, 2}^{( \al ; \a )} ( \bla ) \text{,} \non 
\end{align}
where 
\begin{align}
&\De _{4, 1}^{( \al ; \a )} ( \bla ) = E \Big[ - {n_{\cdot } / 2 - \al \over 3} + {X_{\cdot , \cdot } + Y_{\cdot } + Z \over 3} \log {X_{\cdot , \cdot } + Y_{\cdot } + Z + n_{\cdot } / 2 - 1 \over X_{\cdot , \cdot } + Y_{\cdot } + Z + \al - 1} \Big] \quad \text{and} \non \\
&\De _{4, 2}^{( \al ; \a )} ( \bla ) = E \Big[ {X_{\cdot , \cdot } + Y_{\cdot } \over 2} \log {X_{\cdot , \cdot } + Y_{\cdot } + a_{\cdot } - 1 \over X_{\cdot , \cdot } + Y_{\cdot } + n_{\cdot } / 2 - 1} + \sum_{i = 1}^{m} {X_{i, \cdot } + Y_i \over 2} \log {X_{i, \cdot } + Y_i + n_i / 2 - 1 \over X_{i, \cdot } + Y_i + a_i - 1} \Big] \text{.} \non 
\end{align}

For $\De _{4, 1}^{( \al ; \a )} ( \bla )$, we have 
\begin{align}
\De _{4, 1}^{( \al ; \a )} ( \bla ) &< E \Big[ - {n_{\cdot } / 2 - \al \over 3} + {X_{\cdot , \cdot } + Y_{\cdot } + Z \over 3}  {n_{\cdot } / 2 - \al \over X_{\cdot , \cdot } + Y_{\cdot } + Z + \al - 1} \Big] \non \\
&= E \Big[ - {n_{\cdot } / 2 - \al \over 3} {\al - 1 \over X_{\cdot , \cdot } + Y_{\cdot } + Z + \al - 1} \Big] \text{.} \non 
\end{align}
Meanwhile, since $n_i / 2 = n_{\cdot } / (2 m)$ and $a_i = a_{\cdot } / m$ for all $i = 1, \dots , m$ by assumption, 
\begin{align}
\sum_{i = 1}^{m} {X_{i, \cdot } + Y_i \over 2} \log {X_{i, \cdot } + Y_i + n_i / 2 - 1 \over X_{i, \cdot } + Y_i + a_i - 1} = {1 \over m} \sum_{i = 1}^{m} {m ( X_{i, \cdot } + Y_i ) \over 2} \log {m ( X_{i, \cdot } + Y_i ) + n_{\cdot } / 2 - m \over m ( X_{i, \cdot } + Y_i ) + a_{\cdot } - m} \text{.} \non 
\end{align}
Since $w \log \{ (w + n_{\cdot } / 2 - m) / (w + a_{\cdot } - m) \} $ is a concave function of $w \in [0, \infty )$, by Jensen's inequality, 
\begin{align}
&{1 \over m} \sum_{i = 1}^{m} {m ( X_{i, \cdot } + Y_i ) \over 2} \log {m ( X_{i, \cdot } + Y_i ) + n_{\cdot } / 2 - m \over m ( X_{i, \cdot } + Y_i ) + a_{\cdot } - m} \le {X_{\cdot , \cdot } + Y_{\cdot } \over 2} \log {X_{\cdot , \cdot } + Y_{\cdot } + n_{\cdot } / 2 - m \over X_{\cdot , \cdot } + Y_{\cdot} + a_{\cdot } - m} \text{.} \non 
\end{align}
Therefore, 
\begin{align}
\De _{4, 2}^{( \al ; \a )} ( \bla ) &\le E \Big[ {X_{\cdot , \cdot } + Y_{\cdot } \over 2} \log \Big\{ \Big( 1 - {n_{\cdot } / 2 - a_{\cdot } \over X_{\cdot , \cdot } + Y_{\cdot } + n_{\cdot } / 2 - 1} \Big) \Big( 1 + {n_{\cdot } / 2 - a_{\cdot } \over X_{\cdot , \cdot } + Y_{\cdot} + a_{\cdot } - m} \Big) \Big\} \Big] \text{.} \non 
\end{align}
Since 
\begin{align}
&\Big( 1 - {n_{\cdot } / 2 - a_{\cdot } \over X_{\cdot , \cdot } + Y_{\cdot } + n_{\cdot } / 2 - 1} \Big) \Big( 1 + {n_{\cdot } / 2 - a_{\cdot } \over X_{\cdot , \cdot } + Y_{\cdot} + a_{\cdot } - m} \Big) \non \\
&= 1 + {(m - 1) ( n_{\cdot } / 2 - a_{\cdot } ) \over ( X_{\cdot , \cdot } + Y_{\cdot } + n_{\cdot } / 2 - 1) ( X_{\cdot , \cdot } + Y_{\cdot} + a_{\cdot } - m)} \text{,} \non 
\end{align}
it follows that 
\begin{align}
\De _{4, 2}^{( \al ; \a )} ( \bla ) &\le E \Big[ {X_{\cdot , \cdot } + Y_{\cdot } \over 2} {(m - 1) ( n_{\cdot } / 2 - a_{\cdot } ) \over ( X_{\cdot , \cdot } + Y_{\cdot } + n_{\cdot } / 2 - 1) ( X_{\cdot , \cdot } + Y_{\cdot} + a_{\cdot } - m)} \Big] \text{.} \non 
\end{align}
Thus, 
\begin{align}
\De _{4}^{( \al ; \a )} ( \bla ) &< E \Big[ - {n_{\cdot } / 2 - \al \over 3} {\al - 1 \over X_{\cdot , \cdot } + Y_{\cdot } + Z + \al - 1} + {X_{\cdot , \cdot } + Y_{\cdot } \over 2} {(m - 1) ( n_{\cdot } / 2 - a_{\cdot } ) \over ( X_{\cdot , \cdot } + Y_{\cdot } + n_{\cdot } / 2 - 1) ( X_{\cdot , \cdot } + Y_{\cdot} + a_{\cdot } - m)} \Big] \text{.} \non 
\end{align}

By Jensen's inequality, 
\begin{align}
&E \Big[ - {n_{\cdot } / 2 - \al \over 3} {\al - 1 \over X_{\cdot , \cdot } + Y_{\cdot } + Z + \al - 1} \Big] \le - {n_{\cdot } / 2 - \al \over 3} {\al - 1 \over 3 \La + \al - 1} \text{.} \non 
\end{align}
Meanwhile, by the covariance inequality, 
\begin{align}
&E \Big[ {X_{\cdot , \cdot } + Y_{\cdot } \over 2} {(m - 1) ( n_{\cdot } / 2 - a_{\cdot } ) \over ( X_{\cdot , \cdot } + Y_{\cdot } + n_{\cdot } / 2 - 1) ( X_{\cdot , \cdot } + Y_{\cdot} + a_{\cdot } - m)} \Big] \non \\
&\le {(m - 1) ( n_{\cdot } / 2 - a_{\cdot } ) \over 2} E \Big[ {X_{\cdot , \cdot } + Y_{\cdot } \over X_{\cdot , \cdot } + Y_{\cdot } + n_{\cdot } / 2 - 1} \Big] E \Big[ {1 \over X_{\cdot , \cdot } + Y_{\cdot} + a_{\cdot } - m} \Big] \non \\
&\le {(m - 1) ( n_{\cdot } / 2 - a_{\cdot } ) \over 2} {2 \La \over 2 \La + n_{\cdot } / 2 - 1} E \Big[ {1 \over X_{\cdot , \cdot } + Y_{\cdot} + a_{\cdot } - m} \Big] \non \\
&= {(m - 1) ( n_{\cdot } / 2 - a_{\cdot } ) \over 2} {1 \over 2 \La + n_{\cdot } / 2 - 1} E \Big[ {X_{\cdot , \cdot } + Y_{\cdot} \over X_{\cdot , \cdot } + Y_{\cdot} + a_{\cdot } - m - 1} \Big] \non \\
&\le {(m - 1) ( n_{\cdot } / 2 - a_{\cdot } ) \over 2} {1 \over 2 \La + n_{\cdot } / 2 - 1} \text{,} \non 
\end{align}
where the equality follows from Lemma \ref{lem:hudson} and where the last inequality follows since $a_{\cdot } - m - 1 \ge 0$ by assumption. 
Thus, 
\begin{align}
&(3 \La + \al - 1) (2 \La + n_{\cdot } / 2 - 1) \De _{4}^{( \al ; \a )} ( \bla ) \non \\
&< (3 \La + \al - 1) (2 \La + n_{\cdot } / 2 - 1) \Big\{ - {n_{\cdot } / 2 - \al \over 3} {\al - 1 \over 3 \La + \al - 1} + {(m - 1) ( n_{\cdot } / 2 - a_{\cdot } ) \over 2} {1 \over 2 \La + n_{\cdot } / 2 - 1} \Big\} \non \\
&= - \Big\{ {n_{\cdot } / 2 - \al \over 3} ( \al - 1) (2 \La + n_{\cdot } / 2 - 1) - {(m - 1) ( n_{\cdot } / 2 - a_{\cdot } ) \over 2} (3 \La + \al - 1) \Big\} \non \\
&= - \Big[ \Big\{ {2 \over 3} ( n_{\cdot } / 2 - \al ) ( \al - 1) - {3 \over 2} (m - 1) ( n_{\cdot } / 2 - a_{\cdot } ) \Big\} \La \non \\
&\quad + \Big\{ {( n_{\cdot } / 2 - \al ) ( n_{\cdot } / 2 - 1) \over 3} - {(m - 1) ( n_{\cdot } / 2 - a_{\cdot } ) \over 2} \Big\} ( \al - 1) \Big] \text{,} \non 
\end{align}
the right-hand side of which is nonpositive by assumptions (\ref{eq:KL_A1}) and (\ref{eq:KL_A2}). 
This completes the proof. 
\hfill$\Box$

\bigskip

\noindent
{\bf Proof of Lemma \ref{lem:general_Jeff}.} \ \ Fix $D' = 0, 1, \dots , D$. 
Note that the likelihood of $\X ( D' )$ is given by 
\begin{align}
p( \X ( D' ) | \bla ) %
&= \Big( \prod_{d = 0}^{D} \prod_{i_1 = 1}^{n_1} \cdots \prod_{i_d = 1}^{n_d} {1 \over \Xt _{( i_1 , \dots , i_d )} !} \Big) \la ^{\Xt _{(+)}} e^{- (1 + D - D' ) \la } \prod_{d = 1}^{D} \prod_{i_1 = 1}^{n_1} \cdots \prod_{i_d = 1}^{n_d} {\th _{( i_1 , \dots , i_d )}}^{\Xt _{( i_1 , \dots , i_d , +)}} \text{,} \non 
\end{align}
where $\Xt _{( i_1 , \dots , i_d )} = 1(d \ge D' ) X_{( i_1 , \dots , i_d )}$ and $\Xt _{( i_1 , \dots , i_d , +)} = \sum_{d' = d}^{D} \sum_{i_{d + 1} = 1}^{n_{d + 1}} \cdots \sum_{i_{d'} = 1}^{n_{d'}} \Xt _{( i_1 , \dots , i_d , i_{d + 1} , \dots , i_{d'} )}$ for $( i_1 , \dots , i_d ) \in \prod_{d' = 1}^{d} \{ 1, \dots , n_{d'} \} $ for $d = 0, 1, \dots , D$. 
Then we have that 
\begin{align}
&{\pd ^2 \over ( \pd \la )^2} \log p( \X ( D' ) | \bla ) = - {\Xt _{(+)} \over \la ^2} \non 
\end{align}
and that 
\begin{align}
&{\pd ^2 \over ( \pd \th _{( i_1 , \dots , i_{d - 1} , i_d )} ) ( \pd \th _{( i_1 , \dots , i_{d - 1} {i_d}' )} )} \log p( \X ( D' ) | \bla ) = - 1( i_d = {i_d}' ) {\Xt _{( i_1 , \dots , i_{d - 1} , i_d , +)} \over {\th _{( i_1 , \dots , i_{d - 1} , i_d )}}^2} - {\Xt _{( i_1 , \dots , i_{d - 1} , n_d , +)} \over {\th _{( i_1 , \dots , i_{d - 1} , n_d )}}^2} \non 
\end{align}
for $i_d , {i_d}' \in \{ 1, \dots , n_d - 1 \} $ for $( i_1 , \dots , i_{d - 1} ) \in \prod_{d' = 1}^{d - 1} \{ 1, \dots , n_{d'} \} $ for $d = 1, \dots , D$. 
Therefore, the Jeffreys prior corresponding to $\X ( D' )$ is given by 
\begin{align}
\pi _{D'}^{\rm{J}} ( \la , \bth ) &\propto \Big( E_{\bla } \Big[ {\Xt _{(+)} \over \la ^2} \Big] \Big) ^{1 / 2} \prod_{d = 1}^{D} \prod_{i_1 = 1}^{n_1} \cdots \prod_{i_{d - 1} = 1}^{n_{d - 1}} \Big| E_{\bla } \Big[ \diag \Big( {\Xt _{( i_1 , \dots , i_{d - 1} , 1, +)} \over {\th _{( i_1 , \dots , i_{d - 1} , 1)}}^2} , \dots , {\Xt _{( i_1 , \dots , i_{d - 1} , n_d - 1, +)} \over {\th _{( i_1 , \dots , i_{d - 1} , n_d - 1)}}^2} \Big) \non \\
&\quad + \j ^{( n_d - 1)} {\Xt _{( i_1 , \dots , i_{d - 1} , n_d , +)} \over {\th _{( i_1 , \dots , i_{d - 1} , n_d )}}^2} {\j ^{( n_d - 1)}}^{\top } \Big] \Big| ^{1 / 2} \non \\
&\propto \Big( {1 \over \la } \Big) ^{1 / 2} \prod_{d = 1}^{D} \prod_{i_1 = 1}^{n_1} \cdots \prod_{i_{d - 1} = 1}^{n_{d - 1}} \Big[ \Big\{ \la ^{( n_d - 1) / 2} \prod_{d' = 1}^{d - 1} {\th _{( i_1 , \dots , i_{d'} )}}^{( n_d - 1) / 2} \Big\} \non \\
&\quad \times \Big| \diag \Big( {1 \over \th _{( i_1 , \dots , i_{d - 1} , 1)}} , \dots , {1 \over \th _{( i_1 , \dots , i_{d - 1} , n_d - 1)}} \Big) + \j ^{( n_d - 1)} {1 \over \th _{( i_1 , \dots , i_{d - 1} , n_d )}} {\j ^{( n_d - 1)}}^{\top } \Big| ^{1 / 2} \Big] \non \\
&= \Big( {1 \over \la } \Big) ^{1 / 2} \prod_{d = 1}^{D} \prod_{i_1 = 1}^{n_1} \cdots \prod_{i_{d - 1} = 1}^{n_{d - 1}} \Big[ \Big\{ \la ^{( n_d - 1) / 2} \prod_{d' = 1}^{d - 1} {\th _{( i_1 , \dots , i_{d'} )}}^{( n_d - 1) / 2} \Big\} \prod_{i_d = 1}^{n_d} \Big( {1 \over \th _{( i_1 , \dots , i_{d - 1} , i_d )}} \Big) ^{1 / 2} \Big] \non \\
&= \la ^{n_1 \dotsm n_D / 2 - 1}  \prod_{d = 1}^{D} \prod_{i_1 = 1}^{n_1} \cdots \prod_{i_d = 1}^{n_d} {\th _{( i_1 , \dots , i_d )}}^{n_{d + 1} \dotsm n_D / 2 - 1} \text{,} \non 
\end{align}
which is the desired result. 
\hfill$\Box$

\bigskip

\noindent
{\bf Proof of Lemma \ref{lem:general_estimator}.} \ \ %
Note that the posterior density of $( \la , \bth ) | \X ( D' )$ is given by 
\begin{align}
\pi ^{( \a )} ( \la , \bth | \X ) &\propto \la ^{\Xt _{(+)} ( D' ) + a_0 - 1} e^{- (1 + D - D' ) \la } \prod_{d = 1}^{D} \prod_{i_1 = 1}^{n_1} \cdots \prod_{i_d = 1}^{n_d} {\th _{( i_1 , \dots , i_d )}}^{\Xt _{( i_1 , \dots , i_d , +)} ( D' ) + a_d - 1} \text{.} \non 
\end{align}
Then for any $\dbt = ( \cdots ( \dt _{i_1 , \dots , i_D} )_{i_D = 1, \dots , n_D} \cdots )_{i_1 = 1, \dots , n_1} \in (0, \infty )^{n_1 \dotsm n_D}$, 
\begin{align}
&E_{\pi ^{( \a )}} [ \widetilde{L} ^{\rm{KL}} ( \dbt , \bla ) | \X ( D' ) ] \non \\
&= E_{\pi ^{( \a )}} \Big[ \sum_{d = 0}^{D} \sum_{i_1 = 1}^{n_1} \cdots \sum_{i_d = 1}^{n_d} \Big( \dt _{( i_1 , \dots , i_d )} - \la _{( i_1 , \dots , i_d )} - \la _{( i_1 , \dots , i_d )} \log {\dt _{( i_1 , \dots , i_d )} \over \la _{( i_1 , \dots , i_d )}} \Big) \Big| \X ( D' ) \Big] \non \\
&= E_{\pi ^{( \a )}} \Big[ (1 + D) ( \dt - \la ) - \sum_{d = 0}^{D} \sum_{i_1 = 1}^{n_1} \cdots \sum_{i_d = 1}^{n_d} \la _{( i_1 , \dots , i_d )} \Big( \log {\dt \over \la } + \sum_{d' = 1}^{d} \log {\pt _{( i_1 , \dots , i_{d'} )} \over \th _{( i_1 , \dots , i_{d'} )}} \Big) \Big| \X ( D' ) \Big] \non \\
&= E_{\pi ^{( \a )}} \Big[ (1 + D) \Big( \dt - \la - \la \log {\dt \over \la } \Big) - \sum_{d' = 1}^{D} (D - d' + 1) \sum_{i_1 = 1}^{n_1} \cdots \sum_{i_{d'} = 1}^{n_{d'}} \la _{( i_1 , \dots , i_{d'} )} \log {\pt _{( i_1 , \dots , i_{d'} )} \over \th _{( i_1 , \dots , i_{d'} )}} \Big| \X ( D' ) \Big] \non \\
&= E_{\pi ^{( \a )}} \Big[ (1 + D) \Big( \dt - \la - \la \log {\dt \over \la } \Big) \Big| \X ( D' ) \Big] \non \\
&\quad + \sum_{d = 1}^{D} (D - d + 1) \sum_{i_1 = 1}^{n_1} \cdots \sum_{i_{d - 1} = 1}^{n_{d - 1}} E_{\pi ^{( \a )}} \Big[ \Big( \la \prod_{d' = 1}^{d - 1} \th _{( i_1 , \dots , i_{d'} )} \Big) \Big| \X ( D' ) \Big] E_{\pi ^{( \a )}} \Big[ \sum_{i_d = 1}^{n_d} \th _{( i_1 , \dots , i_d )} \log {\th _{( i_1 , \dots , i_d )} \over \pt _{( i_1 , \dots , i_d )}} \Big| \X ( D' ) \Big] \text{,} \non 
\end{align}
where $\dt _{( i_1 , \dots , i_d )} = \sum_{i_{d + 1} = 1}^{n_{d + 1}} \cdots \sum_{i_D = 1}^{n_D} \dt _{i_1 , \dots , i_D}$ for $( i_1 , \dots , i_d ) \in \prod_{d' = 1}^{d} \{ 1, \dots , n_{d'} \} $ for $d = 0, 1, \dots , D$ and where $\pt _{( i_1 , \dots , i_d )} = \dt _{( i_1 , \dots , i_d )} / \dt _{( i_1 , \dots , i_{d - 1} )}$ for $( i_1 , \dots , i_d ) \in \prod_{d' = 1}^{d} \{ 1, \dots , n_{d'} \} $ for $d = 1, \dots , D$. 
Therefore, 
\begin{align}
&\min_{\dbt \in (0, \infty )^{n_1 \dotsm n_D}} E_{\pi ^{( \a )}} [ \widetilde{L} ^{\rm{KL}} ( \dbt , \bla ) | \X ( D' ) ] = E_{\pi ^{( \a )}} [ \widetilde{L} ^{\rm{KL}} ( \blah ^{( \a )} ( \X ( D' )), \bla ) | \X ( D' ) ] \text{,} \non 
\end{align}
where 
\begin{align}
\blah ^{( \a )} ( \X ( D' )) &= \Big( \cdots \Big( E_{\pi ^{( \a )}} [ \la | \X ( D' ) ] \prod_{d = 1}^{D} E_{\pi ^{( \a )}} [ \th _{( i_1 , \dots , i_d )} | \X ( D' ) ] \Big) _{i_D = 1, \dots , n_D} \cdots \Big) _{i_1 = 1, \dots , n_1} \non \\
&= \Big( \cdots \Big( {\Xt _{(+)} ( D' ) + a_0 \over 1 + D - D'} \prod_{d = 1}^{D} {\Xt _{( i_1 , \dots , i_d , +)} ( D' ) + a_d \over \sum_{{i_d}' = 1}^{n_d} \Xt _{( i_1 , \dots , i_{d - 1} , {i_d}' , +)} ( D' ) + n_d a_d} \Big) _{i_D = 1, \dots , n_D} \cdots \Big) _{i_1 = 1, \dots , n_1} \text{,} \non 
\end{align}
and this proves the lemma. 
\hfill$\Box$

\bigskip

\noindent
{\bf Proof of Theorem \ref{thm:general_1}.} \ \ Fix $D' = 1, \dots , D$ and let $\De _1 ( \bla ) = E_{\bla } [ \widetilde{L} ^{\rm{KL}} ( \blah ^{( \a ^{( D_0 )} ( D_0 ))} ( \X ( D' - 1)), \bla ) ] - E_{\bla } [ \widetilde{L} ^{\rm{KL}} ( \dbt , \blah ^{( \a ^{( D_0 )} ( D_0 ))} ( \X ( D' ))) ]$. 
Let $( a_{d}^{\rm{J}} )_{d = 0, 1, \dots , D} = \a ^{( D_0 )} ( D_0 )$. 
Then 
\begin{align}
\De _1 ( \bla ) &= E_{\bla } \Big[ \sum_{d = 0}^{D} \sum_{i_1 = 1}^{n_1} \cdots \sum_{i_d = 1}^{n_d} \Big( \sum_{i_{d + 1} = 1}^{n_{d + 1}} \cdots \sum_{i_D = 1}^{n_D} {\Xt _{(+)} ( D' - 1) + a_{0}^{\rm{J}} \over 2 + D - D'} \prod_{d' = 1}^{D} {\Xt _{( i_1 , \dots , i_{d'} , +)} ( D' - 1) + a_{d'}^{\rm{J}} \over \sum_{{i_{d'}}' = 1}^{n_{d'}} \Xt _{( i_1 , \dots , i_{d' - 1} , {i_{d'}}' , +)} ( D' - 1) + n_{d'} a_{d'}^{\rm{J}}} \non \\
&\quad - \sum_{i_{d + 1} = 1}^{n_{d + 1}} \cdots \sum_{i_D = 1}^{n_D} {\Xt _{(+)} ( D' ) + a_{0}^{\rm{J}} \over 1 + D - D'} \prod_{d' = 1}^{D} {\Xt _{( i_1 , \dots , i_{d'} , +)} ( D' ) + a_{d'}^{\rm{J}} \over \sum_{{i_{d'}}' = 1}^{n_{d'}} \Xt _{( i_1 , \dots , i_{d' - 1} , {i_{d'}}' , +)} ( D' ) + n_{d'} a_{d'}^{\rm{J}}} \non \\
&\quad - \la _{( i_1 , \dots , i_d )} \log \Big[ \Big\{ \sum_{i_{d + 1} = 1}^{n_{d + 1}} \cdots \sum_{i_D = 1}^{n_D} {\Xt _{(+)} ( D' - 1) + a_{0}^{\rm{J}} \over 2 + D - D'} \prod_{d' = 1}^{D} {\Xt _{( i_1 , \dots , i_{d'} , +)} ( D' - 1) + a_{d'}^{\rm{J}} \over \sum_{{i_{d'}}' = 1}^{n_{d'}} \Xt _{( i_1 , \dots , i_{d' - 1} , {i_{d'}}' , +)} ( D' - 1) + n_{d'} a_{d'}^{\rm{J}}} \Big\} \non \\
&\quad / \Big\{ \sum_{i_{d + 1} = 1}^{n_{d + 1}} \cdots \sum_{i_D = 1}^{n_D} {\Xt _{(+)} ( D' ) + a_{0}^{\rm{J}} \over 1 + D - D'} \prod_{d' = 1}^{D} {\Xt _{( i_1 , \dots , i_{d'} , +)} ( D' ) + a_{d'}^{\rm{J}} \over \sum_{{i_{d'}}' = 1}^{n_{d'}} \Xt _{( i_1 , \dots , i_{d' - 1} , {i_{d'}}' , +)} ( D' ) + n_{d'} a_{d'}^{\rm{J}}} \Big\} \Big] \Big) \Big] \non \\
&= E_{\bla } \Big[ (1 + D) \Big( {\Xt _{(+)} ( D' - 1) + a_{0}^{\rm{J}} \over 2 + D - D'} - {\Xt _{(+)} ( D' ) + a_{0}^{\rm{J}} \over 1 + D - D'} \Big) \non \\
&\quad + \sum_{d = 0}^{D} \sum_{i_1 = 1}^{n_1} \cdots \sum_{i_d = 1}^{n_d} \la _{( i_1 , \dots , i_d )} \log \Big[ \Big\{ {\Xt _{(+)} ( D' ) + a_{0}^{\rm{J}} \over 1 + D - D'} \prod_{d' = 1}^{d} {\Xt _{( i_1 , \dots , i_{d'} , +)} ( D' ) + a_{d'}^{\rm{J}} \over \sum_{{i_{d'}}' = 1}^{n_{d'}} \Xt _{( i_1 , \dots , i_{d' - 1} , {i_{d'}}' , +)} ( D' ) + n_{d'} a_{d'}^{\rm{J}}} \Big\} \non \\
&\quad / \Big\{ {\Xt _{(+)} ( D' - 1) + a_{0}^{\rm{J}} \over 2 + D - D'} \prod_{d' = 1}^{d} {\Xt _{( i_1 , \dots , i_{d'} , +)} ( D' - 1) + a_{d'}^{\rm{J}} \over \sum_{{i_{d'}}' = 1}^{n_{d'}} \Xt _{( i_1 , \dots , i_{d' - 1} , {i_{d'}}' , +)} ( D' - 1) + n_{d'} a_{d'}^{\rm{J}}} \Big\} \Big] \Big] \non \\
&= E_{\bla } \Big[ (1 + D) \Big( {a_{0}^{\rm{J}} \over 2 + D - D'} - {a_{0}^{\rm{J}} \over 1 + D - D'} \Big) + (1 + D) \la \log \Big\{ {\Xt _{(+)} ( D' ) + a_{0}^{\rm{J}} \over 1 + D - D'} / {\Xt _{(+)} ( D' - 1) + a_{0}^{\rm{J}} \over 2 + D - D'} \Big\} \non \\
&\quad + \sum_{d' = 1}^{D} \sum_{d = d'}^{D} \sum_{i_1 = 1}^{n_1} \cdots \sum_{i_{d'} = 1}^{n_{d'}} \la _{( i_1 , \dots , i_{d'} )} \log \Big\{ {\Xt _{( i_1 , \dots , i_{d'} , +)} ( D' ) + a_{d'}^{\rm{J}} \over \sum_{{i_{d'}}' = 1}^{n_{d'}} \Xt _{( i_1 , \dots , i_{d' - 1} , {i_{d'}}' , +)} ( D' ) + n_{d'} a_{d'}^{\rm{J}}} \non \\
&\quad / {\Xt _{( i_1 , \dots , i_{d'} , +)} ( D' - 1) + a_{d'}^{\rm{J}} \over \sum_{{i_{d'}}' = 1}^{n_{d'}} \Xt _{( i_1 , \dots , i_{d' - 1} , {i_{d'}}' , +)} ( D' - 1) + n_{d'} a_{d'}^{\rm{J}}} \Big\} \Big] \non \\
&= E_{\bla } \Big[ (1 + D) \Big( {a_{0}^{\rm{J}} \over 2 + D - D'} - {a_{0}^{\rm{J}} \over 1 + D - D'} \Big) + (1 + D) \la \log \Big\{ {\Xt _{(+)} ( D' ) + a_{0}^{\rm{J}} \over 1 + D - D'} / {\Xt _{(+)} ( D' - 1) + a_{0}^{\rm{J}} \over 2 + D - D'} \Big\} \non \\
&\quad + \sum_{d' = 1}^{D' - 1} (1 + D - d' ) \sum_{i_1 = 1}^{n_1} \cdots \sum_{i_{d'} = 1}^{n_{d'}} \la _{( i_1 , \dots , i_{d'} )} \log \Big\{ {\Xt _{( i_1 , \dots , i_{d'} , +)} ( D' ) + a_{d'}^{\rm{J}} \over \sum_{{i_{d'}}' = 1}^{n_{d'}} \Xt _{( i_1 , \dots , i_{d' - 1} , {i_{d'}}' , +)} ( D' ) + n_{d'} a_{d'}^{\rm{J}}} \non \\
&\quad / {\Xt _{( i_1 , \dots , i_{d'} , +)} ( D' - 1) + a_{d'}^{\rm{J}} \over \sum_{{i_{d'}}' = 1}^{n_{d'}} \Xt _{( i_1 , \dots , i_{d' - 1} , {i_{d'}}' , +)} ( D' - 1) + n_{d'} a_{d'}^{\rm{J}}} \Big\} \Big] \text{.} \non 
\end{align}
Since 
\begin{align}
&E_{\bla } \Big[ \sum_{d' = 1}^{D' - 1} (1 + D - d' ) \sum_{i_1 = 1}^{n_1} \cdots \sum_{i_{d'} = 1}^{n_{d'}} \la _{( i_1 , \dots , i_{d'} )} \log \Big\{ {\Xt _{( i_1 , \dots , i_{d'} , +)} ( D' ) + a_{d'}^{\rm{J}} \over \sum_{{i_{d'}}' = 1}^{n_{d'}} \Xt _{( i_1 , \dots , i_{d' - 1} , {i_{d'}}' , +)} ( D' ) + n_{d'} a_{d'}^{\rm{J}}} \non \\
&\quad / {\Xt _{( i_1 , \dots , i_{d'} , +)} ( D' - 1) + a_{d'}^{\rm{J}} \over \sum_{{i_{d'}}' = 1}^{n_{d'}} \Xt _{( i_1 , \dots , i_{d' - 1} , {i_{d'}}' , +)} ( D' - 1) + n_{d'} a_{d'}^{\rm{J}}} \Big\} \Big] \non \\
&= E_{\bla } \Big[ \sum_{d' = 1}^{D' - 1} (1 + D - d' ) \sum_{i_1 = 1}^{n_1} \cdots \sum_{i_{d'} = 1}^{n_{d'}} \la _{( i_1 , \dots , i_{d'} )} \log {\Xt _{( i_1 , \dots , i_{d'} , +)} ( D' ) + a_{d'}^{\rm{J}} \over \Xt _{( i_1 , \dots , i_{d'} , +)} ( D' - 1) + a_{d'}^{\rm{J}}} \non \\
&\quad - \sum_{d' = 0}^{D' - 2} (D - d' ) \sum_{i_1 = 1}^{n_1} \cdots \sum_{i_{d' + 1} = 1}^{n_{d' + 1}} \la _{( i_1 , \dots , i_{d' + 1} )} \log {\sum_{{i_{d' + 1}}' = 1}^{n_{d' + 1}} \Xt _{( i_1 , \dots , i_{d'} , {i_{d' + 1}}' , +)} ( D' ) + n_{d' + 1} a_{d' + 1}^{\rm{J}} \over \sum_{{i_{d' + 1}}' = 1}^{n_{d' + 1}} \Xt _{( i_1 , \dots , i_{d'} , {i_{d' + 1}}' , +)} ( D' - 1) + n_{d' + 1} a_{d' + 1}^{\rm{J}}} \Big] \non \\
&= E_{\bla } \Big[ (2 + D - D' ) \sum_{i_1 = 1}^{n_1} \cdots \sum_{i_{D' - 1} = 1}^{n_{D' - 1}} \la _{( i_1 , \dots , i_{D' - 1} )} \log {\Xt _{( i_1 , \dots , i_{D' - 1} , +)} ( D' ) + a_{D' - 1}^{\rm{J}} \over \Xt _{( i_1 , \dots , i_{D' - 1} , +)} ( D' - 1) + a_{D' - 1}^{\rm{J}}} \non \\
&\quad + \sum_{d' = 1}^{D' - 2} \sum_{i_1 = 1}^{n_1} \cdots \sum_{i_{d'} = 1}^{n_{d'}} \Big\{ (1 + D - d' ) \la _{( i_1 , \dots , i_{d'} )} \log {\Xt _{( i_1 , \dots , i_{d'} , +)} ( D' ) + a_{d'}^{\rm{J}} \over \Xt _{( i_1 , \dots , i_{d'} , +)} ( D' - 1) + a_{d'}^{\rm{J}}} \non \\
&\quad - (D - d' ) \la _{( i_1 , \dots , i_{d'} )} \log {\sum_{{i_{d' + 1}}' = 1}^{n_{d' + 1}} \Xt _{( i_1 , \dots , i_{d'} , {i_{d' + 1}}' , +)} ( D' ) + n_{d' + 1} a_{d' + 1}^{\rm{J}} \over \sum_{{i_{d' + 1}}' = 1}^{n_{d' + 1}} \Xt _{( i_1 , \dots , i_{d'} , {i_{d' + 1}}' , +)} ( D' - 1) + n_{d' + 1} a_{d' + 1}^{\rm{J}}} \Big\} \non \\
&\quad - D \la \log {\sum_{{i_1}' = 1}^{n_1} \Xt _{( {i_1}' , +)} ( D' ) + n_1 a_{1}^{\rm{J}} \over \sum_{{i_1}' = 1}^{n_1} \Xt _{( {i_1}' , +)} ( D' - 1) + n_1 a_{1}^{\rm{J}}} \Big] \non 
\end{align}
and since $n_{d' + 1} a_{d' + 1}^{\rm{J}} = a_{d'}^{\rm{J}}$ for all $d' = 0, \dots , D' - 2$, it follows that 
\begin{align}
\De _1 ( \bla ) &= E_{\bla } \Big[ (1 + D) \Big( {a_{0}^{\rm{J}} \over 2 + D - D'} - {a_{0}^{\rm{J}} \over 1 + D - D'} \Big) + (1 + D) \la \log {2 + D - D' \over 1 + D - D'} \non \\
&\quad + (2 + D - D' ) \sum_{i_1 = 1}^{n_1} \cdots \sum_{i_{D' - 1} = 1}^{n_{D' - 1}} \la _{( i_1 , \dots , i_{D' - 1} )} \log {\Xt _{( i_1 , \dots , i_{D' - 1} , +)} ( D' ) + a_{D' - 1}^{\rm{J}} \over \Xt _{( i_1 , \dots , i_{D' - 1} , +)} ( D' - 1) + a_{D' - 1}^{\rm{J}}} \non \\
&\quad + \sum_{d' = 0}^{D' - 2} \sum_{i_1 = 1}^{n_1} \cdots \sum_{i_{d'} = 1}^{n_{d'}} \Big\{ (1 + D - d' ) \la _{( i_1 , \dots , i_{d'} )} \log {\Xt _{( i_1 , \dots , i_{d'} , +)} ( D' ) + a_{d'}^{\rm{J}} \over \Xt _{( i_1 , \dots , i_{d'} , +)} ( D' - 1) + a_{d'}^{\rm{J}}} \non \\
&\quad - (D - d' ) \la _{( i_1 , \dots , i_{d'} )} \log {\sum_{{i_{d' + 1}}' = 1}^{n_{d' + 1}} \Xt _{( i_1 , \dots , i_{d'} , {i_{d' + 1}}' , +)} ( D' ) + a_{d'}^{\rm{J}} \over \sum_{{i_{d' + 1}}' = 1}^{n_{d' + 1}} \Xt _{( i_1 , \dots , i_{d'} , {i_{d' + 1}}' , +)} ( D' - 1) + a_{d'}^{\rm{J}}} \Big\} \Big] \non \\
&= E_{\bla } \Big[ (1 + D) \Big( {a_{0}^{\rm{J}} \over 2 + D - D'} - {a_{0}^{\rm{J}} \over 1 + D - D'} \Big) + (1 + D) \la \log {2 + D - D' \over 1 + D - D'} \non \\
&\quad + (2 + D - D' ) \sum_{i_1 = 1}^{n_1} \cdots \sum_{i_{D' - 1} = 1}^{n_{D' - 1}} \la _{( i_1 , \dots , i_{D' - 1} )} \log {\Xt _{( i_1 , \dots , i_{D' - 1} , +)} ( D' ) + a_{D' - 1}^{\rm{J}} \over \Xt _{( i_1 , \dots , i_{D' - 1} , +)} ( D' - 1) + a_{D' - 1}^{\rm{J}}} \non \\
&\quad + \sum_{d' = 0}^{D' - 2} \sum_{i_1 = 1}^{n_1} \cdots \sum_{i_{d'} = 1}^{n_{d'}} \la _{( i_1 , \dots , i_{d'} )} \log {\Xt _{( i_1 , \dots , i_{d'} , +)} ( D' ) + a_{d'}^{\rm{J}} \over \Xt _{( i_1 , \dots , i_{d'} , +)} ( D' - 1) + a_{d'}^{\rm{J}}} \Big] \non \\
&= \int_{1 + D - D'}^{2 + D - D'} \Big[ - {(1 + D) a_{0}^{\rm{J}} \over \de ^2} + (1 + D) \la {1 \over \de } \non \\
&\quad - (2 + D - D' ) \sum_{i_1 = 1}^{n_1} \cdots \sum_{i_{D' - 1} = 1}^{n_{D' - 1}} \la _{( i_1 , \dots , i_{D' - 1} )} {\pd \over \pd \de } E_{\bla } [ \log \{ \Zt _{( i_1 , \dots , i_{D' - 1} )} ( \de ) + a_{D' - 1}^{\rm{J}} \} ] \non \\
&\quad - \sum_{d' = 0}^{D' - 2} \sum_{i_1 = 1}^{n_1} \cdots \sum_{i_{d'} = 1}^{n_{d'}} \la _{( i_1 , \dots , i_{d'} )} {\pd \over \pd \de } E_{\bla } [ \log \{ \Zt _{( i_1 , \dots , i_{d'} )} ( \de ) + a_{d'}^{\rm{J}} \} ] \Big] d\de \text{,} \non 
\end{align}
where $\Zt _{( i_1 , \dots , i_{d'} )} ( \de ) \sim {\rm{Po}} ( \de \la _{( i_1 , \dots , i_{d'} )} )$ for $\de \in [1 + D - D' , 2 + D - D' ]$ for $( i_1 , \dots , i_{d'} ) \in \{ 1, \dots , n_1 \} \times \dots \times \{ 1, \dots , n_{d'} \} $ for $d' = 0, \dots , D' - 1$. 

Note that by Lemma \ref{lem:hudson}, 
\begin{align}
&{\pd \over \pd \de } E_{\bla } [ \log \{ \Zt _{( i_1 , \dots , i_{d'} )} ( \de ) + a_{d'}^{\rm{J}} \} ] \non \\
&= E_{\bla } \Big[ \Big\{ {\Zt _{( i_1 , \dots , i_{d'} )} ( \de ) \over \de } - \la _{( i_1 , \dots , i_{d'} )} \Big\} \log \{ \Zt _{( i_1 , \dots , i_{d'} )} ( \de ) + a_{d'}^{\rm{J}} \} \Big] \non \\
&= {1 \over \de } E_{\bla } \Big[ \Zt _{( i_1 , \dots , i_{d'} )} ( \de ) \log {\Zt _{( i_1 , \dots , i_{d'} )} ( \de ) + a_{d'}^{\rm{J}} \over \Zt _{( i_1 , \dots , i_{d'} )} ( \de ) + a_{d'}^{\rm{J}} - 1} \Big] \non 
\end{align}
for all $( i_1 , \dots , i_{d'} ) \in \{ 1, \dots , n_1 \} \times \dots \times \{ 1, \dots , n_{d'} \} $ for all $d' = 0, \dots , D' - 1$. 
Then 
\begin{align}
\De _1 ( \bla ) &= \int_{1 + D - D'}^{2 + D - D'} \Big\{ - {(1 + D) a_{0}^{\rm{J}} \over \de ^2} + (1 + D) \la {1 \over \de } \non \\
&\quad - (2 + D - D' ) \sum_{i_1 = 1}^{n_1} \cdots \sum_{i_{D' - 1} = 1}^{n_{D' - 1}} \la _{( i_1 , \dots , i_{D' - 1} )} {1 \over \de } E_{\bla } \Big[ \Zt _{( i_1 , \dots , i_{D' - 1} )} ( \de ) \log {\Zt _{( i_1 , \dots , i_{D' - 1} )} ( \de ) + a_{D' - 1}^{\rm{J}} \over \Zt _{( i_1 , \dots , i_{D' - 1} )} ( \de ) + a_{D' - 1}^{\rm{J}} - 1} \Big] \non \\
&\quad - \sum_{d' = 0}^{D' - 2} \sum_{i_1 = 1}^{n_1} \cdots \sum_{i_{d'} = 1}^{n_{d'}} \la _{( i_1 , \dots , i_{d'} )} {1 \over \de } E_{\bla } \Big[ \Zt _{( i_1 , \dots , i_{d'} )} ( \de ) \log {\Zt _{( i_1 , \dots , i_{d'} )} ( \de ) + a_{d'}^{\rm{J}} \over \Zt _{( i_1 , \dots , i_{d'} )} ( \de ) + a_{d'}^{\rm{J}} - 1} \Big\} \Big] d\de \non \\
&= \int_{1 + D - D'}^{2 + D - D'} \Big\{ - {(1 + D) a_{0}^{\rm{J}} \over \de ^2} \non \\
&\quad - (2 + D - D' ) \sum_{i_1 = 1}^{n_1} \cdots \sum_{i_{D' - 1} = 1}^{n_{D' - 1}} \la _{( i_1 , \dots , i_{D' - 1} )} {1 \over \de } E_{\bla } \Big[ \Zt _{( i_1 , \dots , i_{D' - 1} )} ( \de ) \log {\Zt _{( i_1 , \dots , i_{D' - 1} )} ( \de ) + a_{D' - 1}^{\rm{J}} \over \Zt _{( i_1 , \dots , i_{D' - 1} )} ( \de ) + a_{D' - 1}^{\rm{J}} - 1} - 1 \Big] \non \\
&\quad - \sum_{d' = 0}^{D' - 2} \sum_{i_1 = 1}^{n_1} \cdots \sum_{i_{d'} = 1}^{n_{d'}} \la _{( i_1 , \dots , i_{d'} )} {1 \over \de } E_{\bla } \Big[ \Zt _{( i_1 , \dots , i_{d'} )} ( \de ) \log {\Zt _{( i_1 , \dots , i_{d'} )} ( \de ) + a_{d'}^{\rm{J}} \over \Zt _{( i_1 , \dots , i_{d'} )} ( \de ) + a_{d'}^{\rm{J}} - 1} - 1 \Big] \Big\} d\de \non \\
&< \int_{1 + D - D'}^{2 + D - D'} \Big\{ - {(1 + D) a_{0}^{\rm{J}} \over \de ^2} \non \\
&\quad + (2 + D - D' ) \sum_{i_1 = 1}^{n_1} \cdots \sum_{i_{D' - 1} = 1}^{n_{D' - 1}} \la _{( i_1 , \dots , i_{D' - 1} )} {1 \over \de } E_{\bla } \Big[ {a_{D' - 1}^{\rm{J}} \over \Zt _{( i_1 , \dots , i_{D' - 1} )} ( \de ) + a_{D' - 1}^{\rm{J}}} \Big] \non \\
&\quad + \sum_{d' = 0}^{D' - 2} \sum_{i_1 = 1}^{n_1} \cdots \sum_{i_{d'} = 1}^{n_{d'}} \la _{( i_1 , \dots , i_{d'} )} {1 \over \de } E_{\bla } \Big[ {a_{d'}^{\rm{J}} \over \Zt _{( i_1 , \dots , i_{d'} )} ( \de ) + a_{d'}^{\rm{J}}} \Big] \Big\} d\de \text{.} \non 
\end{align}
Since $a_{d'}^{\rm{J}} \ge 1$ for all $d' = 0, \dots , D - 1$ by assumption, it follows from Lemma 2.1 of Fourdrinier, Strawderman and Wells (2018) that 
\begin{align}
\De _1 ( \bla ) &< \int_{1 + D - D'}^{2 + D - D'} \Big\{ - {(1 + D) a_{0}^{\rm{J}} \over \de ^2} \non \\
&\quad + (2 + D - D' ) \sum_{i_1 = 1}^{n_1} \cdots \sum_{i_{D' - 1} = 1}^{n_{D' - 1}} \la _{( i_1 , \dots , i_{D' - 1} )} {1 \over \de } E_{\bla } \Big[ {a_{D' - 1}^{\rm{J}} \over \de \la _{( i_1 , \dots , i_{D' - 1} )} + a_{D' - 1}^{\rm{J}} - 1} \Big] \non \\
&\quad + \sum_{d' = 0}^{D' - 2} \sum_{i_1 = 1}^{n_1} \cdots \sum_{i_{d'} = 1}^{n_{d'}} \la _{( i_1 , \dots , i_{d'} )} {1 \over \de } E_{\bla } \Big[ {a_{d'}^{\rm{J}} \over \de \la _{( i_1 , \dots , i_{d'} )} + a_{d'}^{\rm{J}} - 1} \Big] \Big\} d\de \non \\
&\le \int_{1 + D - D'}^{2 + D - D'} \Big\{ - {(1 + D) a_{0}^{\rm{J}} \over \de ^2} + (2 + D - D' ) \sum_{i_1 = 1}^{n_1} \cdots \sum_{i_{D' - 1} = 1}^{n_{D' - 1}} {a_{D' - 1}^{\rm{J}} \over \de ^2} + \sum_{d' = 0}^{D' - 2} \sum_{i_1 = 1}^{n_1} \cdots \sum_{i_{d'} = 1}^{n_{d'}} {a_{d'}^{\rm{J}} \over \de ^2} \Big\} d\de = 0 \text{,} \non 
\end{align}
and this completes the proof. 
\hfill$\Box$

\bigskip

\noindent
{\bf Proof of Theorem \ref{thm:general_2}.} \ \ Fix $D' = 1, \dots , D_0$. 
Let $\De _2 ( \bla ) = E_{\bla } [ \widetilde{L} ^{\rm{KL}} ( \blah ^{( \a ^{( D_0 )} ( D' - 1))} ( \X ), \bla ) ] - E_{\bla } [ \widetilde{L} ^{\rm{KL}} ( \blah ^{( \a ^{( D_0 )} ( D' ))} ( \X ), \bla ) ]$. 
Let $( a_d ( D' - 1))_{d = 0, 1, \dots , D} = \a ^{( D_0 )} ( D' - 1)$ and $( a_d ( D' ))_{d = 0, 1, \dots , D} = \a ^{( D_0 )} ( D' )$. 
Then 
\begin{align}
\De _2 ( \bla ) &= E_{\bla } \Big[ \sum_{d = 0}^{D} \sum_{i_1 = 1}^{n_1} \cdots \sum_{i_d = 1}^{n_d} \Big( \sum_{i_{d + 1} = 1}^{n_{d + 1}} \cdots \sum_{i_D = 1}^{n_D} {X_{(+)} + a_0 ( D' - 1) \over 1 + D} \prod_{d' = 1}^{D} {X_{( i_1 , \dots , i_{d'} , +)} + a_{d'} ( D' - 1) \over \sum_{{i_{d'}}' = 1}^{n_{d'}} X_{( i_1 , \dots , i_{{d'} - 1} , {i_{d'}}' , +)} + n_{d'} a_{d'} ( D' - 1)} \non \\
&\quad - \sum_{i_{d + 1} = 1}^{n_{d + 1}} \cdots \sum_{i_D = 1}^{n_D} {X_{(+)} + a_0 ( D' ) \over 1 + D} \prod_{d' = 1}^{D} {X_{( i_1 , \dots , i_{d'} , +)} + a_{d'} ( D' ) \over \sum_{{i_{d'}}' = 1}^{n_{d'}} X_{( i_1 , \dots , i_{{d'} - 1} , {i_{d'}}' , +)} + n_{d'} a_{d'} ( D' )} \non \\
&\quad - \la _{( i_1 , \dots , i_d )} \log \Big[ \Big\{ \sum_{i_{d + 1} = 1}^{n_{d + 1}} \cdots \sum_{i_D = 1}^{n_D} {X_{(+)} + a_0 ( D' - 1) \over 1 + D} \prod_{d' = 1}^{D} {X_{( i_1 , \dots , i_{d'} , +)} + a_{d'} ( D' - 1) \over \sum_{{i_{d'}}' = 1}^{n_{d'}} X_{( i_1 , \dots , i_{{d'} - 1} , {i_{d'}}' , +)} + n_{d'} a_{d'} ( D' - 1)} \Big\} \non \\
&\quad / \Big\{ \sum_{i_{d + 1} = 1}^{n_{d + 1}} \cdots \sum_{i_D = 1}^{n_D} {X_{(+)} + a_0 ( D' ) \over 1 + D} \prod_{d' = 1}^{D} {X_{( i_1 , \dots , i_{d'} , +)} + a_{d'} ( D' ) \over \sum_{{i_{d'}}' = 1}^{n_{d'}} X_{( i_1 , \dots , i_{{d'} - 1} , {i_{d'}}' , +)} + n_{d'} a_{d'} ( D' )} \Big\} \Big] \Big) \Big] \non \\
&= E_{\bla } \Big[ \sum_{d = 0}^{D} {X_{(+)} + a_0 ( D' - 1) \over 1 + D} - \sum_{d = 0}^{D} {X_{(+)} + a_0 ( D' ) \over 1 + D} + \sum_{d = 0}^{D} \la \log \Big\{ {X_{(+)} + a_0 ( D' ) \over 1 + D} / {X_{(+)} + a_0 ( D' - 1) \over 1 + D} \Big\} \non \\
&\quad + \sum_{d = 0}^{D} \sum_{i_1 = 1}^{n_1} \cdots \sum_{i_d = 1}^{n_d} \la _{( i_1 , \dots , i_d )} \log \Big[ \Big\{ \prod_{d' = 1}^{d} {X_{( i_1 , \dots , i_{d'} , +)} + a_{d'} ( D' ) \over \sum_{{i_{d'}}' = 1}^{n_{d'}} X_{( i_1 , \dots , i_{{d'} - 1} , {i_{d'}}' , +)} + n_{d'} a_{d'} ( D' )} \Big\} \non \\
&\quad / \Big\{ \prod_{d' = 1}^{d} {X_{( i_1 , \dots , i_{d'} , +)} + a_{d'} ( D' - 1) \over \sum_{{i_{d'}}' = 1}^{n_{d'}} X_{( i_1 , \dots , i_{{d'} - 1} , {i_{d'}}' , +)} + n_{d'} a_{d'} ( D' - 1)} \Big\} \Big] \Big] \non \\
&= E_{\bla } \Big[ - \{ a_0 ( D' ) - a_0 ( D' - 1) \} + (1 + D) \la \log {X_{(+)} + a_0 ( D' ) \over X_{(+)} + a_0 ( D' - 1)} \non \\
&\quad + \sum_{d' = 1}^{D} \sum_{d = d'}^{D} \sum_{i_1 = 1}^{n_1} \cdots \sum_{i_d = 1}^{n_d} \la _{( i_1 , \dots , i_d )} \log \Big\{ {X_{( i_1 , \dots , i_{d'} , +)} + a_{d'} ( D' ) \over \sum_{{i_{d'}}' = 1}^{n_{d'}} X_{( i_1 , \dots , i_{{d'} - 1} , {i_{d'}}' , +)} + n_{d'} a_{d'} ( D' )} \non \\
&\quad / {X_{( i_1 , \dots , i_{d'} , +)} + a_{d'} ( D' - 1) \over \sum_{{i_{d'}}' = 1}^{n_{d'}} X_{( i_1 , \dots , i_{{d'} - 1} , {i_{d'}}' , +)} + n_{d'} a_{d'} ( D' - 1)} \Big\} \Big] \text{.} \non 
\end{align}
Note that for all $d = 1, \dots , D$, we have $a_d ( D' - 1) = a_d ( D' )$ if $d \in [ D_0 , \infty ) \cup ((- \infty , D_0 - 1] \cap [ D' , \infty ))$. 
Then 
\begin{align}
\De _2 ( \bla ) &= E_{\bla } \Big[ - \{ a_0 ( D' ) - a_0 ( D' - 1) \} + (1 + D) \la \log {X_{(+)} + a_0 ( D' ) \over X_{(+)} + a_0 ( D' - 1)} \non \\
&\quad + \sum_{d' = 1}^{D' - 1} \sum_{d = d'}^{D} \sum_{i_1 = 1}^{n_1} \cdots \sum_{i_d = 1}^{n_d} \la _{( i_1 , \dots , i_d )} \log \Big\{ {X_{( i_1 , \dots , i_{d'} , +)} + a_{d'} ( D' ) \over \sum_{{i_{d'}}' = 1}^{n_{d'}} X_{( i_1 , \dots , i_{{d'} - 1} , {i_{d'}}' , +)} + n_{d'} a_{d'} ( D' )} \non \\
&\quad / {X_{( i_1 , \dots , i_{d'} , +)} + a_{d'} ( D' - 1) \over \sum_{{i_{d'}}' = 1}^{n_{d'}} X_{( i_1 , \dots , i_{{d'} - 1} , {i_{d'}}' , +)} + n_{d'} a_{d'} ( D' - 1)} \Big\} \Big] \non \\
&= E_{\bla } \Big[ - \{ a_0 ( D' ) - a_0 ( D' - 1) \} + (1 + D) \la \log {X_{(+)} + a_0 ( D' ) \over X_{(+)} + a_0 ( D' - 1)} \non \\
&\quad + \sum_{d' = 1}^{D' - 1} \sum_{d = d'}^{D} \sum_{i_1 = 1}^{n_1} \cdots \sum_{i_{d'} = 1}^{n_{d'}} \la _{( i_1 , \dots , i_{d'} )} \log \Big\{ {X_{( i_1 , \dots , i_{d'} , +)} + a_{d'} ( D' ) \over \sum_{{i_{d'}}' = 1}^{n_{d'}} X_{( i_1 , \dots , i_{{d'} - 1} , {i_{d'}}' , +)} + n_{d'} a_{d'} ( D' )} \non \\
&\quad / {X_{( i_1 , \dots , i_{d'} , +)} + a_{d'} ( D' - 1) \over \sum_{{i_{d'}}' = 1}^{n_{d'}} X_{( i_1 , \dots , i_{{d'} - 1} , {i_{d'}}' , +)} + n_{d'} a_{d'} ( D' - 1)} \Big\} \Big] \non \\
&= E_{\bla } \Big[ - \{ a_0 ( D' ) - a_0 ( D' - 1) \} + (1 + D) \la \log {X_{(+)} + a_0 ( D' ) \over X_{(+)} + a_0 ( D' - 1)} \non \\
&\quad + \sum_{d' = 1}^{D' - 1} (1 + D - d' ) \sum_{i_1 = 1}^{n_1} \cdots \sum_{i_{d'} = 1}^{n_{d'}} \la _{( i_1 , \dots , i_{d'} )} \log \Big\{ {X_{( i_1 , \dots , i_{d'} , +)} + a_{d'} ( D' ) \over X_{( i_1 , \dots , i_{d'} , +)} + a_{d'} ( D' - 1)} \non \\
&\quad / {\sum_{{i_{d'}}' = 1}^{n_{d'}} X_{( i_1 , \dots , i_{{d'} - 1} , {i_{d'}}' , +)} + n_{d'} a_{d'} ( D' ) \over \sum_{{i_{d'}}' = 1}^{n_{d'}} X_{( i_1 , \dots , i_{{d'} - 1} , {i_{d'}}' , +)} + n_{d'} a_{d'} ( D' - 1)} \Big\} \Big] \text{.} \non 
\end{align}
Therefore, 
\begin{align}
\De _2 ( \bla ) &= E_{\bla } \Big[ - \{ a_0 ( D' ) - a_0 ( D' - 1) \} + (1 + D) \la \log {X_{(+)} + a_0 ( D' ) \over X_{(+)} + a_0 ( D' - 1)} \non \\
&\quad + \sum_{d' = 1}^{D' - 1} (1 + D - d' ) \sum_{i_1 = 1}^{n_1} \cdots \sum_{i_{d'} = 1}^{n_{d'}} \la _{( i_1 , \dots , i_{d'} )} \log {X_{( i_1 , \dots , i_{d'} , +)} + a_{d'} ( D' ) \over X_{( i_1 , \dots , i_{d'} , +)} + a_{d'} ( D' - 1)} \non \\
&\quad - \sum_{d' = 0}^{D' - 2} (D - d' ) \sum_{i_1 = 1}^{n_1} \cdots \sum_{i_{d' + 1} = 1}^{n_{d' + 1}} \la _{( i_1 , \dots , i_{d' + 1} )} \log {\sum_{{i_{d' + 1}}' = 1}^{n_{d' + 1}} X_{( i_1 , \dots , i_{d'} , {i_{d' + 1}}' , +)} + n_{d' + 1} a_{d' + 1} ( D' ) \over \sum_{{i_{d' + 1}}' = 1}^{n_{d' + 1}} X_{( i_1 , \dots , i_{d'} , {i_{d' + 1}}' , +)} + n_{d' + 1} a_{d' + 1} ( D' - 1)} \Big] \non \\
&= E_{\bla } \Big[ - \{ a_0 ( D' ) - a_0 ( D' - 1) \} \non \\
&\quad + (2 + D - D' ) \sum_{i_1 = 1}^{n_1} \cdots \sum_{i_{D' - 1} = 1}^{n_{D' - 1}} \la _{( i_1 , \dots , i_{D' - 1} )} \log {X_{( i_1 , \dots , i_{D' - 1} , +)} + a_{D' - 1} ( D' ) \over X_{( i_1 , \dots , i_{D' - 1} , +)} + a_{D' - 1} ( D' - 1)} \non \\
&\quad + \sum_{d' = 0}^{D' - 2} (1 + D - d' ) \sum_{i_1 = 1}^{n_1} \cdots \sum_{i_{d'} = 1}^{n_{d'}} \la _{( i_1 , \dots , i_{d'} )} \log {X_{( i_1 , \dots , i_{d'} , +)} + a_{d'} ( D' ) \over X_{( i_1 , \dots , i_{d'} , +)} + a_{d'} ( D' - 1)} \non \\
&\quad - \sum_{d' = 0}^{D' - 2} (D - d' ) \sum_{i_1 = 1}^{n_1} \cdots \sum_{i_{d' + 1} = 1}^{n_{d' + 1}} \la _{( i_1 , \dots , i_{d' + 1} )} \log {\sum_{{i_{d' + 1}}' = 1}^{n_{d' + 1}} X_{( i_1 , \dots , i_{d'} , {i_{d' + 1}}' , +)} + n_{d' + 1} a_{d' + 1} ( D' ) \over \sum_{{i_{d' + 1}}' = 1}^{n_{d' + 1}} X_{( i_1 , \dots , i_{d'} , {i_{d' + 1}}' , +)} + n_{d' + 1} a_{d' + 1} ( D' - 1)} \Big] \text{.} \non 
\end{align}
Since $n_{d' + 1} a_{d' + 1} ( D' ) = a_{d'} ( D' )$ and $n_{d' + 1} a_{d' + 1} ( D' - 1) = a_{d'} ( D' - 1)$ for all $d' = 0, \dots , D' - 2$ and since $X_{( i_1 , \dots , i_{d'} , +)} \sim {\rm{Po}} ((1 + D - d' ) \la _{( i_1 , \dots , i_{d'} )} )$ and $\sum_{{i_{d' + 1}}' = 1}^{n_{d' + 1}} X_{( i_1 , \dots , i_{d'} , {i_{d' + 1}}' , +)} \sim {\rm{Po}} ((D - d' ) \la _{( i_1 , \dots , i_{d'} )} )$ for all $( i_1 , \dots , i_{d'} , i_{d' + 1} ) \in \{ 1, \dots , n_1 \} \times \dots \times \{ 1, \dots , n_{d'} \} \times \{ 1, \dots , n_{d' + 1} \} $ for all $d' = 0, \dots , D' - 2$, it follows that 
\begin{align}
\De _2 ( \bla ) &= E_{\bla } \Big[ - \{ a_0 ( D' ) - a_0 ( D' - 1) \} \non \\
&\quad + (2 + D - D' ) \sum_{i_1 = 1}^{n_1} \cdots \sum_{i_{D' - 1} = 1}^{n_{D' - 1}} \la _{( i_1 , \dots , i_{D' - 1} )} \log {X_{( i_1 , \dots , i_{D' - 1} , +)} + a_{D' - 1} ( D' ) \over X_{( i_1 , \dots , i_{D' - 1} , +)} + a_{D' - 1} ( D' - 1)} \Big] \non \\
&\quad + \sum_{d' = 0}^{D' - 2} \sum_{i_1 = 1}^{n_1} \cdots \sum_{i_{d'} = 1}^{n_{d'}} \int_{D - d'}^{1 + D - d'} {\pd \over \pd \de } E_{\bla } \Big[ \de \la _{( i_1 , \dots , i_{d'} )} \log {Z_{( i_1 , \dots , i_{d'} )} ( \de ) + a_{d'} ( D' ) \over Z_{( i_1 , \dots , i_{d'} )} ( \de ) + a_{d'} ( D' - 1)} \Big] d\de \text{,} \label{tgeneral_2p1} 
\end{align}
where $Z_{( i_1 , \dots , i_{d'} )} ( \de ) \sim {\rm{Po}} ( \de \la _{( i_1 , \dots , i_{d'} )} )$ for $\de \in [D - d' , 1 + D - d' ]$ for $( i_1 , \dots , i_{d'} ) \in \{ 1, \dots , n_1 \} \times \dots \times \{ 1, \dots , n_{d'} \} $ for $d' = 0, \dots , D' - 2$. 
If $D' = 1$, we obtain 
\begin{align}
\De _2 ( \bla ) &= E_{\bla } \Big[ - ( n_1 a_{D_0}^{( D_0 )} - a_{D_0}^{( D_0 )} ) + (1 + D) \la \log {X_{(+)} + n_1 a_{D_0}^{( D_0 )} \over X_{(+)} + a_{D_0}^{( D_0 )}} \Big] \non \\
&< E_{\bla } \Big[ - ( n_1 a_{D_0}^{( D_0 )} - a_{D_0}^{( D_0 )} ) + (1 + D) \la {n_1 a_{D_0}^{( D_0 )} - a_{D_0}^{( D_0 )} \over X_{(+)} + a_{D_0}^{( D_0 )}} \Big] \non \\
&\le - ( n_1 a_{D_0}^{( D_0 )} - a_{D_0}^{( D_0 )} ) + (1 + D) \la {n_1 a_{D_0}^{( D_0 )} - a_{D_0}^{( D_0 )} \over (1 + D) \la + a_{D_0}^{( D_0 )} - 1} < 0 \text{,} \non 
\end{align}
where the equality follows since 
\begin{align}
a_0 (1) = n_1 a_{D_0}^{( D_0 )} \quad \text{and} \quad a_0 (0) = a_{D_0}^{( D_0 )} \text{,} \non 
\end{align}
where the first and third inequalities follow since $n_1 \ge 2$ by assumption, and where the second inequality follows from Lemma 2.1 of Fourdrinier, Strawderman and Wells (2018). 

Now, suppose that $D' \ge 2$ and fix $d' = 0, \dots , D' - 2$ and $( i_1 , \dots , i_{d'} ) \in \{ 1, \dots , n_1 \} \times \dots \times \{ 1, \dots , n_{d'} \} $. 
Then for all $\de \in [D - d' , 1 + D - d' ]$, we have 
\begin{align}
&{\pd \over \pd \de } E_{\bla } \Big[ \de \la _{( i_1 , \dots , i_{d'} )} \log {Z_{( i_1 , \dots , i_{d'} )} ( \de ) + a_{d'} ( D' ) \over Z_{( i_1 , \dots , i_{d'} )} ( \de ) + a_{d'} ( D' - 1)} \Big] \non \\
&= \la _{( i_1 , \dots , i_{d'} )} E_{\bla } \Big[ \log {Z_{( i_1 , \dots , i_{d'} )} ( \de ) + a_{d'} ( D' ) \over Z_{( i_1 , \dots , i_{d'} )} ( \de ) + a_{d'} ( D' - 1)} \Big] \non \\
&\quad + \de \la _{( i_1 , \dots , i_{d'} )} E_{\bla } \Big[ \Big\{ {Z_{( i_1 , \dots , i_{d'} )} ( \de ) \over \de } - \la _{( i_1 , \dots , i_{d'} )} \Big\} \log {Z_{( i_1 , \dots , i_{d'} )} ( \de ) + a_{d'} ( D' ) \over Z_{( i_1 , \dots , i_{d'} )} ( \de ) + a_{d'} ( D' - 1)} \Big] \non \\
&= \la _{( i_1 , \dots , i_{d'} )} E_{\bla } \Big[ \log {Z_{( i_1 , \dots , i_{d'} )} ( \de ) + a_{d'} ( D' ) \over Z_{( i_1 , \dots , i_{d'} )} ( \de ) + a_{d'} ( D' - 1)} \non \\
&\quad + Z_{( i_1 , \dots , i_{d'} )} ( \de ) \log \Big\{ {Z_{( i_1 , \dots , i_{d'} )} ( \de ) + a_{d'} ( D' ) \over Z_{( i_1 , \dots , i_{d'} )} ( \de ) + a_{d'} ( D' - 1)} / {Z_{( i_1 , \dots , i_{d'} )} ( \de ) + a_{d'} ( D' ) - 1 \over Z_{( i_1 , \dots , i_{d'} )} ( \de ) + a_{d'} ( D' - 1) - 1} \Big\} \Big] \non \\
&= \la _{( i_1 , \dots , i_{d'} )} E_{\bla } \Big[ \log \Big\{ 1 + {a_{d'} ( D' ) - a_{d'} ( D' - 1) \over Z_{( i_1 , \dots , i_{d'} )} ( \de ) + a_{d'} ( D' - 1)} \Big\} \non \\
&\quad + Z_{( i_1 , \dots , i_{d'} )} ( \de ) \log \Big[ 1 - {a_{d'} ( D' ) - a_{d'} ( D' - 1) \over \{ Z_{( i_1 , \dots , i_{d'} )} ( \de ) + a_{d'} ( D' - 1) \} \{ Z_{( i_1 , \dots , i_{d'} )} ( \de ) + a_{d'} ( D' ) - 1 \} } \Big] \Big] \text{,} \non 
\end{align}
where the second equality follows from Lemma \ref{lem:hudson}. 
Note that $a_{d'} ( D' ) \ge a_{d'} ( D' - 1) \ge a_{D_0}^{( D_0 )} \ge 2$ by assumption. 
Then for all $\de \in [D - d' , 1 + D - d' ]$, 
\begin{align}
&{\pd \over \pd \de } E_{\bla } \Big[ \de \la _{( i_1 , \dots , i_{d'} )} \log {Z_{( i_1 , \dots , i_{d'} )} ( \de ) + a_{d'} ( D' ) \over Z_{( i_1 , \dots , i_{d'} )} ( \de ) + a_{d'} ( D' - 1)} \Big] \non \\
&\le \la _{( i_1 , \dots , i_{d'} )} E_{\bla } \Big[ {a_{d'} ( D' ) - a_{d'} ( D' - 1) \over Z_{( i_1 , \dots , i_{d'} )} ( \de ) + a_{d'} ( D' - 1)} - {\{ a_{d'} ( D' ) - a_{d'} ( D' - 1) \} Z_{( i_1 , \dots , i_{d'} )} ( \de ) \over \{ Z_{( i_1 , \dots , i_{d'} )} ( \de ) + a_{d'} ( D' - 1) \} \{ Z_{( i_1 , \dots , i_{d'} )} ( \de ) + a_{d'} ( D' ) - 1 \} } \Big] \non \\
&= \la _{( i_1 , \dots , i_{d'} )} E_{\bla } \Big[ {a_{d'} ( D' ) - a_{d'} ( D' - 1) \over Z_{( i_1 , \dots , i_{d'} )} ( \de ) + a_{d'} ( D' - 1)} {a_{d'} ( D' ) - 1 \over Z_{( i_1 , \dots , i_{d'} )} ( \de ) + a_{d'} ( D' ) - 1} \Big] \non \\
&= {1 \over \de } E_{\bla } \Big[ {Z_{( i_1 , \dots , i_{d'} )} ( \de ) \over Z_{( i_1 , \dots , i_{d'} )} ( \de ) + a_{d'} ( D' - 1) - 1} {\{ a_{d'} ( D' ) - a_{d'} ( D' - 1) \} \{ a_{d'} ( D' ) - 1 \} \over Z_{( i_1 , \dots , i_{d'} )} ( \de ) + a_{d'} ( D' ) - 2} \Big] \text{,} \non 
\end{align}
where the second equality follows from Lemma \ref{lem:hudson}, and this implies that 
\begin{align}
&{\pd \over \pd \de } E_{\bla } \Big[ \de \la _{( i_1 , \dots , i_{d'} )} \log {Z_{( i_1 , \dots , i_{d'} )} ( \de ) + a_{d'} ( D' ) \over Z_{( i_1 , \dots , i_{d'} )} ( \de ) + a_{d'} ( D' - 1)} \Big] \non \\
&\le {1 \over \de } E_{\bla } \Big[ {\{ a_{d'} ( D' ) - a_{d'} ( D' - 1) \} \{ a_{d'} ( D' ) - 1 \} \over Z_{( i_1 , \dots , i_{d'} )} ( \de ) + a_{d'} ( D' - 1) - 1} \Big] \non \\
&\le {1 \over \de } {\{ a_{d'} ( D' ) - a_{d'} ( D' - 1) \} \{ a_{d'} ( D' ) - 1 \} \over \de \la _{( i_1 , \dots , i_{d'} )} + a_{d'} ( D' - 1) - 2} \text{,} \non 
\end{align}
where the second inequality follows from Lemma 2.1 of Fourdrinier, Strawderman and Wells (2018). 
Therefore, 
\begin{align}
&\int_{D - d'}^{1 + D - d'} {\pd \over \pd \de } E_{\bla } \Big[ \de \la _{( i_1 , \dots , i_{d'} )} \log {Z_{( i_1 , \dots , i_{d'} )} ( \de ) + a_{d'} ( D' ) \over Z_{( i_1 , \dots , i_{d'} )} ( \de ) + a_{d'} ( D' - 1)} \Big] d\de \non \\
&\le \int_{D - d'}^{1 + D - d'} {1 \over \de } {\{ a_{d'} ( D' ) - a_{d'} ( D' - 1) \} \{ a_{d'} ( D' ) - 1 \} \over \de \la _{( i_1 , \dots , i_{d'} )} + a_{d'} ( D' - 1) - 2} d\de \non \\
&\le {1 \over D - d'} {\{ a_{d'} ( D' ) - a_{d'} ( D' - 1) \} \{ a_{d'} ( D' ) - 1 \} \over (D - d' ) \la _{( i_1 , \dots , i_{d'} )} + a_{d'} ( D' - 1) - 2} \text{.} \label{tgeneral_2p2} 
\end{align}

Meanwhile, since 
\begin{align}
a_0 ( D' ) &= \Big( \prod_{d' = 1}^{D'} n_{d'} \Big) a_{D_0}^{( D_0 )} \quad \text{and} \quad a_0 ( D' - 1) = \Big( \prod_{d' = 1}^{D' - 1} n_{d'} \Big) a_{D_0}^{( D_0 )} \non 
\end{align}
and 
\begin{align}
a_{D' - 1} ( D' ) &= n_{D'} a_{D_0}^{( D_0 )} \quad \text{and} \quad a_{D' - 1} ( D' - 1) = a_{D_0}^{( D_0 )} \text{,} \non 
\end{align}
we have 
\begin{align}
&E_{\bla } \Big[ - \{ a_0 ( D' ) - a_0 ( D' - 1) \} + (2 + D - D' ) \sum_{i_1 = 1}^{n_1} \cdots \sum_{i_{D' - 1} = 1}^{n_{D' - 1}} \la _{( i_1 , \dots , i_{D' - 1} )} \log {X_{( i_1 , \dots , i_{D' - 1} , +)} + a_{D' - 1} ( D' ) \over X_{( i_1 , \dots , i_{D' - 1} , +)} + a_{D' - 1} ( D' - 1)} \Big] \non \\
&= E_{\bla } \Big[ - ( n_{D'} - 1) \Big( \prod_{d' = 1}^{D' - 1} n_{d'} \Big) a_{D_0}^{( D_0 )} + (2 + D - D' ) \sum_{i_1 = 1}^{n_1} \cdots \sum_{i_{D' - 1} = 1}^{n_{D' - 1}} \la _{( i_1 , \dots , i_{D' - 1} )} \log {X_{( i_1 , \dots , i_{D' - 1} , +)} + n_{D'} a_{D_0}^{( D_0 )} \over X_{( i_1 , \dots , i_{D' - 1} , +)} + a_{D_0}^{( D_0 )}} \Big] \non \\
&< E_{\bla } \Big[ - ( n_{D'} - 1) \Big( \prod_{d' = 1}^{D' - 1} n_{d'} \Big) a_{D_0}^{( D_0 )} + (2 + D - D' ) \sum_{i_1 = 1}^{n_1} \cdots \sum_{i_{D' - 1} = 1}^{n_{D' - 1}} \la _{( i_1 , \dots , i_{D' - 1} )} {n_{D'} a_{D_0}^{( D_0 )} - a_{D_0}^{( D_0 )} \over X_{( i_1 , \dots , i_{D' - 1} , +)} + a_{D_0}^{( D_0 )}} \Big] \non \\
&\le - ( n_{D'} - 1) \Big( \prod_{d' = 1}^{D' - 1} n_{d'} \Big) a_{D_0}^{( D_0 )} + (2 + D - D' ) \sum_{i_1 = 1}^{n_1} \cdots \sum_{i_{D' - 1} = 1}^{n_{D' - 1}} \la _{( i_1 , \dots , i_{D' - 1} )} {n_{D'} a_{D_0}^{( D_0 )} - a_{D_0}^{( D_0 )} \over (2 + D - D' ) \la _{( i_1 , \dots , i_{D' - 1} )} + a_{D_0}^{( D_0 )} - 1} \non \\
&= - \sum_{i_1 = 1}^{n_1} \cdots \sum_{i_{D' - 1} = 1}^{n_{D' - 1}} {( n_{D'} - 1) a_{D_0}^{( D_0 )} ( a_{D_0}^{( D_0 )} - 1) \over (2 + D - D' ) \la _{( i_1 , \dots , i_{D' - 1} )} + a_{D_0}^{( D_0 )} - 1} \text{,} \non 
\end{align}
where the first inequality follows since $n_{D'} \ge 2$ by assumption and where the second inequality follows from Lemma 2.1 of Fourdrinier, Strawderman and Wells (2018). 
Thus, by Jensen's inequality, 
\begin{align}
&E_{\bla } \Big[ - \{ a_0 ( D' ) - a_0 ( D' - 1) \} + (2 + D - D' ) \sum_{i_1 = 1}^{n_1} \cdots \sum_{i_{D' - 1} = 1}^{n_{D' - 1}} \la _{( i_1 , \dots , i_{D' - 1} )} \log {X_{( i_1 , \dots , i_{D' - 1} , +)} + a_{D' - 1} ( D' ) \over X_{( i_1 , \dots , i_{D' - 1} , +)} + a_{D' - 1} ( D' - 1)} \Big] \non \\
&< - \sum_{i_1 = 1}^{n_1} \cdots \sum_{i_{d'} = 1}^{n_{d'}} {n_{d' + 1} \dotsm n_{D' - 1} \over n_{d' + 1} \dotsm n_{D' - 1}} \sum_{i_{d' + 1} = 1}^{n_{d' + 1}} \cdots \sum_{i_{D' - 1} = 1}^{n_{D' - 1}} {( n_{D'} - 1) a_{D_0}^{( D_0 )} ( a_{D_0}^{( D_0 )} - 1) \over (2 + D - D' ) \la _{( i_1 , \dots , i_{D' - 1} )} + a_{D_0}^{( D_0 )} - 1} \non \\
&\le - \sum_{i_1 = 1}^{n_1} \cdots \sum_{i_{d'} = 1}^{n_{d'}} {( n_{d' + 1} \dotsm n_{D' - 1} ) ( n_{D'} - 1) a_{D_0}^{( D_0 )} ( a_{D_0}^{( D_0 )} - 1) \over \{ (2 + D - D' ) / ( n_{d' + 1} \dotsm n_{D' - 1} ) \} \sum_{i_{d' + 1} = 1}^{n_{d' + 1}} \cdots \sum_{i_{D' - 1} = 1}^{n_{D' - 1}} \la _{( i_1 , \dots , i_{D' - 1} )} + a_{D_0}^{( D_0 )} - 1} \non \\
&= - \sum_{i_1 = 1}^{n_1} \cdots \sum_{i_{d'} = 1}^{n_{d'}} {( n_{d' + 1} \dotsm n_{D' - 1} )^2 ( n_{D'} - 1) a_{D_0}^{( D_0 )} ( a_{D_0}^{( D_0 )} - 1) \over (2 + D - D' ) \la _{( i_1 , \dots , i_{d'} )} + ( n_{d' + 1} \dotsm n_{D' - 1} ) ( a_{D_0}^{( D_0 )} - 1)} \non \\
&= - \sum_{i_1 = 1}^{n_1} \cdots \sum_{i_{d'} = 1}^{n_{d'}} {\{ a_{d'} ( D' ) - a_{d'} ( D' - 1) \} a_{d'} ( D' - 1) ( a_{D_0}^{( D_0 )} - 1) / a_{D_0}^{( D_0 )} \over (2 + D - D' ) \la _{( i_1 , \dots , i_{d'} )} + a_{d'} ( D' - 1) ( a_{D_0}^{( D_0 )} - 1) / a_{D_0}^{( D_0 )}} \label{tgeneral_2p3} 
\end{align}
for all $d' = 0, \dots , D' - 2$, where the second equality follows since 
\begin{align}
a_{d'} ( D' ) &= ( n_{d' + 1} \dotsm n_{D'} ) a_{D_0}^{( D_0 )} \quad \text{and} \quad a_{d'} ( D' - 1) = ( n_{d' + 1} \dotsm n_{D' - 1} ) a_{D_0}^{( D_0 )} \text{.} \non 
\end{align}

Finally, if $D' \ge 2$, then combining (\ref{tgeneral_2p1}), (\ref{tgeneral_2p2}), and (\ref{tgeneral_2p3}), we obtain 
\begin{align}
\De _2 ( \bla ) &< - \sum_{d' = 0}^{D' - 2} \sum_{i_1 = 1}^{n_1} \cdots \sum_{i_{d'} = 1}^{n_{d'}} \Big[ {1 \over D' - 1} {\{ a_{d'} ( D' ) - a_{d'} ( D' - 1) \} a_{d'} ( D' - 1) ( a_{D_0}^{( D_0 )} - 1) / a_{D_0}^{( D_0 )} \over (2 + D - D' ) \la _{( i_1 , \dots , i_{d'} )} + a_{d'} ( D' - 1) ( a_{D_0}^{( D_0 )} - 1) / a_{D_0}^{( D_0 )}} \non \\
&\quad - {1 \over D - d'} {\{ a_{d'} ( D' ) - a_{d'} ( D' - 1) \} \{ a_{d'} ( D' ) - 1 \} \over (D - d' ) \la _{( i_1 , \dots , i_{d'} )} + a_{d'} ( D' - 1) - 2} \Big] \text{.} \non 
\end{align}
The right-hand side of the above inequality is nonpositive since 
\begin{align}
&{D - d' \over D' - 1} {a_{d'} ( D' - 1) ( a_{D_0}^{( D_0 )} - 1) / a_{D_0}^{( D_0 )} \over (2 + D - D' ) \la _{( i_1 , \dots , i_{d'} )} + a_{d'} ( D' - 1) ( a_{D_0}^{( D_0 )} - 1) / a_{D_0}^{( D_0 )}} - {a_{d'} ( D' ) - 1 \over (D - d' ) \la _{( i_1 , \dots , i_{d'} )} + a_{d'} ( D' - 1) - 2} \non \\
&\ge {D - d' \over D' - 1} {a_{d'} ( D' - 1) ( a_{D_0}^{( D_0 )} - 1) / a_{D_0}^{( D_0 )} \over (D - d' ) \la _{( i_1 , \dots , i_{d'} )} + a_{d'} ( D' - 1)} - {n_{D'} a_{d'} ( D' - 1) \over (D - d' ) \la _{( i_1 , \dots , i_{d'} )} + a_{d'} ( D' - 1)} {a_{d'} ( D' - 1) \over a_{d'} ( D' - 1) - 2} \non \\
&\ge {a_{d'} ( D' - 1) \over (D - d' ) \la _{( i_1 , \dots , i_{d'} )} + a_{d'} ( D' - 1)} \Big( {2 + D - D' \over D' - 1} {a_{D_0}^{( D_0 )} - 1 \over a_{D_0}^{( D_0 )}} - n_{D'} {n_{D' - 1} a_{D_0}^{( D_0 )} \over n_{D' - 1} a_{D_0}^{( D_0 )} - 2} \Big) \ge 0 \non 
\end{align}
for all $( i_1 , \dots , i_{d'} ) \in \{ 1, \dots , n_1 \} \times \dots \times \{ 1, \dots , n_{d'} \} $ for all $d' = 0, \dots , D' - 2$ by assumption. 
This completes the proof. 
\hfill$\Box$

\subsection{Additional results}
\label{subsec:additional_results} 
Parts (i) and (ii) of the following proposition show that the risk functions, under the Kullback-Leibler divergence, of arbitrary plug-in and Bayesian predictive density estimators can be expressed using the risk functions, under the corresponding entropy loss, of the corresponding point estimators. 
These results are direct extensions of well-known results in simpler settings; see, for example, Robert (1996) and Komaki (2006, 2015). 

\begin{prp}
\label{prp:connection} 
Let $D$, $n_1 , \dots , n_D$, and $\bla $  be as in Section \ref{sec:general}. 
Let $\la _{( i_1 , \dots , i_d )}$ and $\th _{( i_1 , \dots , i_d )}$ be as in Section \ref{sec:general} for $( i_1 , \dots , i_d ) \in \prod_{d' = 1}^{d} \{ 1, \dots , n_{d'} \} $ for $d = 0, 1, \dots , D$. 
Let 
\begin{align}
(( \cdots ( r_{( i_1 , \dots , i_d )} ) _{i_d = 1, \dots , n_d} \cdots ) _{i_1 = 1, \dots , n_1} ) _{d = 0, 1, \dots , D} \in [0, \infty )^{\sum_{d = 0}^{D} \prod_{d' = 1}^{d} n_{d'}} \non 
\end{align}
and 
\begin{align}
(( \cdots ( s_{( i_1 , \dots , i_d )} ) _{i_d = 1, \dots , n_d} \cdots ) _{i_1 = 1, \dots , n_1} ) _{d = 0, 1, \dots , D} \in [0, \infty )^{\sum_{d = 0}^{D} \prod_{d' = 1}^{d} n_{d'}} \non 
\end{align}
be known constants. 
Suppose that 
\begin{align}
\X = (( \cdots ( X_{( i_1 , \dots , i_d )} )_{i_d = 1, \dots , n_d} \cdots )_{i_1 = 1, \dots , n_1} )_{d = 0, 1, \dots , D} \non 
\end{align}
and 
\begin{align}
\Y = (( \cdots ( Y_{( i_1 , \dots , i_d )} )_{i_d = 1, \dots , n_d} \cdots )_{i_1 = 1, \dots , n_1} )_{d = 0, 1, \dots , D} \non 
\end{align}
are independent Poisson variables with means 
\begin{align}
(( \cdots ( r_{( i_1 , \dots , i_d )} \la _{( i_1 , \dots , i_d )} ) _{i_d = 1, \dots , n_d} \cdots ) _{i_1 = 1, \dots , n_1} ) _{d = 0, 1, \dots , D} \non 
\end{align}
and 
\begin{align}
(( \cdots ( s_{( i_1 , \dots , i_d )} \la _{( i_1 , \dots , i_d )} ) _{i_d = 1, \dots , n_d} \cdots ) _{i_1 = 1, \dots , n_1} ) _{d = 0, 1, \dots , D} \non 
\end{align}
and with likelihood functions $p( \X | \bla )$ and $p( \Y | \bla )$. 
\begin{itemize}
\item[{\rm{(i)}}]
Let $\blah ( \X ) = ( \cdots ( \lah _{i_1 , \dots , i_d} ( \X ))_{i_D = 1, \dots , n_D} \cdots )_{i_1 = 1, \dots , n_1} \in (0, \infty )^{n_1 \dotsm n_D}$ be an estimator of $\bla $ based on $\X$. 
Then the expected Kullback--Leibler divergence $E_{\bla } [ \log \{ p( \Y | \bla ) / p( \Y | \blah ( \X )) \} ]$ is given by 
\begin{align}
&E_{\bla } [ \log \{ p( \Y | \bla ) / p( \Y | \blah ( \X )) \} ] \non \\
&= E_{\bla } \Big[ \sum_{d = 0}^{D} \sum_{i_1 = 1}^{n_1} \cdots \sum_{i_d = 1}^{n_d} s_{( i_1 , \dots , i_d )} \la _{( i_1 , \dots , i_d )} \Big\{ {\lah _{( i_1 , \dots , i_d )} ( \X ) \over \la _{( i_1 , \dots , i_d )}} - 1 - \log {\lah _{( i_1 , \dots , i_d )} ( \X ) \over \la _{( i_1 , \dots , i_d )}} \Big\} \Big] \text{,} \non 
\end{align}
where $\lah _{( i_1 , \dots , i_d )} ( \X ) = \sum_{i_{d + 1} = 1}^{n_{d + 1}} \cdots \sum_{i_D = 1}^{n_D} \lah _{i_1 , \dots , i_D} ( \X )$ for $( i_1 , \dots , i_d ) \in \prod_{d' = 1}^{d} \{ 1, \dots , n_{d'} \} $ for $d = 0, 1, \dots , D$. 
\item[{\rm{(ii)}}]
Let $\pi ( \bla ) d\bla $ be a prior for $\bla $. 
Let $\ph ^{( \pi )} ( \y ; \X ) = E_{\pi } [ p( \y | \bla ) | \X ]$, $\y \in {\mathbb{N} _0}^{\sum_{d = 0}^{D} \prod_{d' = 1}^{d} n_{d'}}$, be the predictive density estimator for $\Y $ based on $\pi $ and $\X $. 
Then the expected Kullback--Leibler divergence $E_{\bla } [ \log \{ p( \Y | \bla ) / \ph ^{( \pi )} ( \Y ; \X ) \} ]$ is given by 
\begin{align}
&E_{\bla } [ \log \{ p( \Y | \bla ) / \ph ^{( \pi )} ( \Y ; \X ) \} ] \non \\
&= \int_{0}^{1} E_{\bla } \Big[ \sum_{d = 0}^{D} \sum_{i_1 = 1}^{n_1} \cdots \sum_{i_d = 1}^{n_d} {\pd t_{( i_1 , \dots , i_d )} \over \pd \ta } ( \ta ) \la _{( i_1 , \dots , i_d )} \Big\{ {\lah _{( i_1 , \dots , i_d )} ( \Z ( \ta )) \over \la _{( i_1 , \dots , i_d )}} - 1 - \log {\lah _{( i_1 , \dots , i_d )} ( \Z ( \ta )) \over \la _{( i_1 , \dots , i_d )}} \Big\} \Big] d\ta \text{,} \non 
\end{align}
where $t_{( i_1 , \dots , i_d )} \colon [0, 1] \to [0, \infty )$ is a nondecreasing function satisfying $t_{( i_1 , \dots , i_d )} (0) = r_{( i_1 , \dots , i_d )}$ and $t_{( i_1 , \dots , i_d )} (1) = s_{( i_1 , \dots , i_d )} + r_{( i_1 , \dots , i_d )}$ for $( i_1 , \dots , i_d ) \in \prod_{d' = 1}^{d} \{ 1, \dots , n_{d'} \} $ for $d = 0, 1, \dots , D$, where $\Z ( \ta ) \sim \prod_{d = 0}^{D} \prod_{i_1 = 1}^{n_1} \dotsc \prod_{i_d = 1}^{n_d} {\rm{Po}} ( t_{( i_1 , \dots , i_d )} ( \ta ) \la _{( i_1 , \dots , i_d )} )$ for $\ta \in [0, 1]$, and where $\lah _{( i_1 , \dots , i_d )} ( \Z ( \ta )) = E_{\pi } [ \la _{( i_1 , \dots , i_d )} | \Z ( \ta ) ]$ for $\ta \in [0, 1]$ for $( i_1 , \dots , i_d ) \in \prod_{d' = 1}^{d} \{ 1, \dots , n_{d'} \} $ for $d = 0, 1, \dots , D$. 
\end{itemize}
\end{prp}

\noindent
{\bf Proof%
.} \ \ For part (i), we have 
\begin{align}
&E_{\bla } [ \log \{ p( \Y | \bla ) / p( \Y | \blah ( \X )) \} ] \non \\
&= E_{\bla } \Big[ \sum_{d = 0}^{D} \sum_{i_1 = 1}^{n_1} \cdots \sum_{i_d = 1}^{n_d} \Big[ Y_{( i_1 , \dots , i_d )} \log {\la _{( i_1 , \dots , i_d )} \over \lah _{( i_1 , \dots , i_d )} ( \X )} - s_{( i_1 , \dots , i_d )} \{ \la _{( i_1 , \dots , i_d )} - \lah _{( i_1 , \dots , i_d )} ( \X ) \} \Big] \Big] \non \\
&= E_{\bla } \Big[ \sum_{d = 0}^{D} \sum_{i_1 = 1}^{n_1} \cdots \sum_{i_d = 1}^{n_d} s_{( i_1 , \dots , i_d )} \la _{( i_1 , \dots , i_d )} \Big\{ {\lah _{( i_1 , \dots , i_d )} ( \X ) \over \la _{( i_1 , \dots , i_d )}} - 1 - \log {\lah _{( i_1 , \dots , i_d )} ( \X ) \over \la _{( i_1 , \dots , i_d )}} \Big\} \Big] \text{.} \non 
\end{align}
For part (ii), note that 
\begin{align}
&\ph ^{( \pi )} ( \Y ; \X ) \non \\
&= \Big( \prod_{d = 0}^{D} \prod_{i_1 = 1}^{n_1} \cdots \prod_{i_d = 1}^{n_d} {{s_{( i_1 , \dots , i_d )}}^{Y_{( i_1 , \dots , i_d )}} \over Y_{( i_1 , \dots , i_d )} !} \Big) E_{\pi } \Big[ \prod_{d = 0}^{D} \prod_{i_1 = 1}^{n_1} \cdots \prod_{i_d = 1}^{n_d} {\la _{( i_1 , \dots , i_d )}}^{Y_{( i_1 , \dots , i_d )}} e^{- s_{( i_1 , \dots , i_d )} \la _{( i_1 , \dots , i_d )}} \Big| \X \Big] \non \\
&= \Big( \prod_{d = 0}^{D} \prod_{i_1 = 1}^{n_1} \cdots \prod_{i_d = 1}^{n_d} {{s_{( i_1 , \dots , i_d )}}^{Y_{( i_1 , \dots , i_d )}} \over Y_{( i_1 , \dots , i_d )} !} \Big) \non \\
&\quad \times \frac{ \int_{(0, \infty )^{n_1 \dotsm n_D}} \big[ \prod_{d = 0}^{D} \prod_{i_1 = 1}^{n_1} \cdots \prod_{i_d = 1}^{n_d} \{ {\la _{( i_1 , \dots , i_d )}}^{Y_{( i_1 , \dots , i_d )} + X_{( i_1 , \dots , i_d )}} e^{- ( s_{( i_1 , \dots , i_d )} + r_{( i_1 , \dots , i_d )} ) \la _{( i_1 , \dots , i_d )}} \} \big] \pi ( \bla ) d\bla }{ \int_{(0, \infty )^{n_1 \dotsm n_D}} \big[ \prod_{d = 0}^{D} \prod_{i_1 = 1}^{n_1} \cdots \prod_{i_d = 1}^{n_d} \{ {\la _{( i_1 , \dots , i_d )}}^{X_{( i_1 , \dots , i_d )}} e^{- r_{( i_1 , \dots , i_d )} \la _{( i_1 , \dots , i_d )}} \} \big] \pi ( \bla ) d\bla } \text{.} \non 
\end{align}
Then 
\begin{align}
&E_{\bla } [ \log \{ p( \Y | \bla ) / \ph ^{( \pi )} ( \Y ; \X ) \} ] \non \\
&= E_{\bla } \Big[ \log \prod_{d = 0}^{D} \prod_{i_1 = 1}^{n_1} \cdots \prod_{i_d = 1}^{n_d} \{ {\la _{( i_1 , \dots , i_d )}}^{Y_{( i_1 , \dots , i_d )}} e^{- s_{( i_1 , \dots , i_d )} \la _{( i_1 , \dots , i_d )}} \} \non \\
&\quad - \log \frac{ \int_{(0, \infty )^{n_1 \dotsm n_D}} \big[ \prod_{d = 0}^{D} \prod_{i_1 = 1}^{n_1} \cdots \prod_{i_d = 1}^{n_d} \{ {\la _{( i_1 , \dots , i_d )}}^{Y_{( i_1 , \dots , i_d )} + X_{( i_1 , \dots , i_d )}} e^{- ( s_{( i_1 , \dots , i_d )} + r_{( i_1 , \dots , i_d )} ) \la _{( i_1 , \dots , i_d )}} \} \big] \pi ( \bla ) d\bla }{ \int_{(0, \infty )^{n_1 \dotsm n_D}} \big[ \prod_{d = 0}^{D} \prod_{i_1 = 1}^{n_1} \cdots \prod_{i_d = 1}^{n_d} \{ {\la _{( i_1 , \dots , i_d )}}^{X_{( i_1 , \dots , i_d )}} e^{- r_{( i_1 , \dots , i_d )} \la _{( i_1 , \dots , i_d )}} \} \big] \pi ( \bla ) d\bla } \Big] \non \\
&= \int_{0}^{1} \Big( \sum_{d = 0}^{D} \sum_{i_1 = 1}^{n_1} \cdots \sum_{i_d = 1}^{n_d} {\pd t_{( i_1 , \dots , i_d )} \over \pd \ta } ( \ta ) ( \la _{( i_1 , \dots , i_d )} \log \la _{( i_1 , \dots , i_d )} - \la _{( i_1 , \dots , i_d )} ) \non \\
&\quad - {\pd \over \pd \ta } E_{\bla } \Big[ \log \int_{(0, \infty )^{n_1 \dotsm n_D}} \Big[ \prod_{d = 0}^{D} \prod_{i_1 = 1}^{n_1} \cdots \prod_{i_d = 1}^{n_d} \{ {\la _{( i_1 , \dots , i_d )}}^{Z_{( i_1 , \dots , i_d )} ( \ta )} e^{- t_{( i_1 , \dots , i_d )} ( \ta ) \la _{( i_1 , \dots , i_d )}} \} \Big] \pi ( \bla ) d\bla \Big] \Big) d\ta \text{,} \non 
\end{align}
where $(( \cdots ( Z_{( i_1 , \dots , i_d )} ( \ta ))_{i_d = 1, \dots , n_d} \cdots )_{i_1 = 1, \dots , n_1} )_{d = 1, \dots , D} = \Z ( \ta )$ for $\ta \in [0, 1]$. 
Note that for all $\ta \in [0, 1]$, 
\begin{align}
&{\pd \over \pd \ta } E_{\bla } \Big[ \log \int_{(0, \infty )^{n_1 \dotsm n_D}} \Big[ \prod_{d = 0}^{D} \prod_{i_1 = 1}^{n_1} \cdots \prod_{i_d = 1}^{n_d} \{ {\la _{( i_1 , \dots , i_d )}}^{Z_{( i_1 , \dots , i_d )} ( \ta )} e^{- t_{( i_1 , \dots , i_d )} ( \ta ) \la _{( i_1 , \dots , i_d )}} \} \Big] \pi ( \bla ) d\bla \Big] \non \\
&= E_{\bla } \Big[ \sum_{d = 0}^{D} \sum_{i_1 = 1}^{n_1} \cdots \sum_{i_d = 1}^{n_d} \Big[ {\pd t_{( i_1 , \dots , i_d )} \over \pd \ta } ( \ta ) \Big\{ {Z_{( i_1 , \dots , i_d )} ( \ta ) \over t_{( i_1 , \dots , i_d )} ( \ta )} - \la _{( i_1 , \dots , i_d )} \Big\} \non \\
&\quad \times \log \int_{(0, \infty )^{n_1 \dotsm n_D}} \Big[ \prod_{d' = 0}^{D} \prod_{{i_1}' = 1}^{n_1} \cdots \prod_{{i_d}' = 1}^{n_d} \{ {\la _{( {i_1}' , \dots , {i_d}' )}}^{Z_{( {i_1}' , \dots , {i_d}' )} ( \ta )} e^{- t_{( {i_1}' , \dots , {i_d}' )} ( \ta ) \la _{( {i_1}' , \dots , {i_d}' )}} \} \Big] \pi ( \bla ) d\bla \Big] \non \\
&\quad + \frac{ 1 }{ \int_{(0, \infty )^{n_1 \dotsm n_D}} \big[ \prod_{d' = 0}^{D} \prod_{{i_1}' = 1}^{n_1} \cdots \prod_{{i_d}' = 1}^{n_d} \{ {\la _{( {i_1}' , \dots , {i_d}' )}}^{Z_{( {i_1}' , \dots , {i_d}' )} ( \ta )} e^{- t_{( {i_1}' , \dots , {i_d}' )} ( \ta ) \la _{( {i_1}' , \dots , {i_d}' )}} \} \big] \pi ( \bla ) d\bla } \non \\
&\quad \times \int_{(0, \infty )^{n_1 \dotsm n_D}} \Big( \sum_{d = 0}^{D} \sum_{i_1 = 1}^{n_1} \cdots \sum_{i_d = 1}^{n_d} {\pd t_{( i_1 , \dots , i_d )} \over \pd \ta } ( \ta ) (- \la _{( i_1 , \dots , i_d )} ) \non \\
&\quad \times \Big[ \prod_{d' = 0}^{D} \prod_{{i_1}' = 1}^{n_1} \cdots \prod_{{i_d}' = 1}^{n_d} \{ {\la _{( {i_1}' , \dots , {i_d}' )}}^{Z_{( {i_1}' , \dots , {i_d}' )} ( \ta )} e^{- t_{( {i_1}' , \dots , {i_d}' )} ( \ta ) \la _{( {i_1}' , \dots , {i_d}' )}} \} \Big] \pi ( \bla ) \Big) d\bla \Big] \non \\
&= E_{\bla } \Big[ \sum_{d = 0}^{D} \sum_{i_1 = 1}^{n_1} \cdots \sum_{i_d = 1}^{n_d} {\pd t_{( i_1 , \dots , i_d )} \over \pd \ta } ( \ta ) \la _{( i_1 , \dots , i_d )} ( \ta ) \non \\
&\quad \times \log \frac{ \int_{(0, \infty )^{n_1 \dotsm n_D}} \la _{( i_1 , \dots , i_d )} \big[ \prod_{d' = 0}^{D} \prod_{{i_1}' = 1}^{n_1} \cdots \prod_{{i_d}' = 1}^{n_d} \{ {\la _{( {i_1}' , \dots , {i_d}' )}}^{Z_{( {i_1}' , \dots , {i_d}' )} ( \ta )} e^{- t_{( {i_1}' , \dots , {i_d}' )} ( \ta ) \la _{( {i_1}' , \dots , {i_d}' )}} \} \big] \pi ( \bla ) d\bla }{ \int_{(0, \infty )^{n_1 \dotsm n_D}} \big[ \prod_{d' = 0}^{D} \prod_{{i_1}' = 1}^{n_1} \cdots \prod_{{i_d}' = 1}^{n_d} \{ {\la _{( {i_1}' , \dots , {i_d}' )}}^{Z_{( {i_1}' , \dots , {i_d}' )} ( \ta )} e^{- t_{( {i_1}' , \dots , {i_d}' )} ( \ta ) \la _{( {i_1}' , \dots , {i_d}' )}} \} \big] \pi ( \bla ) d\bla } \non \\
&\quad - \sum_{d = 0}^{D} \sum_{i_1 = 1}^{n_1} \cdots \sum_{i_d = 1}^{n_d} {\pd t_{( i_1 , \dots , i_d )} \over \pd \ta } ( \ta ) \lah _{( i_1 , \dots , i_d )} ( \Z ( \ta )) \Big] \non \\
&= E_{\bla } \Big[ \sum_{d = 0}^{D} \sum_{i_1 = 1}^{n_1} \cdots \sum_{i_d = 1}^{n_d} {\pd t_{( i_1 , \dots , i_d )} \over \pd \ta } ( \ta ) \{ \la _{( i_1 , \dots , i_d )} ( \ta ) \log \lah _{( i_1 , \dots , i_d )} ( \Z ( \ta )) - \lah _{( i_1 , \dots , i_d )} ( \Z ( \ta )) \} \Big] \text{,} \non 
\end{align}
where the second equality follows from Lemma \ref{lem:hudson}. 
Then it follows that 
\begin{align}
&E_{\bla } [ \log \{ p( \Y | \bla ) / \ph ^{( \pi )} ( \Y ; \X ) \} ] \non \\
&= \int_{0}^{1} \Big[ \sum_{d = 0}^{D} \sum_{i_1 = 1}^{n_1} \cdots \sum_{i_d = 1}^{n_d} {\pd t_{( i_1 , \dots , i_d )} \over \pd \ta } ( \ta ) ( \la _{( i_1 , \dots , i_d )} \log \la _{( i_1 , \dots , i_d )} - \la _{( i_1 , \dots , i_d )} ) \non \\
&\quad - E_{\bla } \Big[ \sum_{d = 0}^{D} \sum_{i_1 = 1}^{n_1} \cdots \sum_{i_d = 1}^{n_d} {\pd t_{( i_1 , \dots , i_d )} \over \pd \ta } ( \ta ) \{ \la _{( i_1 , \dots , i_d )} ( \ta ) \log \lah _{( i_1 , \dots , i_d )} ( \Z ( \ta )) - \lah _{( i_1 , \dots , i_d )} ( \Z ( \ta )) \} \Big] \Big] d\ta \non \\
&= \int_{0}^{1} E_{\bla } \Big[ \sum_{d = 0}^{D} \sum_{i_1 = 1}^{n_1} \cdots \sum_{i_d = 1}^{n_d} {\pd t_{( i_1 , \dots , i_d )} \over \pd \ta } ( \ta ) \la _{( i_1 , \dots , i_d )} \Big\{ {\lah _{( i_1 , \dots , i_d )} ( \Z ( \ta )) \over \la _{( i_1 , \dots , i_d )}} - 1 - \log {\lah _{( i_1 , \dots , i_d )} ( \Z ( \ta )) \over \la _{( i_1 , \dots , i_d )}} \Big\} \Big] d\ta \text{.} \non 
\end{align}
This completes the proof. 
\hfill$\Box$

\section*{Acknowledgments}
Research of the author was supported in part by JSPS KAKENHI Grant Number JP20J10427 from Japan Society for the Promotion of Science.

\end{document}